\documentclass[12pt, leqno]{amsart}

\usepackage{graphicx}
\usepackage{amssymb,amsmath,amsthm,amscd}
\usepackage{mathrsfs,mathdots}
\usepackage{enumerate,stmaryrd}
\usepackage{upgreek}
\usepackage[usenames,dvipsnames]{color}
\usepackage[colorlinks=true, pdfstartview=FitV,
 linkcolor=blue,citecolor=blue,urlcolor=blue]{hyperref}
\usepackage[all]{xy}

\usepackage{verbatim}
\usepackage{oldgerm}
\usepackage{yhmath}
\usepackage[normalem]{ulem}  
\usepackage{xspace}
\usepackage{xcolor,cancel}
\usepackage{tikz}
\tikzset{dynkdot/.style={circle,draw,scale=.38}}

\usepackage{marginnote}

\numberwithin{equation}{section} \allowdisplaybreaks[4]

\usepackage[colorinlistoftodos]{todonotes}
\newcommand{\nc}{\newcommand}

\newcommand{\ls}{(\hspace{-0.3ex}(}
\newcommand{\rs}{)\hspace{-0.3ex})}

\renewcommand{\le}{\leqslant}

\renewcommand{\ge}{\geqslant}

\nc{\op}{\operatorname}

\theoremstyle{plain}
\newtheorem{lemma}{Lemma}[section]
\newtheorem{proposition}[lemma]{Proposition}
\newtheorem{theorem}[lemma]{Theorem}
\newtheorem{maintheorem}{Main Theorem}

\newtheorem*{problem}{Problem}

\newtheorem{sublemma}{Sublemma}

\theoremstyle{definition}
\newtheorem{remark}[lemma]{Remark}
\newtheorem{conjecture}[lemma]{Conjecture}
\newtheorem{example}[lemma]{Example}

\newtheorem{definition}[lemma]{Definition}
\newtheorem{corollary}[lemma]{Corollary}
\newtheorem*{convention}{Convention}

\nc{\Prop}{\begin{proposition}} \nc{\enprop}{\end{proposition}}
\nc{\Th}{\begin{theorem}} \nc{\enth}{\end{theorem}}
\nc{\Lemma}{\begin{lemma}} \nc{\enlemma}{\end{lemma}}
\nc{\Cor}{\begin{corollary}} \nc{\encor}{\end{corollary}}
\nc{\Rem}{\begin{remark}} \nc{\enrem}{\end{remark}}
\nc{\Def}{\begin{definition}} \nc{\edf}{\end{definition}}
\nc{\Sub}{\begin{sublemma}} \nc{\ensub}{\end{sublemma}}
\nc{\Prob}{\begin{problem}} \nc{\enprob}{\end{problem}}

\nc{\shc}{\mathcal{C}}

\newcommand{\Q}{\mathbb{Q}}

\newcommand{\frakC}{\mathfrak{C}}

\newcommand{\tb}{{\widetilde{b}}}
\newcommand{\ta}{{\widetilde{a}}}

\newcommand{\C}{\mathbb{C}}

\newcommand{\Seq}{\Sigma}

\newcommand{\dT}{\mathrm{T}}

\newcommand{\A}{\mathcal{A}}
\newcommand{\im}{\imath}
\newcommand{\jm}{\jmath}

\newcommand{\Ca}{\mathscr{C}}
\nc{\F}{\mathcal{F}}
\newcommand{\fD}{\mathsf{D}}
\newcommand{\uxi}{\upvartheta}
\newcommand{\N}{\ell}
\newcommand{\D}{\mathscr{D}}
\nc{\dual}{\D}
\newcommand{\Dd}{\mathcal{D}}

\nc{\HOM}{\on{H\textsc{om}}}
\newcommand{\M}{\mathcal{M}}
\newcommand{\W}{\mathsf{W}}
\newcommand{\fd}{\mathsf{d}}

\newcommand{\DC}{{\mathsf{D}}}

\newcommand{\LL}{\mathcal{L}}
\newcommand{\RR}{\mathcal{R}}

\newcommand{\bS}[1]{\ms{1mu}{\mathsf{S}}_{#1} }
\newcommand{\cusp}{\ms{1mu}{\mathsf{S}}}
\newcommand{\cuspS}[1]{\cusp_{#1} }

\newcommand{\B}{\mathbf{B}}

\newcommand{\sfJ}{\mathsf{J}}

\newcommand{\wvee}{\widetilde{\vee}}
\newcommand{\Dynkin}{\triangle}

\newcommand{\qm}{{q}_{\mspace{1mu}\raisebox{-.5ex}{${\scriptstyle{\mathrm{sh}}}$}}}

\newcommand{\Z}{\mathbb{\ms{1mu}Z\ms{1mu}}}
\newcommand{\seteq}{\mathbin{:=}}
\newcommand{\hd}{{\operatorname{hd}}}
\newcommand{\soc}{{\operatorname{soc}}}
\nc{\ov}[1]{\overline{#1}} \nc{\wit}[1]{\widetilde{#1}}

\newcommand{\hI}{\widehat{I}}

\newcommand{\fdi}{\ms{2mu}\fd_{\im}}
\newcommand{\fdj}{\ms{2mu}\fd_{\jm}}
\newcommand{\fdii}{\fd_{\im}}
\newcommand{\fdjj}{\fd_{\jm}}
\newcommand{\oi}{{{\im}}}

\nc{\dis}[1]{\vert\ms{1.8mu}{#1}\ms{2mu}\vert}
\nc{\disp}[1]{\dis{#1}_\phi}

\nc{\Wlmj}[3]{\W_{#2,#3}^{(#1)}} \nc{\Mkl}[2]{\M_\ttww(#1,#2)}
\nc{\rmat}[1]{{\mathbf
r}_{\mspace{-2mu}\raisebox{-.5ex}{${\scriptstyle{#1}}$}}}
\nc{\po}[1]{p_{\mspace{-2mu}\raisebox{-.5ex}{${\scriptstyle{#1}}$}}}
\nc{\ddfrac}[2]{  \dfrac{#2}{#1} }
\newcommand{\on}{\operatorname}
\nc{\de}{\on{\textfrak{d}}} \nc{\tL}{\widetilde{\Lambda}}
\nc{\tl}{\widetilde{\lambda}} \nc{\mqs}{(-q^2)}
\nc{\Cquiver}{\upsigma}

\nc{\mut}[1]{{\mu}_{\mspace{-2mu}\raisebox{-.5ex}{${\scriptstyle{#1}}$}}}

\newcommand{\g}{{\mathfrak{g}}}
\newcommand{\fing}{{\mathsf{g}}}
\newcommand{\finn}{{\mathsf{n}}}
\newcommand{\An}{A_q(\finn)}
\newcommand{\Anw}{A_q(\finn(w))}
\newcommand{\up}{\mathrm{up}}
\newcommand{\h}{{\mathfrak{h}}}
\newcommand{\n}{{\mathfrak{n}}}
\newcommand{\isoto}[1][]{\mathop{\xrightarrow%
[{\raisebox{.3ex}[0ex][.3ex]{$\scriptstyle{#1}$}}]%
{{\raisebox{-.6ex}[0ex][-.6ex]{$\mspace{2mu}\sim\mspace{2mu}$}}}}}
\newcommand{\longisoto}[1][]{\mathop{\xrightarrow%
[{\raisebox{.3ex}[0ex][.3ex]{$\scriptstyle{\hs{1ex}#1\hs{1ex}}$}}]%
{{\raisebox{-.6ex}[0ex][-.6ex]{$\hs{1ex}\sim\hs{1ex}$}}}}}

\newcommand{\id}{\on{id}}

\newcommand{\soplus}{\mathop{\mbox{\normalsize$\bigoplus$}}\limits}

\newcommand{\sotimes}{\mathop{\mbox{\normalsize$\bigotimes$}}\limits}

\newcommand{\ww}{ \textbf{\textit{w}}}
\newcommand{\ttww}{{\widetilde{\ww}} }

\nc{\nconv}{\mathop{\mbox{\large $\odot$}}}
\nc{\nnconv}{\mathop{\mbox{\large $\star$}}}

\nc{\lb}{\llbracket} \nc{\rb}{\rrbracket}

\newcommand{\ko}{{{\mathbf{k}}}}
\nc{\la}{\lambda} \nc{\La}{\Lambda}

\nc{\finite}{{\mathrm{fin}}} \nc{\gf}{{\g_\finite}}
\nc{\tLa}{\widetilde{\Lambda}} \nc{\ve}{\varepsilon}
\nc{\ep}{\epsilon} \nc{\vp}{\varphi} \nc{\lan}{\langle}
\nc{\ran}{\rangle} \nc{\Uqg}{U_q(\g)} \nc{\Aqg}{A_q(\g)}
\nc{\Aqn}{A_q(\n)} \nc{\ual}{\upalpha\ms{1mu}}
\nc{\sPi}{\mathsf{\Pi}}
\nc{\Pif}{\Pi_\finite} 
\nc{\sLa}{\mathsf{\Lambda}\ms{1mu}} \nc{\al}{\alpha} \nc{\be}{\beta}
\nc{\ga}{\gamma} \nc{\wt}{\operatorname{wt}}
\nc{\ch}{\operatorname{ch}}

\nc{\norm}{{\mathrm{norm}}} \nc{\aff}{{\mathrm{aff}}}
\nc{\Maf}{M_\aff} \nc{\ev}{{\mathrm{even}}} \nc{\od}{{\mathrm{odd}}}
\nc{\Sev}{\Seq^{\ev}} \nc{\Sod}{\Seq^{\od}} \nc{\Spl}{\Seq^{+}}
\nc{\Smi}{\Seq^{-}} \nc{\low}{{\mathrm{low}}}
\nc{\upper}{{\mathrm{up}}} \nc{\one}{{\bf{1}}}
\nc{\To}[1][{\hspace{2ex}}]{\xrightarrow{\,#1\,}}
\nc{\te}{\tilde{e}} \nc{\tw}{{\underline{w}}}
\nc{\hw}{{\mathfrak{s}}}

\nc{\hhw}{\widehat{\tw}_0}
\nc{\tww}{\ww} \nc{\tuu}{{\mathsf{u}}} \nc{\tel}{\tilde{e}^\low}
\nc{\teu}{\tilde{e}^\upper} \nc{\tf}{\tilde{f}}
\nc{\tfl}{\tilde{f}^\low} \nc{\tfu}{\tilde{f}^\upper}
\nc{\tE}{\widetilde{E}} \nc{\tF}{\widetilde{F}}
\nc{\tFF}{\widetilde{\F}} \nc{\tB}{\widetilde{B}}
\nc{\tz}{\tilde{z}}
\nc{\tQ}{\hspace{-.2ex}\textbf{\textit{Q}}}
\nc{\Ft}{\F^\dT}
\nc{\Seed}{\mathcal{S}}

\newenvironment{magenta}{\relax\color{magenta}}{\hspace*{.5ex}\relax}

\newcommand{\bem}{\begin{magenta}}
\newcommand{\eem}{\end{magenta}}

\nc{\ble}[1]{\underline{#1}}

\nc{\cor}{\mathbf{k}} \nc{\tens}{\mathop\otimes}
\nc{\gmod}{\mbox{-$\mathrm{gmod}$}}
\nc{\gMod}{\mbox{-$\mathrm{gMod}$}} \nc{\Md}{\mbox{-$\mathrm{Mod}$}}
\nc{\md}{\mbox{-$\mathrm{mod}$}} \nc{\uqm}{\mathscr{C}_\g}

\nc{\proj}{\mbox{-$\mathrm{proj}$}}
\nc{\gproj}{\mbox{-$\mathrm{gproj}$}}
\nc{\smod}{\mbox{-$\mathrm{mod}$}}
\nc{\nmod}{\mbox{-$\mathrm{nilmod}$}}

\nc{\seed}{\mathscr{S}}

\newcommand{\cmA}{\mathsf{A}}
\newcommand{\cmf}{\mathsf{C}} 

\newcommand{\rmQ}{\mathrm{Q}}

\nc{\Rnorm}{\mathrm{R}^{\mathrm{norm}}}
\nc{\Runiv}{\mathrm{R}^{\ms{1mu}\mathrm{univ}}}
\nc{\Rren}{\mathrm{R}^{\ms{1mu}\mathrm{ren}}} 
\nc{\col}{\colon} \nc{\epiTo}[1][]{\xymatrix{\ar@{->>}[r]^-{{#1}}&}}
\nc{\epito}{\twoheadrightarrow}
\nc{\monoTo}[1][]{\xymatrix{\ar@{>->}[r]^-{{#1}}&}}
\nc{\monogets}[1][]{\xymatrix{&\ar@{_{(}->}[l]^-{{#1}}}}
\nc{\sym}{\mathfrak{S}} \nc{\rl}{\mathsf{Q}} \nc{\prl}{\rl_+}
\nc{\crl}{\mathsf{Q}^\vee} \nc{\pcrl}{\crl_+} \nc{\Qq}{{\Q(q)}}
\nc{\wl}{\mathsf{P}}   
\nc{\wlf}{\mathsf{P}_\finite}   
\nc{\Oint}{\mathcal{O}_{{\mathrm{int}}}}
\newcommand{\scbul}{{\,\raise1pt\hbox{$\scriptscriptstyle\bullet$}\,}}

\nc{\conv}{\mathop{\mathbin{\mbox{\large $\circ$}}}}
\newcommand{\hconv}{\mathbin{\scalebox{.9}{$\nabla$}}}
\newcommand{\sconv}{\mathbin{\scalebox{.9}{$\Delta$}}}
\newcommand{\hconvs}{\mathbin{\scalebox{.6}{$\nabla$}}}

\nc{\scrA}{\mathscr{A}}

\nc{\pv}{  \to\updownarrow\gets }
\nc{\nv}{  \longleftrightarrow {\raise -1pt\hbox{$\hspace{-2ex}\begin{matrix}\downarrow \\[-1ex] \uparrow\end{matrix}$}} }

\newcommand{\Hom}{\operatorname{Hom}}

\renewcommand{\Im}{\op{Im}}

\newcommand{\ex}{{\mathrm{ex}}}
\nc{\K}{\mathsf{K}} \nc{\Kex}{{\K}^{\mathrm{ex}}}
\nc{\J}{\mathsf{J}}
\nc{\Uex}{\Uppsi_{\mathrm{ex}}}
\nc{\Kfr}{{\K}^{\mathrm{f\mspace{.01mu}r}}}
\nc{\cl}{{\ms{1mu}\mathrm{cl}\ms{1mu}}} \nc{\ben}{\begin{enumerate}}
\nc{\ee}{\end{enumerate}} \nc{\bnum}{\begin{enumerate}[{\rm(i)}]}
\nc{\bna}{\begin{enumerate}[{\rm(a)}]} \nc{\bc}{\begin{cases}}
\nc{\ec}{\end{cases}}

\newenvironment{myequation}
{\relax\setlength{\arraycolsep}{1pt}\begin{eqnarray}}
{\end{eqnarray}}
\newenvironment{myequationn}
{\relax\setlength{\arraycolsep}{1pt}\begin{eqnarray*}}
{\end{eqnarray*}}

\nc{\eq}{\begin{myequation}} \nc{\eneq}{\end{myequation}}
\nc{\eqn}{\begin{myequationn}} \nc{\eneqn}{\end{myequationn}}
\nc{\hs}{\hspace*} \nc{\hl}{\hspace{-.5ex}}

\newenvironment{myarray}[1]{\relax\setlength{\arraycolsep}{.5pt}
\renewcommand{\arraystretch}{1.3}
\begin{array}{#1}}{\end{array}\relax}

\newcommand{\ba}{\begin{myarray}}
\newcommand{\ea}{\end{myarray}}

\nc{\noi}{\noindent} \nc{\ang}[1]{\langle{#1}\rangle}
\nc{\fr}{{\mathrm{fr}}} \nc{\qt}[1]{\quad\text{#1}\quad}
\nc{\ol}{\overline} \nc{\true}{\delta} \nc{\ms}{\mspace}
\renewcommand{\mod}{\ms{3mu}\mathbin{\mathrm{mod}}\ms{1mu}}
\nc{\vs}{\vspace*} \nc{\bl}{\bigl(} \nc{\br}{\bigr)}
\nc{\bep}{\ol{\ep}} \nc{\bal}{\,\ol{\al}}
\nc{\qtq}[1][{and}]{\quad\text{#1}\quad}
\nc{\set}[2]{\left\{{#1}\mid{#2}\right\}} \nc{\rmo}{{\rm(}}
\nc{\rmf}{{\rm)}\xspace} \nc{\Proof}{\begin{proof}}
\nc{\QED}{\end{proof}}
\nc{\monoto}[1][]{\xymatrix@C=2ex{\ar@{>->}[r]^-{{#1}}&}\ms{-8mu}}
\nc{\etens}{\boxtimes} \nc{\height}[1]{\vert{#1}\vert}
\nc{\Lrev}{L^{\bek{\rev}\ek}}
\nc{\rev}{\mathrm{rev}} \nc{\fw}{\Lambda} \nc{\uqpg}{U_q'(\g)}
\nc{\tp}{\ms{1.5mu}{\widetilde{p}}\ms{2mu}}
\nc{\Deg}{\mathrm{\ms{1mu}Deg\ms{1mu}}} \nc{\Bg}{\mathcal{G}}

\nc{\wb}[1]{\mbox{$\rule[-1.1ex]{0ex}{2ex}#1$}}

\nc{\bwr}{\mbox{\large$\wr$}} \nc{\vphi}{\varphi}
\nc{\G}{\mathcal{G}} \nc{\tD}{\widetilde{\mathrm{De}}\mathrm{g}}
\nc{\Li}{\La^\infty} \nc{\Di}{\Deg^\infty}
\nc{\zero}{\ms{2mu}\mathrm{zero}\ms{2mu}} \nc{\cwl}{\wl^\vee}
\nc{\rc}{renormalizing coefficient\xspace}
\nc{\cz}{{\cor[z^{\pm1}]}}
\nc{\ake}[1][2ex]{\rule[-.5ex]{0ex}{#1}}
\nc{\akete}[1][0ex]{\rule[{#1}]{0ex}{1ex}}
\nc{\akew}[1][2ex]{\rule[-1ex]{#1}{0ex}} \nc{\rd}{{}^*\ms{-3mu}}
\nc{\st}[1]{\{{#1}\}}
\nc{\seq}[1]{ \boldsymbol{(} {#1} \boldsymbol{)}   }
\nc{\corh}{\widehat{\cor}}
\nc{\czt}{\cz^\times} \nc{\eps}{\varepsilon} \nc{\rr}{rationally
renormalizable\xspace} \nc{\QHA}{\mathrm{QHA}} \nc{\Ker}{\on{Ker}}
\nc{\usq}{U'_q(\g)} \nc{\Wf}{\W_\finite} \nc{\If}{I_\finite}
\nc{\ord}{\mathrm{ord}}
\nc{\Qd}{\mathscr{Q}} \nc{\e}{\mathrm{e}} \nc{\snoi}{\smallskip\noi}
\nc{\mnoi}{\medskip\noi}
\nc{\ci}{\mathfrak{c}}
\newcommand{\tc}{{\widetilde{\ci}}}
\nc{\Iff}{I_\fing} \nc{\CM}[1][s]{M[{#1},0\}}
\nc{\CMp}[1][s]{M[{#1},0\}'} \nc{\cM}[1][s]{\mathsf{M}_{#1}}
\nc{\fM}{\mathsf{M}} \nc{\cMp}[1][s]{\mathsf{M}'_{#1}}
\nc{\Vi}{\mathrm{Vi}} \nc{\Vo}{\mathrm{Vo}}
\nc{\GLS}{\mathrm{Q}_{\mathrm{GLS}}}
\nc{\HL}{\mathrm{Q}_{\mathrm{HL}}} \nc{\vpi}{\varpi}
\nc{\bpi}{\ol{\pi}}
\nc{\dC}{\ms{1mu}\mathsf{h}^\vee} 
\nc{\ca}{completely $\Uplambda$-admissible\xspace}
\nc{\KR}{Kirillov-Reshetikhin\xspace}
\nc{\qA}{quasi-admissible\xspace}
\nc{\Ad}{affine determinantial\xspace}
\nc{\AD}{Affine determinantial\xspace}
\nc{\SW}{Schur-Weyl\xspace}
\nc{\xiz}{\xi^{\mathrm{zz}}} \nc{\Kin}{\K_{\mathrm{in}}}
\nc{\Kb}{\K_{\mathrm{b}}} \nc{\Kout}{\K_{\mathrm{out}}} \nc{\rE}{
\mathsf{E} } \nc{\rW}{ \mathcal{W} } \nc{\catCO}{\Ca^0_\g}
\nc{\sig}{\upsigma(\g)} \nc{\sigZ}{\upsigma_0(\g)}
\nc{\cato}{\catCO} \nc{\cat}{\Ca} \nc{\catg}{\Ca_\g}
\nc{\Set}{{\mathrm{Set}}} \nc{\BHL}{\tB_{\mathrm{HL}}}
\nc{\R}{\mathbb{R}} \nc{\BGLS}{\tB_{\mathrm{GLS}}}
\nc{\hIg}{{\hI_\g}} \nc{\KK}{;\K,\Kex}
\nc{\KKC}{;\K(\frakC),\Kex(\frakC)} \nc{\Fd}{\F_\Dd}
\nc{\KKCp}{;\K(\frakC'),\Kex(\frakC')}
\nc{\hF}{\widehat{\F}} \nc{\Rc}{R_\cmf} \nc{\pbw}{PBW-pair\xspace}
\nc{\pbws}{PBW-pairs\xspace} \nc{\Refl}{\mathscr{S}}
\nc{\Reflinv}{{\Refl}^{-1}} \nc{\Rt}{\mathsf{L}}
\nc{\qa}{$\Uplambda$-quasi-admissible\xspace}
\nc{\Lad}{$\Uplambda$-admissible\xspace}
\nc{\htens}{\hconv}
\nc{\Isf}{I_\fing}
\nc{\scb}{\scalebox}
\nc{\vrho}{\varrho}
\nc{\dsum}{\displaystyle\sum}
\nc{\dtens}{\mathop{\mbox{\scb{1.1}{\akete[-.6ex]\normalsize$\bigotimes$}}}\limits}
\nc{\cq}{\Check{q}}
\nc{\ud}[1]{\underline{#1}}
\nc{\indlim}[1][]{\mathop{\varinjlim}\limits_{#1}}
\nc{\me}{monoidally mutation equivalent\xspace}
\nc{\eqv}[1][\ell]{\equiv_{#1}}

\newenvironment{magem}{\relax\color{magenta}}{\relax}

\newcommand{\bema}{\begin{magem}}
\newcommand{\ema}{\end{magem}}

\newlength{\mylength}
\title[Monoidal categorification and  quantum affine algebras II]{Monoidal categorification and \\  quantum affine algebras II}

\author[M. Kashiwara]{Masaki Kashiwara}
\thanks{The research of M.\ Kashiwara
was supported by Grant-in-Aid for Scientific Research (B) 20H01795,
Japan Society for the Promotion of Science.}
\address[M. Kashiwara]{%
Kyoto University Institute for Advanced Study, Research Institute
for Mathematical Sciences, Kyoto University, Kyoto 606-8502, Japan
\& Korea Institute for Advanced Study, Seoul 02455, Korea }
\email[M. Kashiwara]{masaki@kurims.kyoto-u.ac.jp}

\author[M. Kim]{Myungho Kim}
\address[M. Kim]{Department of Mathematics, Kyung Hee University, Seoul 02447, Korea}
\email[M. Kim]{mkim@khu.ac.kr}
\thanks{The research of M.\ Kim was supported by the National Research Foundation of
Korea (NRF) Grant funded by the Korea government(MSIP)
(NRF-2017R1C1B2007824  and NRF-2020R1A5A1016126).}

\author[S.-j. Oh]{Se-jin Oh}
\thanks{ The research of S.-j.\ Oh was supported by the Ministry of Education of the Republic of Korea and the National Research Foundation of Korea (NRF-2019R1A2C4069647).}
\address[S.-j. Oh]{Department of Mathematics, Ewha Womans University, Seoul 03760, Korea}
\email[S.-j. Oh]{sejin092@gmail.com}

\author[E. Park]{Euiyong Park}
\thanks{The research of E.\ Park was supported by the National Research Foundation of Korea (NRF) Grant funded by the Korea Government(MSIP)(NRF-2020R1F1A1A01065992 and NRF-2020R1A5A1016126).}
\address[E. Park]{Department of Mathematics, University of Seoul, Seoul 02504, Korea}
\email[E. Park]{epark@uos.ac.kr}

\keywords{Quantum affine algebra, Monoidal categorification,
R-matrices, Cluster algebra}

\subjclass[2010]{17B37, 13F60, 18D10}

\date{March 18, 2021}

\begin{document}

\begin{abstract}
We introduce a new family of real simple modules over the quantum affine algebras, called the \Ad modules, which contains the
\KR (KR)-modules as a special subfamily, and then
prove T-systems among them
which generalize the T-systems among KR-modules and unipotent quantum minors in the quantum unipotent coordinate algebras simultaneously.
We develop new combinatorial tools:
admissible chains of $i$-boxes which produce commuting families of
\Ad modules, and box moves which describe the T-system in a combinatorial way.
Using these results,
we prove that various module categories over the quantum affine algebras provide monoidal categorifications of cluster algebras.
As special cases,
Hernandez-Leclerc categories $\Ca_\g^0$ and $\Ca_\g^-$
provide monoidal categorifications of the cluster algebras
 for an arbitrary quantum affine algebra.
\end{abstract}

\maketitle \tableofcontents

\section{Introduction}
The notion of \emph{monoidal categorification} is
proposed by Hernandez-Leclerc in \cite{HL10}. This notion offers a framework for proving Laurent positivity
and linear independence for  a cluster algebra $\mathscr{A}$,
when $\mathscr{A}$ is isomorphic to the Grothendieck ring of
  a  monoidal   category $\Ca$.
Conversely,
it provides quite interesting information on the monoidal category $\Ca$:
once we show that $\Ca$ is
a monoidal categorification of a cluster algebra,  we acquire a family $\mathscr{F}$ of \emph{real prime simple objects} in $\Ca$ (identified with
the \emph{cluster variables} in  $\mathscr{A}$)
whose certain groupings in $\mathscr{F}$ (identified with  the \emph{clusters} in $\mathscr{A}$) consist of \emph{mutually commuting} real prime simple objects.

\smallskip

The purpose of this paper is to prove that various module categories over quantum affine algebras
provide monoidal categorifications of  cluster algebras.

\smallskip

A cluster algebra $\mathscr{A}$,  introduced by Fomin and Zelevinsky in \cite{FZ02}, is a commutative  $\Z$-subalgebra of
$\Z[X_k^{\pm1}\mid{k \in \K}]$  with a distinguished set of generators, called the cluster variables,
which is grouped into overlapping subsets, called the clusters. The clusters are defined inductively by a procedure called \emph{mutation} from the initial cluster $\{ X_k\}_{k \in \K}$, which is performed via an
exchange matrix $\tB$.
The notion of cluster algebras  is  extended to  a
quantum version, \emph{quantum cluster algebras} $\mathscr{A}_q$ in $\Z[q^{\pm 1/2}][X_k^{\pm1}]_{k \in \K}$ by Berenstein and Zelevinsky in \cite{BZ05}, which are not commutative any more, but their cluster variables
in each cluster are $q$-commutative.
The  $q$-commutativity is controlled by a $\Z$-valued $\K \times \K$-matrix $L$.
From their introductions,  numerous connections and applications
have been discovered in  various fields of mathematics (see \cite{FZ02,BZ05,L10,RJM14,W14,GP17} and references therein).

\smallskip

For each Kac-Moody algebra $\g$ of affine type, let $U_q'(\g)$ be the corresponding quantum affine algebra and  let $\Ca_\g$
be the category of finite-dimensional integrable modules over $U_q'(\g)$. Since the category $\Ca_\g$ has interesting properties including
monoidality and rigidity, it has been intensively studied since 1990's, in various aspects of point of view. To name a few,
the complete classification of simple modules in $\Ca_\g$
is obtained in terms of Drinfel'd polynomials (\cite{CP91,CP94,CP95} for the untwisted cases and \cite{CP98} for the twisted cases)
and it is proved in \cite{AK97, Kas02, VV02} that every simple module in $\Ca_\g$ can be obtained as the head of a tensor product of \emph{fundamental modules}.
Also,
it is proved in \cite{Nak,Nak04} (for simply-laced untwisted affine types)
and \cite{Her06,Her10} (for general types) that the {\it $q$-characters} of
\emph{\KR\ {\rm(KR)}\ modules}, a special class of modules in $\Ca_\g$,
are solutions of the \emph{T-system}
which is closely related to discrete dynamical systems
arising from the thermodynamic Bethe-ansatz, Y-system (see \cite{KNS,IIKNS}).

The first result on monoidal categorifications
was established in \cite{HL10,Nak11}
for relatively small  monoidal subcategories
$\Ca_\g^1$ of $\Ca_\g$ for untwisted simply-laced affine type $\g$. One of the main ideas of those papers is interpreting  the T-system among KR-modules
as exchange relations of a cluster algebra by mutations.
(See \cite{Qin17} for the monoidal categorifications of related categories $\Ca_\g^N$ $(N\in\Z_{\ge1})$
and see \cite{IIKKN1,IIKKN2} for the relation between cluster algebras and T-systems.)

One of other successful instances on the monoidal categorification is given in \cite{KKKO18} by using the monoidal subcategories $\mathcal{C}_w$ in $R\gmod$
which \emph{categorify} the quantum unipotent coordinate algebras.
Here $R$ denotes a $\Z$-graded algebra, called  \emph{quiver Hecke algebra}, introduced  independently
by Khovanov-Lauda \cite{KL09} and Rouquier \cite{R08,R11},
and $R\gmod$ denotes the category of finite-dimensional graded modules over $R$.
We shall explain the result in a more precise way.
In \cite{GLS11, GLS13S},  Gei\ss, Leclerc and Schr{\"o}er showed that the quantum unipotent coordinate algebra $A_q(\mathfrak{n}(w))$, associated with a
symmetric quantum group $U_q(\mathsf{g})$ and its Weyl
group element $w$, has a skew-symmetric quantum cluster algebra structure.
To see this, they (i) used  a system of quantum determinantial identities
among \emph{unipotent quantum minors} $D_\tw[a,b]$ in $A_q(\mathfrak{n}(w))$, called also T-system, (ii) constructed an initial quiver $\GLS(\tw)$ arising from a choice of a \emph{reduced expression} $\tw$ of $w$
and (iii) employed the representation theory of preprojective algebras related to $A_q(\mathfrak{n}(w))$.
Here, the T-system also plays the role of exchange relation of the quantum cluster algebra by mutations.
In \cite{KKKO18},
it is proved that $\mathcal{C}_w$ provides a monoidal categorification of the quantum cluster algebra $A_q(\mathfrak{n}(w))$ by showing
\begin{enumerate}[{\rm (i)}]
\item \label{it: MC i} $\Z[q^{\pm 1/2}] \otimes_{\Z[q^{\pm1}]} K(\shc_w) \simeq A_{q^{1/2}}(\mathfrak{n}(w)) \seteq
\Z[q^{\pm 1/2}] \otimes_{\Z[q^{\pm1}]} A_q(\mathfrak{n}(w))$,
\item there exists a \emph{quantum monoidal seed} $\seed=(\{ M_i \}_{i \in \K},L,\tB, D)$ in $\shc_w$ consisting of
\bna
\item
a commuting family $\{ M_i\}_{i\in \K}$ of real simple modules in $\shc_w$,
called \emph{determinantial modules}, corresponding to the $q$-commuting family of quantum \emph{flag} minors $D_\tw[0,b]$ via the isomorphism in~\eqref{it: MC i},
\item the matrix $L=(-\Lambda(M_i,M_j))_{i,j\in \K}$ where $\Lambda(M_i,M_j)$ denotes
the homogeneous degree of the \emph{$R$-matrix}
between $M_i$ and $M_j$,
\item
 the  incident matrix $\tB$ of  the quiver $\GLS(\tw)$,
 and
\item the set $D$ of \emph{weights} of $M_i$'s in the root lattice of $\mathsf{g}$
\ee
such that $[\seed]\seteq ( \{ q^{m_i}[M_i] \}_{i \in \K}, L,\tB)$ is a quantum seed of $A_{q^{1/2}}(\mathfrak{n}(w))$ for some $m_i \in \tfrac{1}{2}\Z$,
\item $\seed$ admits successive mutations in all directions in $\Kex$, which is \emph{revealed} as a consequence of the \emph{$\Uplambda$-admissibility} of $\seed$ and the isomorphism $A_{q^{1/2}}(\mathfrak{n}(w)) \simeq \mathscr{A}([\seed])$  \cite[Theorem 7.1.3]{KKKO18}.
\end{enumerate}
Here the $\Uplambda$-admissibility of $\seed$ means that, for every $k \in \Kex$,
\begin{eqnarray}
&&\ba{l}
\text{there exists a simple module $M_k'$ satisfying the following properties:}\\
\hs{4ex}\parbox{74ex}{
\bna
\item \label{it: quasi-ad} $M_k'$ corresponds to the mutated cluster variable $X_k'$ in the direction $k$,
\item \label{it: Lam-quasi-ad} $1=\de(M_k,M_k')\seteq \bl \Lambda(M_k,M_k')+ \Lambda(M_k',M_k) \br/2$,
\item \label{it: Lam-ad} $M'_k$ is real and
commutes with $M_i$ for all $i \in \K \setminus \st{ k }$.
\ee
}\ea\label{eq: admissibility}
\end{eqnarray}
The $\Z$-grading structure of $R$,
the quantum cluster algebra structure of $A_q(\mathfrak{n}(w))$ and the primeness of
cluster variables (\cite{GLS13}) play important role
in the proof of the monoidal categorification by $\mathcal{C}_w$ in \cite{KKKO18}.

\smallskip

In \cite{KKK18A}, Kang-Kashiwara-Kim constructed a functor $\mathcal{F}$, called \emph{quantum affine \SW duality functor}, from $R\gmod$
to $\Ca_\g$ by observing the singularities of the $U_q'(\g)$-module homomorphism  $\Rnorm_{M,N_z}$, called the \emph{normalized R-matrix}, for modules $M,N$ in $\Ca_\g$.
The quantum affine \SW duality functor thus makes a bridge
between the representation theory of quiver Hecke algebras
and the one of quantum affine algebras.

On the other hand, a monoidal full subcategory $\Ca_Q$ of $\Ca_\g$ is introduced in \cite{HL11} for simply-laced untwisted affine $\g$,
which is defined by using the combinatorics of the Auslander-Reiten quiver of
a Dynkin quiver $Q$. Then the notion of $\Ca_Q$ is extended to
all affine types by generalizing
Dynkin quivers $Q$ to $\rmQ$-data $\Qd$  in \cite{KKKO16D,OhSuh19,OT19,FO20} (see \S~\ref{subsec: Q-datum}).
In \cite{KKK15B,KKKO16D,KO18,OT19},
the \SW duality functor $\mathcal{F}_\Qd\col R_\gf\gmod\to\Ca_\Qd $ is constructed for
a  $\rmQ$-datum $\Qd$, and $\mathcal{F}_\Qd$ sends simple modules in $R_\gf\gmod$ to simple modules in the monoidal full subcategory $\Ca_\Qd \subset \Ca_\g$ bijectively. Here, we associate
a simply-laced finite-dimensional simple Lie algebra
 $\gf$ to each affine Lie algebra $\g$
(see \eqref{Table: root system} and \eqref{Table: root system for twisted}),
and $R_\gf$ is the quiver Hecke algebra associated with $\gf$.

Thus one can conclude that $\Ca_\Qd$ provides a monoidal categorification of the coordinate ring $\C[N]$
of the maximal unipotent group $N$ associated with $\gf$ (see also \cite{HL11,HO20,FHOO}).

\smallskip

In $\Ca_\g$, there are other interesting and important monoidal full subcategories (see \S\,\ref{subsec:subcategory} for more details).
\ben[{\rm(A)}]
\item \label{it:catA}
The subcategory $\Ca_\g^0$ is defined in \cite{HL10} for a simply-laced untwisted affine $\g$ satisfying the following property:  for any simple module $V$ in $\Ca_\g$,
it decomposes as a tensor product of parameter shifts of simple modules in $\Ca_\g^0$. Thus sometimes  $\Ca_\g^0$ is referred to as the \emph{skeleton} subcategory.
\item \label{it:catB} The subcategory $\Ca_\g^-$ is defined in \cite{HL16} for
an untwisted affine $\g$ which contains all simple modules in $\Ca_\g^0$ up to parameter shifts.
It is proved in \cite{HL16} that the Grothendieck ring $K(\Ca_\g^-)$ has a cluster algebra structure with an initial cluster
consisting of KR-modules,
 by using  T-systems  among KR-modules.
Note that the  definition of the subcategories $\Ca_\g^0$ and $\Ca_\g^-$ are also extended to all affine types.

\item \label{it:catC} The subcategory $\Ca_\Qd$
has a remarkable property as a subcategory of $\Ca_\g^0$:  for any fundamental module $V$ in $\Ca_\g^0$,
there exists a unique fundamental module $U$ in $\Ca_\Qd$ and $k \in \Z$ such that $V \simeq \D^k(U)$. Here $\D^k(U)$ denotes the $k$-th repeated dual of $U$.
Thus  $\Ca_\Qd$ is understood as a \emph{heart} subcategory of $\Ca_\g^0$.
\item \label{it:catD} For $N \in \Z_{\ge 1}$, the subcategory $\Ca_\g^N$ is defined in \cite{HL10} (denoted by $\Ca_{N-1}$ in \cite{HL10})
for a simply-laced untwisted affine $\g$, which is generated by
$(|\If| \times N)$-many fundamental modules.
Here $\If$ is the index set of
simple roots of $\gf$.
\ee
The subcategories in \eqref{it:catA}--\eqref{it:catD} are also referred to as
\emph{Hernandez-Leclerc categories}.

\medskip
On the other hand,
the authors of this paper have recently developed interesting results
on  $\Ca_\g$ which are briefly summarized as follows.

\snoi
(I) In \cite{KKOP19C,KKOP20A}, the $\Z$-valued invariants $\Lambda(M,N)$, $\de(M,N)$ and $\Lambda^\infty(M,N)$
for a pair of modules $M$ and $N$ in $\Ca_\g$ are introduced by analyzing
the $R$-matrices associated with $M \otimes N_z$. By using them,
we associate a simply-laced finite root system of type $\gf$
to each $\Ca_\g$.
They can be understood as quantum affine analogues of $\Z$-graded structure and weights of modules in
$R\gmod$, respectively.
Contrary to  $R\gmod$, the category  $\Ca_\g$ is a rigid monoidal category and
these invariants enjoy interesting properties
that result from rigidity.  (See \S\,\ref{subsec:invariants}.)
Also in \cite{KKOP19C}, a \emph{criterion} for a subcategory $\shc$ in $\Ca_\g$ to provide a monoidal categorification of a cluster algebra is established (see Theorem~\ref{th:main KKOP19C}).

\snoi
(II) In \cite{KKOP20C}, the authors generalize the notion of $\rmQ$-datum $\Qd$ one step further to the notion of \emph{\rmo complete\rmf duality datum $\Dd$}
(see \cite[Definition 4.7]{KKOP20C} and Definition~\ref{def:complete dd}).
It induces a \SW duality functor
$\F_\Dd\col R_{\fing}\gmod\to\Ca_\g^0$ which sends simple modules in $R_{\fing}\gmod$ to simple modules in $\Ca_\g^0$ injectively,
where $\fing$ is the simply-laced finite-dimensional simple Lie algebra determined by $\Dd$ (see \S~\ref{subsec: QAWS}).
Also it is proved that the category $\Ca_\Dd$,
the image of $\F_\Dd$, also enjoys the similar properties to those of
$\Ca_\Qd$ in the following sense: for each complete duality datum $\Dd$ and a reduced expression $\tw_0$ of the longest element $w_0 \in W_{\fing}$,
$\Ca_\g^0$ is generated by the images $\{ \bS{k} \}_{1 \le k \le \ell(w_0)}$
of $\{ \DC_{\tw_0}[k] \}_{1 \le k \le \ell(w_0)}$ in $R_{\fing}\gmod$ by $\F_\Dd$ (Definition~\ref{def: root module}),
and their repeated duals $\{ \bS{k} \}_{k \in \Z}$.
More precisely,
every simple module $V$ in $\Ca_\g^0$ can be obtained as the head of an ordered tensor product of $\bS{k}$'s.
Thus $\{ \bS{k} \}_{k \in \Z}$ plays the same role as fundamental modules in $\Ca_\g^0$, and we call them the \emph{affine cuspidal modules}
(Definition~\ref{def:pbw})
associated with the complete \emph{\pbw} $(\Dd,\hhw)$.
Here, a \pbw $(\Dd,\hhw)$ is a pair of a  duality datum $\Dd$
and $\hhw \seteq \seq{\im_k}_{k \in \Z} \in I_{\fing }{}^\Z$ which
is a sequence in the index set $I_\fing$ of  simple roots of $\fing$
such that $s_{\im_{k+1}}\cdots s_{\im_{k+\ell}}$
is a reduced expression of the longest element $w_0$ of the Weyl group
for every $k$ $(\ell\seteq\ell(w_0))$ (see \S\,\ref{subsec:pbw}).

Furthermore, it is proved that the invariants defined in both categories, such as $\Lambda(M,N)$ and $\de(M,N)$, are preserved under the functor $\F_\Dd$ (see \S\,\ref{sec: SW dailty and T-system} for notations and details).

\medskip

With these recently developed results at hand,
we show in this paper that various subcategories of $\Ca_\g$ provide monoidal categorifications of cluster algebras.
In our results,
there is \emph{no} restriction on the affine type of $\g$; i.e.,  our results hold for an \emph{arbitrary} $U_q'(\g)$.

 The main results of this paper can be summarized in the following three main theorems:
\ben
\item[(\textbf{MT1})]
We give a \emph{vast generalization of T-system}  which implies the determinantial identities among quantum unipotent minors
and the functional relations among $q$-characters of KR-modules simultaneously.

\item[(\textbf{MT2})]
We develop two combinatorial notions:  \emph{admissible chains of $i$-boxes} which provide
commuting families of affine determinantial modules, and \emph{box moves} which describe T-system in a combinatorial way.

\item[(\textbf{MT3})]
We study a family of subcategories $\Ca_{\g}^{[a,b],\Dd,\hhw}$ of $\catCO$, introduced in \cite{KKOP20C},
which contains Hernandez-Leclerc categories, and prove
that $\Ca_{\g}^{[a,b],\Dd,\hhw}$  provides a monoidal
categorification of a cluster algebra, when  $\Dd$ is the duality datum arising from a $\rmQ$-datum $\Qd$.
\ee

Here is a remark on the main theorems. Although a criterion
providing a monoidal categorification of a cluster algebra has been
established in \cite{KKOP19C}, applying the criterion to
subcategories of $\Ca_\g$ is a quite different problem. More
precisely, it is not easy to find a monoidal seed that
satisfies the conditions of the criterion. To overcome this
difficulty, we use the new invariants and the PBW theory to study
affine determinantial modules (see~\eqref{eq: aff det module} below) and their $T$-systems in a general
setting, and then obtain (\textbf{MT1}) and (\textbf{MT2}).  These
theorems allow us to find desired monoidal seeds,
which provides (\textbf{MT3}). Some of these results have been
already announced in \cite{KKOP20D}.

\smallskip
We first define \emph{affine determinantial modules} as follows.
Let $(\Dd,\hhw)$ be a \pbw.
Then we define the {\em affine determinantial module}
\begin{align}\label{eq: aff det module}
 M[a,b] \seteq \hd(\bS{b} \tens \bS{b^-} \tens \cdots \tens \bS{a^+} \tens \bS{a} )
\end{align}
for any interval $[a,b]$ such that $a \le b$ and $\im_a=\im_b$, called an \emph{$i$-box}.
Here,
we set $s^+\seteq\min\{t \mid s<t ,\; \im_t=\im_s \}$ and
$s^-\seteq\max\{t\mid t<s, \; \im_t=\im_s \}$ for $s \in \Z$.
It is proved in Theorem~\ref{thm: commuting ab} that each affine determinantial module is real simple.

\begin{maintheorem}[{\rm [Theorem~\ref{th:Tsystem}]}]  \label{thm:main1} For an \emph{arbitrary} \pbw $(\Dd,\hhw)$  and an $i$-box $[a,b]$, we have an exact sequence
\begin{align}\label{eq: T-system amond det}
0 \to  \sotimes_{ \substack{ \jm \in I_\fing;\; d(\im_a,\jm)=1}}\hs{-1ex} M[a(\jm)^+,b(\jm)^-]  \to   M[a^+,b]
\tens M[a,b^-]    \to   M[a,b] \tens M[a^+,b^-] \to 0,
\end{align}
where  $d(\im,\jm)$  denotes the distance between $\im$ and $\jm \in I_{\fing }$ in the Dynkin diagram $\Dynkin$ of $\fing$,
$s(\jm)^+\seteq\min\{t \mid s \le  t,\; \im_t=\jm \}$ and $s(\jm)^-\seteq\max\{t \mid t \le s,\; \im_t=\jm \}$ for $s \in \Z$ and $\jm \in I_{ \fing}$. \\[1ex]
The exact sequence \eqref{eq: T-system amond det} is also called
\emph{T-system}.
\end{maintheorem}

When an $i$-box $[a,b]$ is contained in  $[1,\ell]$, $M[a,b]$ is isomorphic to $\F_\Dd(\DC_{\tw_0}[a,b])$
and hence ~\eqref{eq: T-system amond det} can be interpreted as a T-system among quantum unipotent minors.
When a \pbw $(\Dd,\hhw)$ is \emph{associated with a
$\rmQ$-datum},
the affine determinantial module $M[a,b]$ is a KR-module
for an \emph{arbitrary} $i$-box $[a,b]$ (Theorem~\ref{Thm: KR as AD}),
and hence we can interpret~\eqref{eq: T-system amond det} as
a well-known T-system among KR-modules.

However, when the $i$-box $[a,b]$ is \emph{not} contained in $[1,\ell]$
nor is the  \pbw $(\Dd,\hhw)$  associated with
any $\rmQ$-datum, the \Ad modules and exact sequences among those modules were not investigated  before
as far as the authors know.

Thus Main Theorem~\ref{thm:main1} can be understood as a vast generalization of the T-system in $R\gmod$ and the T-system
among KR-modules simultaneously.

\medskip
We also develop the combinatorics to construct a commuting family $\cM[](\frakC)$ of \Ad modules and to describe T-system  among \Ad modules as exchange relations in the cluster algebra by mutations. More precisely,
we define the notion of an \emph{admissible chain of $i$-boxes}
$\frakC=(\ci_1,\ldots,\ci_{l})$ as a sequence of $i$-boxes $\ci_k$
(see Definition~\ref{def: ad ch of i box}),
so that $\tc_k   \seteq \bigcup_{1 \le j \le k} \ci_j$
is an interval with $|\tc_k| = k$ for $1 \le k \le l$.
We say that $\tc_l=[a,b]$  is the \emph{range} of $\frakC$
\begin{maintheorem}   [{\rm [Theorem~\ref{thm: admissible chain commuting}, Lemma~\ref{lem: finite sequence are T-equi} and Proposition~\ref{prop:mutation must be}]}] \label{thm:main2}
Let $(\Dd,\hhw)$ be an \emph{arbitrary} \pbw and let $[a,b]$ be any interval such that $a \le b \in \Z \sqcup \{ \pm \infty\}$.
\bna
\item For any admissible chain $\frakC=(\ci_k)_{1 \le k \le b-a+1}$ of $i$-boxes
 with the range $[a,b]$,
$\cM[](\frakC)\seteq \{ M(\ci_k)\}_{1 \le k \le  b-a+1}$ forms a commuting family of \Ad modules.
\item For any admissible chains $\frakC$ and $\frakC'$ with the same range, we can \emph{obtain} $\cM[](\frakC')$ from $\cM[](\frakC)$ by applying a sequence of T-systems
described in terms of newly introduced notion, called \emph{box moves}.
\ee
\end{maintheorem}

Next we construct an initial exchange matrix (equivalently, an initial quiver) and the initial cluster variable modules which are expected to give a cluster algebra structure on $K(\Ca_\g^{[a,b],\Dd,\hhw})$ for $-\infty\le a\le b<+\infty$.
Here, $\Ca_\g^{[a,b],\Dd,\hhw}$ denotes
the monoidal subcategory of $\Ca_\g^0$ generated by $\{ \bS{k} \}_{a \le k \le b}$,
and  the categories explained in~\eqref{it:catA}--\eqref{it:catD}
are  special cases of $\Ca^{[a,b],\Dd,\hhw}_{\g}$.
We  present such data by adapting   the combinatorics given in \cite{GLS11}:  for the subsequence $\seq{\im_k}_{a \le k \le b}$ of $\hhw$, we take the initial cluster variable modules $\{ \cM \}_{a \le s \le b}$
and the initial quiver $\GLS\seteq \GLS(\seq{\im_k}_{a \le k \le b})$ with the vertices $[a,b]$ as follows:
\begin{align*}
\cM \seteq M[s,b(\im_s)^-]  \quad \text{ and } \quad s \to t \quad  \begin{cases}
\text{if $s^-<t^-<s<t$ and $d(\im_s,\im_t)=1$}, \\
\text{or}\\
\text {if $t = s^-$},
\end{cases}
\end{align*}
for $a\le s \ne t \le b$ (see~\eqref{quiver:Omegaminus}).

Since the sequence $\bl [b-s+1,b(\im_{b-s+1})^-] \br_{1 \le s \le b-a+1}$ is an admissible chain $\frakC_-$ of $i$-boxes,
$ \cM[](\frakC_-) \seteq \{ \cM \}_{a \le s \le b } $
is a commuting family of \Ad modules.

Now, for any \pbw $(\Dd,\hhw)$ and any $i$-box $[a,b]$, we
obtain a monoidal seed
\begin{align} \label{eq: la mon seed}
  \seed = ( \{ \cM \}_{a \le s \le b}, \tB_{\GLS} )
\end{align}
where $\tB_{\GLS}$ denotes the incident matrix of $\GLS$.
Then we prove in Theorem~\ref{thm:triangle} that
$$\text{ $\seed$ is  \emph{$\Uplambda$-admissible} (see~\eqref{eq: admissibility})}$$

Recall that the subcategory $\Ca_\g^-$ is a special case of $\Ca^{[a,b],\Dd,\hhw}_{\g}$ and its Grothendieck ring $K(\Ca_\g^-)$
has a cluster algebra structure. Using the cluster algebra structure,
we have the following theorem:
\begin{maintheorem} [{\rm  [Theorem~\ref{th:Main}], [Proposition~\ref{prop: mutation equivalent}]}]\label{thm:main3}
Let $(\Dd,\hhw)$ be a  PBW-pair with $\Dd=\Dd_\Qd$ for some $\rmQ$-datum $\Qd$.  For any admissible chain $\frakC=(\ci_k)_{1\le k \le l }$ for $l \seteq b-a+1 \in \Z_{\ge 1} \sqcup \{ \infty \}$
with the range $\tc_l = [a,b]$ $(a \le b \in \Z \sqcup \{ \pm\infty \})$,
there exists an $\Uplambda$-admissible monoidal seed $\seed$
such that
\bna
\item its set of cluster variable modules is $\cM[](\frakC) \seteq \{ M(\ci_k) \}_{1 \le k\le l}$,
\item its set of frozen variable modules is
$$\{ M[a(\im)^+,b(\im)^-] \mid \im \in \If,   \;
-\infty<a\le a(\im)^+\le b(\im)^-\le b<+\infty \},$$
\item
$ K(\Ca_\g^{[a,b],\Dd,\hhw})$ has a cluster algebra structure with
the initial seed $[\seed]\seteq\st{\;[M(\ci_k)]\;}_{1\le k\le l}$,
and
$\Ca_\g^{[a,b],\Dd,\hhw}$ provides a monoidal categorification of the cluster algebra
$\mathscr{A}([\seed])$.
\end{enumerate}
 In particular, the categories $\Ca_\g^{0}$ and $\Ca_\g^{-}$ provide monoidal categorifications of the cluster algebras.
\end{maintheorem}

For a pair $\frakC$ and $\frakC'$ of admissible chains of $i$-boxes,
the monoidal seed $\cM[](\frakC')$ is obtained from
$\cM[](\frakC)$ by successive mutations.
Note that $\cM[](\frakC)$ consists of KR-modules  when $\tw_0$ in Main Theorem~\ref{thm:main3} is \emph{adapted} to $\Qd$.

To prove Main Theorem~\ref{thm:main3}, we first prove that the quiver $\GLS$ for $[-\infty,\infty]$ and $\Qd$-adapted $\tw_0$ coincides with the quiver $\HL$ in \cite{HL16} by analyzing a sequence $\seq{ (\im_k,p_k) }_{k \in \Z}$
in $\If \times \Z$ associated with a $\rmQ$-datum, called an \emph{admissible sequence} (Proposition~\ref{prop:HL=GLS}). 
Then we prove Main Theorem~\ref{thm:main3} for $\Ca_\g^-$ by using the criterion established in \cite{KKOP19C}.
Finally, we deduce the result on $\Ca_\g^{[a,b],\Dd,\hhw}$ for any interval $[a,b]$
and $\tw_0$ from the one on $\Ca_\g^-$  by developing various methods including mutation equivalence of \Lad seeds.

We conjecture that Main Theorem~\ref{thm:main3} still holds when we weaken the complete \pbw $(\Dd,\hhw)$ with $\Dd=\Dd_\Qd$ in the statement to any \pbw  $(\Dd,\hhw)$ (see Conjecture~\ref{conj: Lambda-ad}).

\mnoi
\textbf{Organization} This paper is organized as follows.
In Section \ref{sec: Qauntum affine},  we give the necessary background on the quantum affine algebras, their representations, the invariants related to R-matrices and the root system
associated with $\Ca_\g$.
In Section \ref{sec: Quiver Hecke algebra}, we give the necessary background on the quiver Hecke algebras, the quantum unipotent coordinate rings and
T-systems among the determinantial modules.
In Section \ref{sec: SW dailty and T-system}, we recall the quantum affine \SW duality functor $\F$ and the results in \cite{KKOP20C}.
Then we define the \Ad modules, study their commuting condition and prove T-system among \Ad modules.
In Section \ref{sec: Admissible chain}, we develop the combinatorics by introducing the notions of
admissible chain of $i$-boxes and box moves.
In Section \ref{sec: subcategories}, we review the notion of $\rmQ$-datum and introduce the notion of admissible sequence in $\If \times \Z$.
Investigating them, we show that
Hernandez-Leclerc categories \eqref{it:catA}--\eqref{it:catD} and KR-modules are special cases of $\Ca_\g^{[a,b],\Dd,\hhw}$ and
\Ad modules, respectively.
In Section \ref{sec: monoidal categorification}, we first review the cluster algebras and the criterion on monoidal categorification by
a monoidal subcategory $\shc$ of $\Ca_\g^0$ established in \cite{KKOP19C}. Then we study the properties of monoidal seeds of several kinds
and prove that the monoidal seed $\seed$  in~\eqref{eq: la mon seed} is \Lad. In the last part,  we prove the
$\GLS$ coincides with the quiver $\HL$ under certain condition.
In Section \ref{sec: main}, we prove Main theorem~\ref{thm:main3}.

\snoi
\textbf{Acknowledgments} The second, third and fourth authors gratefully acknowledge for
the hospitality of RIMS (Kyoto University) during their visit in 2020.

\begin{convention}  Throughout this paper, we keep the following conventions.
\begin{enumerate}
\item For a statement $\mathtt{P}$, $\delta(\mathtt{P})$ is $1$ or $0$ according that
$\mathtt{P}$ is true or not.  In particular, we set $\delta_{i,j}
\seteq \delta(i = j)$ (Kronecker's delta).
\item For a field $\cor$, $a\in\cor$ and $f(z)\in\cor(z)$,
we denote by $\zero_{z=a}f(z)$ the order of zero of $f(z)$ at $z=a$.
\item For $k,l \in \Z$ and $s \in \Z_{\ge1}$, we write $k \equiv_s l$ if $s$ divides $k-l$ and $k \not\equiv_s l$, otherwise.
\item For an object $M$ of an abelian category with finite length, we denote by $\hd(M)$ the head of $M$ and
by $\soc(M)$ the socle of $M$.
\item For a finite set $A$, we denote by $|A|$ the number of elements in $A$.
\item $\ord(\sigma)$ denotes the order of $\sigma$
for an element $\sigma$ of a finite group.
\item
For vertices $i,j$ in a simply-laced Dynkin diagram, $d(i,j)$
denotes the number of edges between $i$ and $j$.
\item $\sym_n$ stands for the symmetric group of degree $n$.
\end{enumerate}
\end{convention}

\vskip 2em

\section{Review on Quantum affine algebras} \label{sec: Qauntum affine}
In this section, we will briefly review the definition of quantum
affine algebras and their representation theory. Then, we will
recall the invariants related to $R$-matrices
which were recently introduced in \cite{KKOP19C}. We refer
to \cite{KKOP19C} for more
details.

\subsection{Quantum affine algebras} Let $(A,P,\sPi,P^\vee,\sPi^\vee)$   be an \emph{affine Cartan datum} consisting of
an \emph{affine Cartan matrix} $A=({\rm a}_{ij})_{i,j\in I}$ with a
finite index set $I$, a \emph{weight lattice} $P$, a set of
\emph{simple roots} $\sPi=\{\ual_i \ | \ i \in I \}\subset P$, a
\emph{coweight lattice} $P^\vee\seteq \Hom_{\Z}(P,\Z)$ and a set of
\emph{simple coroots} $\sPi^\vee=\{h_i \ | \ i \in I \}\subset
P^\vee$. We have $\lan h_i,\ual_j \ran = {\rm a}_{ij}$ for all $i,j
\in I$ where $\lan \ , \  \ran\col P^\vee \times P \to \Z$ is the
canonical pairing. We choose $\st{\sLa_i}_{i \in I}$ such that $\lan
h_j ,\sLa_i \ran=\delta_{i,j}$ for $i,j \in I$ and call them the
\emph{fundamental weights}.

 We also take the \emph{imaginary root} $\updelta = \sum_{i \in I} \mathsf{u\ms{1mu}}_i \ual_i$ and
the \emph{center} $c=\sum_{i \in I} \mathsf{c}_i h_i$ such that $\{
\la \in \soplus_{i\in I}\Z\ual_i\mid   \lan h_i,\la \ran = 0 \text{
for every } i \in I \}=\Z\updelta$ and $\{ h \in \soplus_{i\in I}\Z
h_i \mid  \lan h,\ual_i \ran = 0 \text{ for every } i \in I \}=\Z c$
(see \cite[Chapter 4]{Kac}). We set $P_\cl  \seteq P / (P\cap\Q
\updelta)$ and call it the \emph{classical weight lattice}. We
choose $\rho \in P$ (resp.\ $\rho^\vee \in P^\vee$) such that $\lan
h_i,\rho \ran=1$ (resp.\ $\lan \rho^\vee,\ual_i\ran =1$) for all $i
\in I$.
Set $\h \seteq \Q \tens_{\Z} P^\vee$. Then there exists a
non-degenerate symmetric bilinear form $( \ , \ )$ on $\h^*$
satisfying
$$\lan h_i,\la \ran= \dfrac{2(\ual_i,\la)}{(\ual_i,\ual_i)} \ \text{ and } \ \lan c,\la \ran = (\updelta,\la)  \ \ \text{ for any $i \in I$ and $\la \in \h^*$}.$$

We denote by $\g$ the \emph{affine Kac-Moody algebra} associated
with $(A,P,\sPi,P^\vee,\sPi^\vee)$ and by $W \seteq \lan s_i \ | \ i
\in I \ran \subset GL(\h^*)$ the \emph{Weyl group} of $\g$ where
$$  s_i \la =  \la - \lan h_i,\la \ran \ual_i \quad \text{ for } \la \in P.$$
We will use the convention in~\cite{Kac} to choose $0\in I$ except
$A^{(2)}_{2n}$-type, in which case we take the longest simple root
as $\ual_0$, and $B_2^{(1)}$, $A_3^{(2)}$, and $E_k^{(1)}$
($k=6,7,8$) types, in which we take the following Dynkin diagrams
$\Dynkin_\g$:
\begin{equation} \label{Eq: DD}
\begin{aligned}
\akew[1ex]& A^{(2)}_{2n} :   \xymatrix@C=4ex@R=3ex{
  *{\circ}<3pt> \ar@{<=}[r]_<{n \ } & *{ \circ }<3pt> \ar@{-}[r]_<{n-1}  & *{ \circ }<3pt> \ar@{-}[r]_<{n-2} & \cdots \ar@{-}[r]_<{ }   &*{\circ}<3pt> \ar@{-}[l]^<{ \ \ 1} &*{ \circ }<3pt>  \ar@{=>}[l]^<{ \ \ 0}  } \hs{3ex}
B^{(1)}_{2} :   \xymatrix@C=4ex@R=3ex{
  *{\circ}<3pt> \ar@{=>}[r]_<{0 \ } & *{ \circ }<3pt>  \ar@{<=}[l]^<{ \ \ 2} &*{ \circ }<3pt>  \ar@{=>}[l]^<{ \ \ 1}  } \hs{3ex}
A^{(2)}_{3} :   \xymatrix@C=4ex@R=3ex{
  *{\circ}<3pt> \ar@{<=}[r]_<{0 \ } & *{ \circ }<3pt>  \ar@{=>}[l]^<{ \ \ 2} &*{ \circ }<3pt>  \ar@{<=}[l]^<{ \ \ 1}  }
\\
&
 E^{(1)}_6 :  \raisebox{2.3em}{\xymatrix@C=3.4ex@R=2ex{ && *{\circ}<3pt>\ar@{-}[d]^<{0}  \\ && *{\circ}<3pt>\ar@{-}[d]^<{2} \\
*{ \circ }<3pt> \ar@{-}[r]_<{1}  &*{\circ}<3pt> \ar@{-}[r]_<{3} &*{
\circ }<3pt> \ar@{-}[r]_<{4} &*{\circ}<3pt> \ar@{-}[r]_<{5}
&*{\circ}<3pt> \ar@{-}[l]^<{\ \ 6}}} \qquad
E^{(1)}_7 : \raisebox{1.3em}{\xymatrix@C=3.4ex@R=3ex{ && & *{\circ}<3pt>\ar@{-}[d]^<{2} \\
*{ \circ }<3pt> \ar@{-}[r]_<{0}  & *{ \circ }<3pt> \ar@{-}[r]_<{1}
&*{\circ}<3pt> \ar@{-}[r]_<{3} &*{ \circ }<3pt> \ar@{-}[r]_<{4}
&*{\circ}<3pt> \ar@{-}[r]_<{5} &*{\circ}<3pt> \ar@{-}[r]_<{6}
&*{\circ}<3pt>
\ar@{-}[l]^<{ \ \ 7} } } \allowdisplaybreaks \\
& E_8^{(1)}  :  \raisebox{1.3em}{\xymatrix@C=3.4ex@R=3ex{ && *{\circ}<3pt>\ar@{-}[d]^<{2} \\
*{ \circ }<3pt> \ar@{-}[r]_<{1}  &*{\circ}<3pt> \ar@{-}[r]_<{3} &*{
\circ }<3pt> \ar@{-}[r]_<{4} &*{\circ}<3pt> \ar@{-}[r]_<{5}
&*{\circ}<3pt> \ar@{-}[r]_<{6} &*{\circ}<3pt> \ar@{-}[r]_<{7}
&*{\circ}<3pt> \ar@{-}[l]^<{ \ \ 8} &*{\circ}<3pt> \ar@{-}[l]^<{0} }
}
\end{aligned}
\end{equation}

We define $\g_0$ to be the subalgebra of $\g$ generated by the
\emph{Chevalley generators} $e_i$, $f_i$  and $h_i$ for $i \in I_0
\seteq I \setminus \{ 0 \}$ and $W_0$ to be the subgroup of $W$
generated by $s_i$ for $i \in I_0$. Note that $\g_0$ is a
finite-dimensional simple Lie algebra and $W_0$ contains the longest
element $w_0$.

Let $q$ be an indeterminate and  let $\ko$ be the algebraic closure
of the subfield $\C(q)$ in the algebraically closed field
$\corh\seteq\displaystyle\bigcup\nolimits_{m >0}\C \ls q^{1/m} \rs$. For $i\in I$, we set
$q_i = q^{(\ual_i,\ual_i)/2}$.

Let us denote
by $U_q'(\g)$ the \emph{quantum affine algebra} associated with an
affine Cartan datum $(A,P,\sPi,P^\vee,\sPi^\vee)$ generated by
$e_i,f_i,q_i^{\pm h_i}$ $(i \in I)$ over $\ko$.

\subsection{Finite-dimensional representations} We say that a module $M$ over $U_q'(\g)$ is \emph{integrable} if (i) $M$ decomposes into $P_{\cl}$-weight spaces; i.e., $M = \bigoplus_{\la \in P_\cl} M_\la$ where
$M_\la = \{ u \in M \ | \ q_i^{h_i}u = q_i^{\lan h_i,\la \ran} \}$,
and (ii) the action of $e_i$ and $f_i$ on $M$ is nilpotent $(i \in I)$.
We denote by $\Ca_\g$ the category of
finite-dimensional integrable $U_q'(\g)$-modules.

A simple module $M$ in $\uqm$ contains a non-zero vector $u$ of
weight $\lambda\in P_\cl$ such that (i) $\langle h_i,\lambda \rangle
\ge 0$ for all $i \in I_0$, (ii) all the weights of $M$ are contained
in $\lambda - \sum_{i \in I_0} \Z_{\ge 0} \cl(\ual_i)$, where
$\cl\colon P\to P_\cl$ denotes the canonical projection. Such a
$\la$ is unique and $u$ is unique up to a constant multiple. We call
$\lambda$ the \emph{dominant extremal weight} of $M$ and $u$ a
\emph{dominant extremal weight vector} of $M$.

For an indeterminate $z$ and a $U_q'(\g)$-module $M$, let us denote
by  $M_z $ the $U_q'(\g)$-module $\ko[z^{\pm1}]\tens M$
defined by
$$ e_i(u_z)=z^{\delta_{i,0}}(e_iu)_z, \quad f_i(u_z)=z^{-\delta_{i,0}}(f_iu)_z, \quad q_i^{h_i}(u_z)=(q_i^{h_i}u)_z,$$
where $u_z$ denotes $1\tens z$ for $u\in M$.

For $x \in \ko^\times$, we define $$M_x \seteq M_z / (z -x)M_z.$$ We
call $x$ a {\it spectral parameter}. Note that, for a module $M$ in
$\uqm$ and $x \in \ko^\times$, $M_x$ is also contained in $\uqm$.
The functor $\mathsf{T}_x$ defined by $\mathsf{T}_x(M)=M_x$ is an
autofunctor of $\uqm$ which commutes with tensor products.

For each $i \in I_0$, we set
\[
\varpi_i \seteq
\gcd(\mathsf{c}_0,\mathsf{c}_i)^{-1}\cl(\mathsf{c}_0\sLa_i-\mathsf{c}_i
\sLa_0) \in P_{\cl}.
\]
Then there exists a unique simple module $V(\varpi_i)$ in $\uqm$,
called the \emph{fundamental module} of weight $\varpi_i$,
satisfying  certain  conditions
(see, e.g,~\cite[\S 5.2]{Kas02}).

For a $U_q'(\g)$-module $M$, we denote by $\overline{M}=\{ \bar{u}
\mid u \in M \}$ the $U_q'(\g)$-module whose module structure is
given as $x \bar{u} \seteq \overline{\ms{2mu}\overline{x}
u\ms{2mu}}$ for $x \in U_q'(\g)$. Here
$\overline{\akew[1ex]\ake[1.5ex]}$ is the ring automorphism of
$\usq$ such that $\ol{q}=q^{-1}$, $\ol{e_i}=e_i$, $\ol{f_i}=f_i$ and
$\ol{q^h}=q^{-h}$.

Then we have
\begin{align*}
\overline{M_a} \simeq (\overline{M})_{\,\overline{a}}, \qquad\qquad
\overline{M \otimes N} \simeq \overline{N} \otimes \overline{M}.
\end{align*}

Set
$\pi_i=\gcd(\mathsf{c}_0,\mathsf{c}_i)^{-1}(\mathsf{c}_0\sLa_i-\mathsf{c}_i
\sLa_0) \in P$. Let $m_i$ be a positive integer such that
\[
W \pi_i\cap \bl\pi_i+\Z\updelta\br=\pi_i+\Z m_i\updelta.
\]
Note that $m_i=(\ual_i,\ual_i)/2$ in the case where $\g$ is the dual
of an untwisted affine algebra, and $m_i=1$ otherwise. Then, for
$x,y\in \ko^\times$, we have $V(\varpi_i)_x \simeq V(\varpi_i)_y$ if
and only if $x^{m_i}=y^{m_i}$ (\cite[\S 1.3]{AK97}).

\smallskip

For simple modules $M$ and $N$ in $\Ca_\g$, we say that $M$ and $N$
{\em commute} if $M\tens N\simeq N\tens M$. We say that they
\emph{strongly commute} (or $M$ \emph{strongly commutes with} $N$)
if $M \tens N$ is simple. Note that $M$ and $N$ commute as soon as
they strongly commute. We say that a simple module $L$ in $\Ca_\g$
is \emph{real} if $L$ strongly commutes with itself, \textit{i.e.},
if $L \tens L$ is simple.

Recall that the simple objects in $\Ca_\g$ are parameterized by
$I_0$-tuples of polynomials $\mathcal{P}=\st{\mathcal{P}_{i} (u)}_{i
\in I_0}$, where $\mathcal{P}_{i}(u) \in \ko[u]$ and
$\mathcal{P}_{i}(0)=1$ (\cite{CP95,CP98}, see also
\cite[(1.5)]{KKK15B} for the choice of a function $I_0 \to \{ \pm1
\}$). We denote by $\mathcal{P}^V_{i}(u)$ the polynomials
corresponding to a simple module $V$, and call
$\mathcal{P}^V_{i}(u)$  the \emph{Drinfel'd polynomials} of $V$.
 The Drinfel'd polynomials $\mathcal{P}^V_{i}(u)$ are determined by the eigenvalues
of the simultaneously commuting actions of some \emph{Drinfel'd
generators} of $\uqpg$ on a subspace of $V$ (see,
e.g., \cite{CP95} for more details).
Note that
$\mathcal{P}^{V_x}_{i}(u)=\mathcal{P}^V_{i}(x^{m_i}u)$ for any $x\in\cor^\times$.

The \emph{Kirillov-Reshetikhin {\rm (KR)} module}, usually denoted
by $W^{(k)}_{m,a}$ for $k \in I_0$, $a \in \ko^\times$ and $m \ge
1$, is a simple module of dominant extremal weight
$m\varpi_k$ in $\Ca_\g$ whose Drinfel'd
polynomials $\mathcal{P} = \st{\mathcal{P}_i }_{i \in I_0}$ are
(see~\cite{Her10,Nak04a,OS08} and \cite[Remark 3.3.1]{N21}
for more detail)
\[
\mathcal{P}_i(u) = \delta_{i,k} (1-au) (1-a{\cq_k}{}^2u) (1-a\cq_k{}^4u)
\cdots (1-a\cq_k{}^{2m-2}u) + (1 - \delta_{i,k}) \quad (i\in I_0),
\]
where
$$
\cq_k =\bc q_k &\text{unless $\g = A_{2n}^{(2)}$ and $k = n$,}\\
q &\text{if $\g = A_{2n}^{(2)}$ and $k = n$.}\ec
$$
Then the \KR module $V(i^m) \seteq W^{(i)}_{m,(-\Check{q}_i)^{1-m}}$
is simple and bar-invariant. With the terminology of
$V(\varpi_i)_z$, the module $V(i^m) $ can be described as follows:
\eq\label{eq: KR1}
&\akew[4ex]&V(i^m) \simeq
\begin{cases}
\hd \left(V(\varpi_i)_{(-q_i)^{m-1}}  \tens \hspace{-.3ex} V(\varpi_i)_{(-q_i)^{m-3}} \tens  \cdots    \tens V(\varpi_i)_{(-q_i)^{1-m}}\right)\hs{-3ex} &\text{ if $\g$ is untwisted,}  \\
\hd \left(V(\varpi_i)_{(-q)^{m-1}}   \tens \hspace{-.3ex}
V(\varpi_i)_{(-q)^{m-3}} \tens  \cdots    \tens
V(\varpi_i)_{(-q)^{1-m}}\right)  &  \text{ otherwise},
\end{cases}
\eneq
for $i \in I_0$ and $m \in \Z_{\ge 1}$ (see ~\cite{Her10}).

Note that the category $\uqm$ is {\it rigid}\,, i.e.,  every module $M$
in $\uqm$ has the right dual $\D M$ and the left dual $\D^{-1}M$.
Hence we have the evaluation morphisms
$$M\tens \D M\to\one\quad \text{and} \quad  \D^{-1} M\tens M\to\one.$$
We extend this to $\D^k$ for $k \in \Z$. In particular, the duals of
$V(i^m)_x$ $(x\in \cor^\times)$ are given as follows:
\begin{equation*}
 \D \bigl(V(i^m)_x \bigr) \simeq   V((i^*)^m)_{p^*x}\qtq
\D^{-1}\bigl( V(i^m)_x\bigr)  \simeq   V((i^*)^m)_{(p^*)^{-1}x}.
\end{equation*}
Here $p^* \seteq (-1)^{\langle \rho^\vee ,\updelta
\rangle}q^{\ang{c,\rho}}$ and $i^*\in I_0$ is defined by
$\ual_{i^*}=-w_0\,\ual_i$.

We set
\begin{align} \label{eq: sig}
\upsigma(\g) \seteq I_0 \times \ko^\times/ \sim
\end{align}
where the equivalence relation is given by $$   (i,x) \sim (j,y) \iff
V(\varpi_i)_x \simeq V(\varpi_i)_y.$$
We denote by $[(i,a)]$ the
equivalence class of $(i,a)$ in $ \sig$.  When no confusion arises,
we simply write $(i,a)$ for the equivalence class $[(i,a)]$.

The set $\sig$ has a graph structure: we join $(i,a)$ and $(j,b)$
when $V(i)_a$ and $V(j)_b$ do not commute. We choose a connected
component $\sigZ$ of $\sig$.

We denote by $\cato$ the full smallest full subcategory of $\catg$
such that it contains $V(\vpi_i)_a$ ($(i,a)\in\sigZ$) and is stable
under taking subquotients, extensions and tensor products (see
\cite{HL10,KKO18} and Definition~\ref{def: Category [a,b]} below).

\subsection{R-matrices and invariants} \label{subsec:invariants}
For modules $M$ and $N$ in $\Ca_\g$, there exists a $\ko(\hspace{-.3ex}(z)\hspace{-.3ex}) \tens \uqpg$-module isomorphism, denoted by $\Runiv_{M \otimes N_z}$ and called the \emph{universal $R$-matrix} of $M$ and $N$:
\begin{align*}
\Runiv_{M,N_z} \colon
\ko(\hspace{-.3ex}(z)\hspace{-.3ex})\tens_{\ko[z^{\pm1}]} (M \tens
N_z) \To \ko\ls z \rs \tens_{\ko[z^{\pm1}]} (N_z\tens M)
\end{align*}
satisfying certain properties (see \cite{Kas02} and \cite[Appendices
A and B]{AK97}).

For non-zero modules $M$ and $N$ in $\Ca_\g$, if there exists $f(z)
\in \ko\ls z \rs^\times $ such that
$$ f(z)\Runiv_{M,N_z}(M \tens N_z) \subset N_z \tens M,$$
then we say that $\Runiv_{M,N_z}$ is \emph{rationally renormalizable}.
 In the rationally renormalizable case, one can choose
$c_{M,N}(z) \in \ko(\hspace{-.3ex}(z)\hspace{-.3ex})^\times$ as
$f(z)$ such that,
 for any $x \in \ko^\times$, the specialization
of $\Rren_{M,N_z} \seteq c_{M,N}(z)\Runiv_{M,N_z} \col M \otimes N_z
\to N_z \otimes M$ at $z=x$
$$  \Rren_{M,N_z}\bigm\vert_{z=x} \colon M \otimes N_x  \to N_x \otimes M$$
does not vanish.  Such $\Rren_{M,N_z}$ and $c_{M,N}(z)$ are unique
up to a multiple of $\cz^\times = \bigsqcup_{\,n \in \Z}\ko^\times
z^n$. We call $c_{M,N}(z)$  the \emph{renormalizing coefficient}. We
write $\rmat{M,N} \seteq \Rren_{M,N_z}\bigm\vert_{z=1}$ and call it
\emph{$R$-matrix}. The $R$-matrix $\rmat{M,N}$ is well-defined up to
a constant multiple when $\Runiv_{M,N_z}$ is rationally
renormalizable. By the definition, $\rmat{M,N}$ never vanishes.

For simple modules $M$ and $N$ in $\Ca_\g$, let $u$ and $v$ be
dominant extremal weight vectors of $M$ and $N$, respectively. Then
there exists $a_{M,N}(z) \in
\ko(\hspace{-.3ex}(z)\hspace{-.3ex})^\times$ such that
$$\Runiv_{M,N_z}(u \tens v_z) = a_{M,N}(z) (v_z \tens u).$$
Then $\Rnorm_{M,N_z}\seteq
a_{M,N}(z)^{-1}\Runiv_{M,N_z}\big\vert_{\;\ko(z)\otimes_{\ko[z^{\pm1}]}
( M \otimes N_z) }$ induces a unique $\ko(z)\tens\uqpg$-module
isomorphism
\begin{equation*}
\Rnorm_{M, N_z} \colon \ko(z)\otimes_{\ko[z^{\pm1}]} \big( M \otimes
N_z\big) \longisoto\ko(z)\otimes_{\ko[z^{\pm1}]}  \big( N_z \otimes
M \big)
\end{equation*}
sending   $ u  \otimes v_z$  to  $v_z\otimes u$.
 Hence, the universal $R$-matrix
$\Runiv_{M,N_z}$ is \rr. We call $a_{M,N}(z)$ the {\it universal
coefficient} of $M$ and $N$, and $\Rnorm_{M,N_z}$ the {\em
normalized $R$-matrix}.
Note that  $\ko(z)\otimes_{\ko[z^{\pm1}]} ( M \otimes N_z)$ is a
simple $\ko(z) \tens \uqpg$-module by \cite[Proposition 9.5]{Kas02}.

Let $d_{M,N}(z) \in \ko[z]$ be a monic polynomial of the smallest
degree such that the image of $d_{M,N}(z) \Rnorm_{M,N_z}(M\tens
N_z)$ is contained in $N_z \otimes M$. We call $d_{M,N}(z)$ the {\em
denominator of $\Rnorm_{M,N_z}$}. Then we have
\begin{equation*}
\Rren_{M,N_z}= d_{M,N}(z)\Rnorm_{M,N_z} \col M \otimes N_z \To N_z
\otimes M \qt{up to a multiple of $\cz^\times$.}
\end{equation*}
Hence, we have
\begin{align*} \Rren_{M,N_z} =a_{M,N}(z)^{-1}d_{M,N}(z)\Runiv_{M,N_z}
\quad \text{and} \quad  c_{M,N}(z)= \dfrac{d_{M,N}(z)}{a_{M,N}(z)}
\end{align*}
up to a multiple of $\ko[z^{\pm1}]^\times$.

Note that for simple modules $M,N$, we have
\begin{align*}
\Hom_{\cor[z^{\pm1}]\tens\uqpg}(M\tens N_z,N_z\tens M)
=\cor[z^{\pm1}]\Rren_{M,N_z}.
\end{align*}

Similarly, there exists a $\cor[z^{\pm1}]\tens\uqpg$-linear
homomorphism $\Rren_{M_z,N}\col M_z\tens N\To N\tens M_z$ such that
\eq \Hom_{\cor[z^{\pm1}]\tens\uqpg}(M_z\tens N,N\tens M_z)
=\cor[z^{\pm1}]\Rren_{M_z,N}. \eneq

\begin{remark}  \hfill
\bna
\item The denominator formulas and universal coefficients were studied and computed between fundamental modules in \cite{AK97, DO94, Fu19, KKK15B, Oh15, OT19}, and between KR-modules in \cite{OT19B}.
\item For $(i,x)$ and $ (j,y) \in \sig$, we put $d$ many arrows  from $(i,x)$ to $(j,y)$, where $d$ is the order of zeros of $d_{ V(\varpi_i), V(\varpi_j) } \allowbreak ( z_{V(\varpi_j)} / z_{V(\varpi_i)}  )$
at $  z_{V(\varpi_j)} / z_{V(\varpi_i)}  = y/x$. Thus, $\sig$ has a
quiver structure. \ee
\end{remark}

We set
$$ \tp \seteq p^{*2}=q^{2\ang{c,\rho}} \quad \text{  and } \quad \vphi(z) \seteq \prod_{s\in\Z_{\ge0}}(1-\tp^sz)
=\sum_{n=0}^\infty\hs{.3ex}\dfrac{(-1)^n\tp^{n(n-1)/2}}{\prod_{k=1}^n(1-\tp^k)}\;z^n
 \in \cor[[z]].$$
\begin{definition}
We define the subset $\G$ of $\cor \ls z \rs^\times$ as follows:
\begin{align*}
\G \seteq \left\{ cz^m \prod_{a \in \ko^\times} \varphi(az)^{\eta_a}
\ \left|  \
\parbox{40ex}{$c \in \ko^\times$, $ m \in \Z$ , \\[.5ex]
$\eta_a \in \Z$ vanishes except finitely many $a$'s}
\right. \right\}.
\end{align*}
\end{definition}
Note that $\G$ forms a group with respect to the multiplication and
$\ko(z)^\times\subset \G$.

\begin{proposition} [{\cite[Proposition 3.2]{KKOP19C}}]
Let $M$ and $N$ be non-zero modules in $\uqm$. \bnum
\item If $\Runiv_{M,N_z}$ is \rr, then $c_{M,N}(z)$ belongs to $\G$.
\item If $M$ and $N$ are simple, then  $a_{M,N}(z)$ as well as $c_{M,N}(z)$ is contained in $\G$.
\end{enumerate}
\end{proposition}

For a subset $S$ of $\Z$, we set
$$ \  \tp^{S} \seteq \{ \tp^k \ | \ k \in S\}.$$

By taking $S=\Z$ or $\Z_{\le 0}$, the following group homomorphisms
from $\G$ to $\Z$ were introduced in \cite[Section 3]{KKOP19C}:
\begin{align*}
\Deg \col   \Bg \to  \Z \quad \text{ and } \quad \Di \col   \Bg \to
\Z,
\end{align*} which are defined by
$$
\Deg(f(z)) = \sum_{a \in \tp^{\,\Z_{\le 0}} }\eta_a - \sum_{a \in
\tp^{\,\Z_{> 0}} } \eta_a
 \quad \text{ and } \quad \Di(f(z)) = \sum_{a \in \tp^{\,\Z}} \eta_a
$$
for $f(z)=cz^m \prod \varphi(az)^{\eta_a} \in \Bg$.

Note that
\begin{align} \label{eq: deg f}
\Deg(f(z))=2\zero_{z=1}f(z) \qt{for
$f(z)\in\cor(z)^\times\subset\G$}
\end{align}
(see \cite[Lemma 3.4]{KKOP19C}).

\begin{definition}[{\cite[Definition 3.6, Definition 3.14]{KKOP19C}}]\label{def: Lams an de}
 Let $M,N \in \uqm$.
\ben
\item  If $\Runiv_{M,N_z}$ is \rr, we define the integers  $\Lambda(M,N)$  and $\Lambda^\infty(M,N)$ as follows:
$$\Lambda(M,N)=\Deg(c_{M,N}(z))
\qtq\Lambda^\infty(M,N)=\Deg^\infty(c_{M,N}(z)).$$
\item For simple modules $M$ and $N$ in $\uqm$, we define $\de(M,N)$ by
$$ \de(M,N)= \dfrac{1}{2}\bl\Lambda(M,N) + \Lambda( \D^{-1} M,N )\br.$$
\end{enumerate}
\end{definition}

\begin{proposition}[{\cite[Proposition 3.16, Corollary 3.19]{KKOP19C}}] \label{prop: de(M,N)}
For simple modules $M$ and $N$ in $\uqm$, we have
\begin{equation} \label{eq: de MN deg}
\de(M,N)=\zero_{z=1}\bl d_{M,N}(z)d_{N,M}(z^{-1})\br.
\end{equation}
In particular, we have
$$ \de(M,N)\in\Z_{\ge0} \quad \text{ and } \quad\de(M,N)   =  \dfrac{1}{2}\bl\Lambda(M,N) + \Lambda(N,M) \br = \de(N,M).$$
\end{proposition}

\begin{corollary}[{\cite[Corollary 3.17]{KKOP19C}}]\label{cor:commde}
Let $M$ and $N$ be simple modules in $\uqm$. Assume that one of them
is real. Then $M$ and $N$ strongly commute if and only if
$\de(M,N)=0$.
\end{corollary}

Interestingly, the invariants $\La$ and $\Li$ are calculated by
$\de$ as follows:

\begin{proposition}[{\cite[Proposition 3.\,22]{KKOP19C}, \cite[Proposition 2.\,16]{KKOP20C}}] \label{prop: de and Lambda}
For simple modules $M$ and $N$ in $\uqm$,  we have the followings:
\bnum
\item $\Lambda(M,N)= \dsum_{k \in \Z} (-1)^{k+\delta(k<0)} \de(M,\D^{k}N)$,
\item $\Lambda^\infty(M,N)= \dsum_{k \in \Z} (-1)^{k} \de(M,\D^{k}N)$,
\item $\zero_{z=1}c_{M,N}(z)=\dsum_{k=0}^\infty(-1)^k\de(M,\D^kN)$.
\end{enumerate}
\end{proposition}

\begin{proposition}[{\cite[Proposition 3.18,  Corollary 3.20]{KKOP19C}}]\label{prop: Lambdas}
For simple modules $M$ and $N$ in $\Ca_\g$, we have the followings:
\bnum
\item $\La(M,N)=\La(\D^{-1}N,M)=\La(N,\D M)$.
\item If $M$ is real, then we have $\La(M,M)=0$.
\end{enumerate}
\end{proposition}
We conjecture that $\La(M,M)=0$ holds for an arbitrary simple module $M$.

\Prop[{\cite[Corollary 3.11]{KKKO15}, \cite[Proposition
2.11]{KKOP19C}}] \label{prop:rcomp}\hfill
\bnum
\item \label{item1} Let $M_k$ be a module in $\uqm$ $(k=1,2,3)$, and let $\varphi_1\col L\to M_2\tens M_3$ and $\varphi_2\col M_1\tens M_2\to L'$
be {\em non-zero} morphisms. Assume further that $M_2$ is a simple
module. Then the composition
\begin{align*}
M_1\tens L \To[M_1\tens \varphi_1] M_1\tens M_2\tens
M_3\To[\varphi_2\tens M_3] L'\tens M_3
\end{align*}
does not vanish.
\item \label{item2} Let $M$, $N_1$ and $N_2$ be non-zero modules in $\uqm$, and assume that
$\Runiv_{N_k,M_z}$ is \rr for $k=1,2$. Then $\Runiv_{N_1\tens
N_2,M_z}$ is \rr, and we have
$$\dfrac{c_{N_1,M}(z)c_{N_2,M}(z)}{c_{N_1\tens N_2, M}(z)}\in \cz.$$
If we assume further that $M$ is simple, then we have
$$   c_{N_1 \otimes N_2,M}(z)\equiv c_{N_1,M}(z) c_{N_2,M}(z) \mod \cz^\times $$
and the following diagram commutes up to a constant multiple{\rm:}
\begin{align*}
\xymatrix@C=12ex { N_1\tens N_2 \tens M \ar[r]_{N_1 \tens
\rmat{N_2,M}}\ar@/^2pc/[rr]|-{\rmat{N_1\tens N_2,\; M}} &N_1\tens
M\tens N_2\ar[r]_{\rmat{N_1,M}\tens\, N_2}& M \tens N_1\tens N_2.
}\label{eq:tensr}
\end{align*}
\ee \enprop

\begin{theorem}[{\cite{KKKO15}}]\label{thm: KKKo15 main} Let $M$ and $N$ be simple modules in $\uqm$ and assume that one of them is real. Then
\bna
\item \label{runiq} $\Hom(M\tens N,N\tens M)=\cor\,\rmat{M,N}$.
\item \label{i2} $M \otimes N$ and $N \otimes M$ have simple socles and simple heads.
\item \label{i3} Moreover, ${\rm Im}(\rmat{M,N})$ is isomorphic to the head of $M \otimes N$ and the socle of $N \otimes M$.
\item \label{i4} $M \otimes N$ is simple whenever its head and its socle
are isomorphic.
\end{enumerate}
\end{theorem}

For modules $M$ and $N$ in $\uqm$, we denote by $M \hconv N$ and $M
\sconv N$ the head and the socle of $M \otimes N$, respectively.

\begin{proposition}[{\cite[Proposition 3.2.17]{KKKO18} (see also \cite[Lemma 7.3]{KO18})}] \label{prop: length 2}
Let $M$ and $N$ be simple modules in $\Ca_\g$. Assume that one of
them is real and $\de(M,N)=1$. Then we have an exact sequence
\[
0 \to M \sconv N \to M\tens N\to M\hconv N \to 0.
\]
In particular, $M\tens N$ has composition length $2$.
\end{proposition}

\begin{definition} [{cf.\ \cite[Definition 2.5]{KK18}}]
A sequence $(L_1,\ldots,L_r)$ of real simple modules  in $\uqm$ is
called a \emph{normal sequence} if the composition of the
$R$-matrices \eqn \rmat{L_1,\ldots,L_r}\seteq
\displaystyle\prod_{1\le i <k \le r} \rmat{L_i,L_k}
=&&(\rmat{L_{r-1},L_r})  \circ \cdots \circ (\rmat{L_2,L_r}\circ
\cdots \circ \rmat{L_2,L_3})  \circ (\rmat{L_1,L_r} \circ \cdots
\circ \rmat{L_1,L_2})
\\
  &&: L_1\tens \cdots \tens L_r \longrightarrow L_r \tens \cdots \tens  L_1
\eneqn does not vanish.
\end{definition}

\begin{lemma} [{\cite[Lemma 4.15, Lemma 4.16]{KKOP19C}}] \label{eq: normal head}
Let $(L_1,\ldots,L_r)$ be a normal sequence of real simple modules in
$\uqm$. Then the image of $\rmat{L_1,\ldots,L_r}$ is simple and
coincides with the head of $L_1\tens \cdots \tens L_r$ and also with
the socle of $ L_r \tens \cdots \tens  L_1$.
Moreover we have
\bnum
\item
$(L_2,\ldots, L_r)$ is a normal sequence and
$\La\bl L_1,\Im(\rmat{L_2,\ldots,L_r})\br=\displaystyle\sum_{k=2}^{r}\La(L_1,L_k)$,
\item
$(L_1,\ldots, L_{r-1})$ is a normal sequence and
$\La\bl \Im(\rmat{L_1,\ldots,L_{r-1}}), L_r\br=\displaystyle\sum_{k=1}^{r-1}\La(L_k,L_r)$.
\ee
\end{lemma}

\begin{lemma} [{\cite[Lemma 4.17]{KKOP19C}}]\label{prop: L N* normal}
For real simple modules $L$, $M$ and $N$ in $\uqm$,  the triple
$(L,M,N)$ is a normal sequence if one of the following three
conditions holds: \bnum
\item $L$ and $M$ strongly commute,
\item $M$ and $N$ strongly commute,
\item $L$ and $\D^{-1}N$ strongly commute.
\ee
\end{lemma}

\begin{definition}[{\cite[Definition 5.1]{KKOP20C}}] Let $(M,N)$ be an ordered pair of simple modules in $\Ca_\g$. We call it {\it unmixed} if
$$   \de(\D M,N) =0,$$
and  {\it strongly unmixed} if
$$   \de(\D^k M,N) =0  \quad \text{ for any } k \in \Z_{\ge1}.$$
\end{definition}

\begin{definition}
A sequence $(M_{s},M_{s-1},\ldots,M_1)$  of real simple modules over
$U_q'(\g)$ is {\it $($strongly$)$ unmixed} if $(M_k,M_i)$ is
(strongly) unmixed for all $s\ge k>i \ge 1$.
\end{definition}

\begin{proposition} [{\cite[Lemma 5.2]{KKOP20C}}] \label{prop: strongly unmixed}
For a strongly unmixed pair $(M,N)$, we have
$$  \La^\infty(M,N) = \La(M,N)  = \La (N,\D M).$$
\end{proposition}

\Prop[{\cite[Lemma 5.3]{KKOP20C}}] \label{prop: s un and normal} Any
unmixed sequence of real simple modules is a normal sequence.
\enprop

\begin{proposition}[{\cite[Proposition 4.2]{KKOP19C}}] \label{pro:subquotient}
For simple modules $L$, $M$ and $N$, we have
\begin{align*}
&\de(S,L)\le\de(M,L)+\de(N,L)
\end{align*}
for any simple subquotient $S$ of $M\otimes N$.
\end{proposition}

\Prop\label{prop: s unmixed} Let $(M_k,M_{k-1},\ldots,M_1)$ be a
strongly unmixed sequence. Then
$$(\hd(M_k \tens M_{k-1} \tens \ldots M_{s}) ,\hd(M_{s-1} \tens  \ldots M_{2} \tens M_{1}))$$
is strongly unmixed for any $ 1 <s \le k$. \enprop

We recall the following criteria of reality of simple modules.

\begin{proposition}[{\cite[Proposition 3.2.20, Corollary 3.2.21]{KKKO18}, \cite[Proposition 4.9]{KKOP19C}}] \label{prop:real}
Let $X,Y,M$ and $N$ be simple modules in $\uqm$.
 Assume that there is an exact sequence
$$0 \to X \to M \tens N \to Y \to 0, $$
and $X \tens N$ and $Y \tens N$ are simple.
\bnum
\item  If  $X \tens N \not\simeq Y \tens N$, then $N$ is a real simple module.
\item \label{it:real2} If $M$ is real, then $N$ is a real simple module.
\end{enumerate}
\end{proposition}

\Lemma[{\cite[Lemma 2.27]{KKOP20C}}]  \label{cor: real 1} Let $M,N$
be a real simple module such that $\de(M,N)\le1$. Then $M \hconv N$
is real. \enlemma Remark that only the case $\de(M,N)=1$ is proved
in \cite[Lemma 2.27]{KKOP20C}. However, the assertion is obvious when
$M$ and $N$ commute.

\Lemma \label{lem: real2}
Let $M$ and $N$ be real simple modules.
We assume that
$M\htens N$ commutes with $M$.
Then $M\htens N$ is real simple.
\enlemma

\Proof
We have a commutative diagram (up to constant multiples):
\eq&&\hs{2ex}\ba{c}
\xymatrix@C=8ex{
M\tens N\tens M\tens N\ar@{->>}[d]\ar@{}[rd]|{\textcircled{\scriptsize A}}&M\tens M\tens N\tens N\ar@{->>}[d]\ar[l]_{\rmat{M,N}}\ar@{}[rd]|(.65){\textcircled{\scriptsize B}}
\ar[r]^{\rmat{M,N}}&M\tens N\tens M\tens N\ar@{->>}[d]\\
(M\htens N)\tens M\tens N\ar@{->>}[d]\ar[r]^-\sim_{\rmat{M\hconvs N,M}}
&
M\tens (M\htens N)\tens N\ar[r]_{\rmat{M\hconvs N,N}}\ar@{>->}[ru]&M\tens N\tens (M\htens N)\ar@{->>}[d]\\
(M\htens N)\tens (M\htens N)\ar[rr]_{\rmat{M\hconvs N,M\hconvs N}}
&&
(M\htens N)\tens (M\htens N).
\ar@{}[llu]|(.5){\textcircled{{\scriptsize C}}}
}\ea\label{eq:commd}
\eneq
Here, the commutativity (up to a constant  multiple) of $\textcircled{\scriptsize A}$ follows from the fact that
$M\tens M\tens N$ has a simple head $M\tens(M\hconv N)$ and that
the composition $M\tens M\tens N\To[\rmat{M,N}]
M\tens N\tens M\To(M\hconv N)\tens M$ does not vanish by
 Proposition~\ref{prop:rcomp}~\eqref{item1}.
The commutativity (up to a constant  multiple) of $\textcircled{\scriptsize B}$
follows from the fact that
the composition $(M\hconv N)\tens N\monoto N\tens M\tens N
\epito N\tens(M\hconv N)$ does not vanish by Proposition~\ref{prop:rcomp}~\eqref{item1}
and that $\dim \Hom\bl(M\hconv N)\tens N,N\tens (M\hconv N)\br=1$
by Theorem~\ref{thm: KKKo15 main}\; ~\eqref{runiq}.

The commutativity (up to a constant multiple) of $\textcircled{\scriptsize C}$
follows from the fact that
the composition $(M\hconv N)\tens M\tens N
\to M\tens N\tens (M\hconv N)\to (M\hconv N)\tens (M\hconv N)$
does not vanish by Proposition~\ref{prop:rcomp}\;~\eqref{item1}.

Thus we obtain the commutativity of the diagram \eqref{eq:commd}.

The composition $M\tens M\tens N\tens N\To[\rmat{M,N}]
M\tens N\tens M\tens N\To(M\hconv N)\tens (M\hconv N)$
is an epimorphism since it is the composition of the epimorphisms
$M\tens M\tens N\tens N\epito (M\hconv N)\tens M\tens N$
and $(M\hconv N)\tens M\tens N\epito(M\hconv N)\tens (M\hconv N)$.
It implies that $\rmat{M\hconvs N,M\hconvs N}=\id_{(M\hconvs N)\tens(M\hconvs N)}$
up to a constant multiple.
\QED

\Lemma\label{lem:3}
Let $L_j$ and  $M_j$ be real simple modules $(j=1,2)$.
Assume that
\bnum\item
$L_j\hconv M_j$ commutes with $L_k$ for $j,k=1,2$,
\item $L_1$ and $L_2$ commute.
\ee
Then we have the followings:
\bna
\item \label{it: real} $L_j\hconv M_j$ is real for $j=1,2$.
\item If $\de(\dual L_j,M_2)=0$ for $j=1,2$, then
$$(L_1\tens L_2)\hconv (M_1\htens M_2)\simeq \bl (L_1\tens L_2)\hconv M_1\br
\htens \label{it: com} M_2
\simeq (L_1\hconv M_1)\hconv (L_2\hconv M_2).$$
\item \label{it: com2}
Assume that $\de(\dual L_j,M_k)=0$ for $j,k=1,2$.
Then $M_1$ and $M_2$ commute if and only if $L_1\hconv M_1$ and $L_2\hconv M_2$ commute.
\ee
\enlemma
\Proof
\eqref{it: real} follows from the preceding lemma.

\snoi
\eqref{it: com}
The first isomorphism follows from the fact that
$(L_1\tens L_2,M_1,M_2)$ is normal.

On the other hand,
$(L_2,L_1\hconv M_1,M_2)$ is normal, and hence
$L_2\tens(L_1\hconv M_1)\tens M_2$ has a simple head.
Since we have epimorphisms
\begin{align*}
&L_2\tens(L_1\hconv M_1)\tens M_2
\simeq(L_1\hconv M_1)\tens L_2\tens M_2\epito
(L_1\hconv M_1)\htens(L_2\htens M_2) \ \text{ and } \\
& L_2\tens(L_1\hconv M_1)\tens M_2
\simeq
\bl(L_2\tens L_1)\htens M_1 \br\tens M_2\epito \bl(L_1\tens L_2)\htens M_1\br
 \htens M_2,
\end{align*}
we obtain
$$(L_1\hconv M_1)\htens(L_2\htens M_2)
\simeq \bl(L_1\tens L_2)\htens M_1\br \htens M_2.$$

\mnoi
\eqref{it: com2}  Assume first that $M_1$ and $M_2$ commute.
Then we have
by (ii)
\eqn
(L_1\hconv M_1)\htens(L_2\htens M_2)
&&\simeq
(L_1\tens L_2)\htens (M_1 \htens M_2)\\
&&\simeq
(L_1\tens L_2)\htens (M_2 \htens M_1)
\simeq
(L_2\htens M_2)\htens (L_1\hconv M_1).
\eneqn
Hence  $L_1\hconv M_1$ and $L_2\hconv M_2$ commute.

Conversely, assume that $L_1\hconv M_1$ and $L_2\hconv M_2$ commute.
Then we have
\eqn
(L_1\tens L_2)\htens (M_1 \htens M_2)
&&\simeq
(L_1\hconv M_1)\htens(L_2\htens M_2)\\
&&\simeq
(L_2\htens M_2)\htens (L_1\hconv M_1)
\simeq
(L_1\tens L_2)\htens (M_2\htens M_1).
\eneqn
Hence we have $M_1\htens M_2\simeq M_2\htens M_1$, which implies that
$M_1$ and $M_2$ commute.
\QED

\begin{proposition} \label{prop: three}
Let $M$, $N$ and $L$ be simple $U_q'(\g)$-modules such that $L$ is
real. Assume that \bnum
\item $\de(\D M,L)=0$, $\de(\D L,N)=0$, and
\item $M \tens L \tens N $ has a simple head.
\end{enumerate}
Then we have
$$  \de(L, \hd(M\tens L \tens N)) = \de(L,M \hconv L) + \de(L,L \hconv N).$$
\end{proposition}

\begin{proof} The condition $\de(\D L,N)=0$ implies that
$  (L, M \hconv L,N) $ is normal by Lemma~\ref{prop: L N* normal}.
Thus we have
\begin{align*}
\La(L,\hd(M\tens L \tens N)) &= \La(L,  (M\hconv L) \hconv N)   \\
&= \La(L,  M\hconv L) +\La(L, N)  = \La(L, M\hconv L) +\La(L,
L\hconv N)
\end{align*}
by Lemma~\ref{eq: normal head}.
Similarly, $  (M, L \hconv
N,L) $ is normal and hence we have
\begin{align*}
\La(\hd(M\tens L \tens N),L) &= \La(  M\hconv (L \hconv N) ,L)  \\
&= \La(M, L) +\La(L \hconv N, L)  = \La(M \hconv L, L) +\La(L \hconv
N, L)
\end{align*}
by the same reason. Summing up two above equations, we have
\[  \de(L, \hd(M\tens L \tens N)) = \de(L,M \hconv L) + \de(L,L \hconv N). \qedhere \]
\end{proof}

\Lemma[{\cite[Proposition 3.11, Proposition
4.5]{KKOP19C}}]\label{lem:Rmat} Let $L, M,N$ be simple modules in
$\catg$ and let $S$ be a simple quotient of $M\tens N$. Then,
$m\seteq\La(L,M)+\La(L,N)-\La(L,S)$ belongs to $\Z_{\ge0}$ and the
diagram
$$
\xymatrix@C=10ex{ L_z\tens M\tens
N\ar[r]^{\Rren_{L_z.M}}\ar@{->>}[d]&M\tens L_z\tens
N\ar[r]^{\Rren_{L_z.N}}
& M\tens N\tens L_z\ar@{->>}[d]\\
L_z\tens S\ar[rr]^{a(z)(z-1)^m\Rren_{L_z.S}}&&S\tens L_z}
$$
commutes for some
$a(z)\in\cor(z)$ which has  neither pole nor
zero at $z=1$.
 \enlemma

\begin{definition} \label{def: root module}
A {\em root module} is a real simple module $L$ such that
\begin{align*}
&\de(L, \D^k(L)) = \delta(k=\pm 1) \quad \text{for any $k\in\Z$.}
\end{align*}
\end{definition}

The following proposition is obtained by
the explicit denominator formulas between fundamental
representations (e.g., see \cite[Appendix]{Oh15}, \cite{OT19} and
\cite{Fu19}).
\Prop Every fundamental representation is a root module.
\enprop

\subsection{Simply-laced finite root system associated with $\Ca_\g^0$} \label{subsec: E(M)}

In \cite{KKOP20A},  we associate to the category $\catCO$ a
simply-laced finite type root system
in a canonical way:
for a simple  module $M \in \Ca_\g$, set $ \rE(M) \in
\Hom_\Set(\sig, \Z)$ (see~\eqref{eq: sig}) by
\begin{align*}
\rE(M)(i,a)= \Li(M, V(\varpi_i)_a) \qquad \text{ for } (i,a) \in
\sig.
\end{align*}
Let
\begin{equation*}
\begin{aligned}
\rW_0 \seteq \{  \rE(M) \mid M \text{ is simple in } \Ca_\g^0 \}
\qtq\Delta_0 \seteq \st{  \rE(V( \varpi_i)_a) \mid (i,a) \in \sigZ }
\subset \rW_0.
\end{aligned}
\end{equation*}

Then, $\rW_0$ is a $\Z$-submodule of $\Hom_\Set(\sig, \Z)$ and
endowed with a symmetric bilinear form $(\ ,\ )$ which satisfies
$$\text{$(\rE(M),\rE(N))=-\Li(M,N)$ for any simple modules $M$, $N$
in $\Ca_\g^0$.}$$

\begin{theorem}[\cite{KKOP20A}]
The pair $(\R\tens \rW_0, \Delta_0 )$ is a root system, the bilinear
form $(\ ,\ )$ is positive definite and invariant by the Weyl group
action, and we have $\Delta_0=\st{\al\in\rW_0\mid
(\al,\al)=2}$.
\end{theorem}

\vskip 2em
\section{Review on Quiver Hecke algebras} \label{sec: Quiver Hecke algebra}
In this section, we briefly recall the definition of quiver Hecke algebras $R$
and the categorification of the quantum unipotent coordinate rings
by $R$.
Then we will review T-system among the determinantial modules.

\subsection{Symmetric quiver Hecke algebras}  In this subsection, we
recall basic notions of symmetric quiver Hecke algebras associated to
a finite simple Lie algebra $\fing$ of simply-laced type. We denote
by $\Pi_\fing=\{ \al_\im \mid \im \in I_\fing \}$  the set of simple
roots of $\fing$ and by $( \ , \ )$ the symmetric bilinear form on the
root lattice $\rl$ of $\fing$ invariant by the Weyl group action and normalized
by $(\al_\im,\al_\im)=2$.
 Set $\rl^+=\dsum_{\im\in\Isf}\Z_{\ge0}\al_\im$.

Take a family of polynomials $(Q_{\im\jm})_{\im,\jm \in I_\fing}$ in
$\ko[u,v]$ which are of the form
$$   Q_{\im\jm} (u,v) = \pm\delta(\im \ne \jm)(u-v)^{-(\al_\im,\al_\jm)} $$
such that $Q_{\im\jm}(u,v)=Q_{\jm\im}(v,u)$.

For each $\be \in \rl^+$ with $|\be|=n$, we set $I^\be_\fing = \{
\nu=( \nu_1,\ldots,\nu_n) \in I^n_\fing \ | \ \sum_{k=1}^n
\al_{\nu_k} = \be \}$. Note that the symmetric group
$\mathfrak{S}_n=\lan \mathfrak{s}_k \ | \ 1 \le k \le n-1 \ran$ acts
on $\Iff^\be$ by place permutations.

The \emph{symmetric quiver Hecke algebra} (also called
\emph{Khovanov-Lauda-Rouquier algebra}) $R(\be)$ at $\be \in \rl^+$
associated to $\fing$ and $(Q_{\im\jm})_{\im,\jm\in  \Iff }$, is the
algebra over $\ko$ generated by the elements $\{ e(\nu) \}_{\nu \in
\Iff^\be}$, $\{ x_k \}_{1\le k\le n}$ and $\{ \tau_m \}_{1\le m\le
n-1}$ satisfying the following defining relations:
\begin{align*}
&e(\nu)e(\nu')=\delta_{\nu,\nu'}e(\nu), \ \ \sum_{\nu \in \Iff^\be} e(\nu)=1, \ \ x_kx_m= x_mx_k, \ \ x_ke(\nu)=e(\nu)x_k, \allowdisplaybreaks\\
& \tau_m e(\nu) = e(\mathfrak{s}_m(\nu))\tau_m,  \ \  \tau_k\tau_m = \tau_m\tau_k \text{ if } |k-m|>1, \ \  \tau_k^2e(\nu)  = Q_{\nu_k,\nu_{k+1}}(x_k,x_{k+1})e(\nu), \allowdisplaybreaks\\
&(\tau_kx_m-x_{\mathfrak{s}_k(m)}\tau_k )e(\nu) =
\begin{cases}
-e(\nu) & \text{ if } m=k, \nu_k=\nu_{k+1}, \\
e(\nu) & \text{ if } m=k+1, \nu_k=\nu_{k+1}, \\
0 & \text{ otherwise},
\end{cases} \allowdisplaybreaks\\
&(\tau_{k+1}\tau_{k}\tau_{k+1}-\tau_{k}\tau_{k+1}\tau_{k})e(\nu) =
\delta_{\nu_k,\nu_{k+2}}
\dfrac{Q_{\nu_k,\nu_{k+1}}(x_k,x_{k+1})-Q_{\nu_k,\nu_{k+1}}(x_{k+2},x_{k+1})
}{x_k-x_{k+2}}e(\nu).
\end{align*}
The algebra $R(\beta)$ is $\Z$-graded with
$$\deg e(\nu)=0, \quad \deg x_ke(\nu)=2, \quad\text{and}\quad \deg \tau_m e(\nu)=-(\alpha_{\nu_m},\alpha_{\nu_{m+1}}).$$

Let us denote by $R(\be)\gMod$ the category of graded
$R(\beta)$-modules with homomorphisms of degree $0$, by
$R(\beta)\gmod$ the full subcategory of  $R(\be)\gMod$ consisting of
finite-dimensional graded $R(\beta)$-modules. For an
$R(\beta)$-module $M$, we set $\wt(M)\seteq -\beta \in \rl^-$. For
the sake of simplicity, we  say that $M$ is an $R$-module instead of
saying that $M$ is a graded $R(\beta)$-module. Let us denote by $q$
the grading shift functor: for $M\in R\gmod \seteq \soplus_{\be \in
\rl^+}R(\be)\gmod$. We have $(qM)_n=M_{n-1}$ by definition.

For $\im \in \Isf$,  $L(\im)$ denotes the $1$-dimensional
simple graded
$R(\al_\im)$-modules $\ko u(\im)$ with the action
$$ \quad x_1u(\im)=0.$$

We also sometimes ignore grading shifts if there is no  risk of
confusion. Hence, for $R$-modules $M$ and $N$,  we sometimes say
that $f\col M\to N$ is a homomorphism if $f\col q^aM\to N$ is a
morphism in $R\gmod$ for some $a\in\Z$. If we want to emphasize that
$f\col q^aM\to N$ is a morphism in $R\gmod$, we say so.  We set
$$\HOM_{R(\beta)}(M,N)\seteq \soplus_{a\in \Z} \Hom_{R(\beta)\gmod}(q^aM,N).$$

For $\beta, \gamma \in \rl^+$,  set $e(\beta,\gamma)\seteq  \sum
_{\nu \in I_\fing^{\beta},\; \nu' \in I_\fing^{\gamma}}
 e(\nu * \nu')  \in  R(\beta+\gamma)$, where $\nu*\nu'$ denotes the concatenation of $\nu$ and $\nu'$.
Then there is a $\cor$-algebra homomorphism $R( \beta)\tens R(
\gamma  )\to e(\beta,\gamma)R( \beta+\gamma)e(\beta,\gamma)$. The
\emph{convolution product} $- \conv - \col R(\beta)\gmod \times
R(\gamma)\gmod \to R(\beta+\gamma)\gmod$
 is a bifunctor given by
$$M\conv N\seteq R(\beta + \gamma) e(\beta,\gamma)
\tens_{R(\beta )\otimes R( \gamma)}(M\otimes N).$$
Set $R\gmod\seteq \soplus_{\beta \in \rl^+} R(\beta)\gmod$.
 Then the category $R\gmod$ is a monoidal category with the tensor product $\conv$ and the unit object $\one\seteq \cor \in R(0) \gmod$.

\smallskip

For $1\le a<|\beta|$,  we define the {\em intertwiner} $\vphi_a\in
R( \beta)$ by
\begin{align*}
  \vphi_a e(\nu)=
\begin{cases}
  \bl\tau_ax_a-x_{a}\tau_a\br e(\nu)
  & \text{if $\nu_a=\nu_{a+1}$,} \\[2ex]
\tau_ae(\nu)& \text{otherwise.}
\end{cases}
\end{align*}

For $m,n\in\Z_{\ge0}$, let us denote by $w[{m,n}]$ the element of
$\sym_{m+n}$  defined by
\begin{align*}
w[{m,n}](k)=\begin{cases}k+n&\text{if $1\le k\le m$,}\\
k-m&\text{if $m<k\le m+n$,}\end{cases}
\end{align*}
and set $\vphi_{[m,n]}\seteq\vphi_{\im_1}\cdots\vphi_{\im_t}$, where
$s_{\im_1}\cdots s_{\im_t}$ is a reduced expression of $w[m,n]$.
 It does not depend on the choice of a reduced expression.

Let $M$ be an $R(\beta)$-module and $N$ an $R(\gamma)$-module. Then
the map $M\tens N\to N\conv M$ given by $u\tens v\longmapsto
\vphi_{[|\gamma|,|\beta|]}(v\tens u)$ induces  an $R( \beta
+\gamma)$-module  homomorphism
\begin{align*}
R_{M,N}\colon M\conv N\To N\conv M.
\end{align*}

Let $z$ be an indeterminate homogeneous of degree $2$, and let
$\psi_z$ be the graded algebra homomorphism
\begin{align*}
\psi_z\col R( \beta )\to \cor[z]\tens R( \beta )
\end{align*}
given by
$$\psi_z(x_k)=x_k+z,\quad\psi_z(\tau_k)=\tau_k, \quad\psi_z(e(\nu))=e(\nu).$$

For an $R( \beta )$-module $M$, we denote by $M_z$ the
$\bigl(\ko[z]\tens R( \beta )\bigr)$-module $\ko[z]\tens M$ with the
action of $R( \beta )$ twisted by $\psi_z$.

For a non-zero $R(\beta)$-module $M$ and a non-zero
$R(\gamma)$-module $N$, let $s$ be the largest non-negative integer
such that the image of $R_{M_z,N}$ is contained in $z^s(N\conv
M_z)$. We denote by $\Rren_{M_z,N}$ the $R(\beta)$-module
homomorphism $z^{-s}R_{M_z,N}$ and call it the \emph{renormalized
R-matrix}. We denote  by $\rmat{M,N}\seteq\Rren_{M_z,N}\bigm\vert_{z=0}
\col M\conv N \to N \conv M$ and call $\rmat{M,N}$ the
\emph{R-matrix}.
 By the definition, $\rmat{M,N}$ never vanishes.

We denote by $\La(M,N)$ the homogeneous degree of $\rmat{M,N}$ and
set
\begin{align*}
\tLa(M,N) \seteq \dfrac{1}{2}\bigl(\La(M,N) +(\wt(M), \wt(N)) \bigr)
,\quad \de(M,N) \seteq \dfrac{1}{2}\bigl(\La(M,N)+\La(N,M) \bigr).
\end{align*}
The invariants $\La$ and $\de$ enjoy the similar properties of $\La$
and $\de$ for $\Ca_\g$. For more details,
we refer to \cite{KKKO18,KKOP19C,KKOP20A}.

Let $\beta\in\rl^+$. For $\im \in I_\fing$ , $1 \le a \le |\beta|$
and $M\in R(\be)\gmod$, we set
\begin{align*}
E_\im M \seteq e_1(\im)M \qtq
E^*_\im M \seteq e_{|\be|}(\im)M
\qt{where $e_a(\im)\seteq\sum_{\nu \in I^\beta_\fing,\nu_a=\im}e(\nu)\in
R(\beta)$.}
\end{align*}
Then, $E_\im$ and $ E^*_\im$ are functors from $R(\beta) \gmod$ to
$R(\beta-\al_\im) \gmod$. For a non-zero module $M\in R(\be)\gmod$,
we also set \eqn
\eps_\im(M)&&=\max\set{n\in\Z_{\ge0}}{E_\im^nM\not=0} \in \Z_{\ge 0} \qtq\\
\eps^*_\im(M)&&=\max\set{n\in\Z_{\ge0}}{{E^*_\im}^nM\not=0} \in
\Z_{\ge 0}. \eneqn

\subsection{Quantum unipotent coordinate ring and T-systems}
Let $U_q(\fing)$ be the quantum group associated to
a finite simply-laced Cartan datum $(\cmf,\wl,\Pi,\wl^\vee,\Pi^\vee)$ of $\fing$,
which is the  associative algebra over $\mathbb Q(q)$ generated by
$e_\im,f_\im$ $(\im \in \Iff)$ and $q^h$ ($h\in\wl^\vee$).

Note that $U_q(\fing)$ admits a weight-space decomposition as
follows:
$$U_q(\fing)= \bigoplus_{\beta \in \rl} U_q(\fing)_{\beta}, \text{ where } U_q(\fing)_{\beta}\seteq\set{ x \in U_q(\fing)}{\text{$q^{h}x q^{-h} =q^{\langle h, \beta \rangle}x$ for any $h \in \wl^\vee$}}.$$

We denote by $U_q^{+}(\fing)$ the $\Q(q)$-subalgebra of $U_q(\fing)$
generated by $\set{e_\im}{\im \in \Iff}$. Let
$U_q^{+}(\fing)_{\Z[q^{\pm1}]}$ be the $\Z[q^{\pm1}]$-subalgebra of
$U_q^{+}(\fing)$ generated by $e_\im^{(n)}\seteq e_\im^n/[n]!$
($\im\in \Iff$, $n\in\Z_{>0}$).

Let $\Delta_\finn$ be the algebra homomorphism $U_q^+(\fing) \to
U_q^+(\fing) \tens U_q^+(\fing)$ given by  $ \Delta_\finn(e_\im) = e_\im
\tens 1 + 1 \tens e_\im$,
where the algebra structure on $U_q^+(\fing)
\tens U_q^+(\fing)$ is defined by
$$(x_1 \tens x_2) \cdot (y_1 \tens y_2) = q^{-(\wt(x_2),\wt(y_1))}(x_1y_1 \tens x_2y_2).$$

Set
$$ A_q(\finn) = \soplus_{\beta \in  \rl^-} A_q(\finn)_\beta \quad \text{ where } A_q(\finn)_\beta \seteq \Hom_{\Q(q)}(U^+_q(\fing)_{-\beta}, \Q(q)).$$
Then  $A_q(\finn)$ is an algebra with the multiplication given by
$(\psi \cdot \theta)(x)= \theta(x_{(1)})\psi(x_{(2)})$, when
$\Delta_\finn(x)=x_{(1)} \tens x_{(2)}$ in  Sweedler's notation.

Let us denote by $\An_{\Z[q^{\pm 1}]}$ the $\Z[q^{\pm 1}]$-submodule
of $\An$
 consisting of $ \psi \in \An$ such that
$ \psi\bl U_q^{+}(\fing)_{\Z[q^{\pm1}]}\br\subset\Z[q^{\pm1}]$. Then
it is a $\Z[q^{\pm 1}]$-subalgebra of $\An$.
 Moreover, the {\em upper global basis}  $\B^\up$ is
a $\Z[q^{\pm1}]$-module basis of $\An_{\Z[q^{\pm 1}]}$. For the
detail of the upper global basis, see \cite[Section 1.3]{KKKO18}.
For $w \in \W_\fing$, let $\Anw_{\Z[q^{\pm 1}]}$ be the $\Z[q^{\pm
1}]$-submodule of $\An_{\Z[q^{\pm 1}]}$
 consisting of  the elements $\psi \in \An_{\Z[q^{\pm 1}]}$ such that
$e_{\im_1}\cdots e_{\im_n} \psi =0$  for any $\beta\in
\bigl(\rl^+\cap w \rl^+\bigr)
 \setminus \{0\}$ and
any sequence $(\im_1,\ldots,\im_n) \in \Iff^\beta$. Then
$\Anw_{\Z[q^{\pm 1}]}$  is a $\Z[q^{\pm 1}]$-subalgebra of
$\An_{\Z[q^{\pm 1}]}$ (\cite[Theorem 2.20]{KKOP18}). We call this
algebra the \emph{quantum unipotent coordinate ring associated with
$w$.} The set $\mathbf B^\up(w) \seteq\mathbf B^\up \cap \Anw$ forms
a $\Z[q^{\pm 1}]$-basis of $\Anw_{\Z[q^{\pm 1}]}$ (\cite{Kimu12}).

For a dominant weight $\la \in \wl^+$, let $V(\la)$ be the
irreducible highest weight $U_q(\fing)$-module with highest weight
vector $u_\la$ of weight $\la$. Let $( \ , \ )_\la$ be the
non-degenerate symmetric bilinear form on $V(\la)$ such that
$(u_\la,u_\la)_\la=1$ and $(xu,v)_\la=(u,\varphi(x)v)_\la$ for $u,v
\in V(\la)$ and $x\in U_q(\fing)$, where $\varphi$ is  the algebra
antiautomorphism on $U_q(\fing)$  defined by $\varphi(e_i)=f_i$,
$\varphi(f_i)=e_i$  and $\varphi(q^h)=q^h$. For each $\mu, \zeta \in
\W_\fing \la$, the \emph{unipotent quantum minor} $D(\mu, \zeta)$ is
an element in $A_q(\finn)$ given by
 $D(\mu,\zeta)(x)=(x u_\mu, u_\zeta)_\la$
for $x\in U_q^+(\fing)$, where $u_\mu$ and $u_{\zeta}$ are the
extremal weight vectors in $V(\la)$ of weight $\mu$ and $\zeta$,
respectively. Then we have $D(\mu,\zeta)\in \mathbf B^\up
\sqcup\st{0}$.

For a reduced expression $\tw=s_{\im_1}s_{\im_2}\cdots s_{\im_l}$ of
$w \in \W_\fing$, define
$\tw_{\le k}\seteq s_{\im_1}\cdots s_{\im_k}$ and
$ \la_k  \seteq \tw_{\le k}\La_{\im_k}$ for $1 \le k \le l$.
Note that $\la_{k^-} =\tw_{\le k-1}\La_{\im_k}$ if $k^- >0$. Here
$$k^-\seteq \max \Bigl(\st{s\in[1,k -1 ]\mid
\im_s=\im_k}\sqcup\st{0}\Bigr).
$$
For $0 \le t \le s \le l$, we set
$$
D_{\tw}[t,s] \seteq \begin{cases}
D( \tw_{\le s}\La_{\im_s},\ \tw_{\le t-1}\La_{\im_t}) & \text{ if } 1 \le  t \le s\le l \text{ and } \im_s=\im_t, \\
D( \tw_{\le s}\La_{\im_s},\ \La_{i_s}) & \text{ if } 0 = t < s \le l, \\
\mathbf{1}  & \text{ if } 0 = t =s.
\end{cases}
$$
Then $D_{\tw}[t,s]$ belongs to $\mathbf B^\up(w)$. The set
$\set{  D_{\tw}[k] \seteq D_{\tw}[k, k]}  {1\le k\le l}$  generates $\Anw_{\Z[q^\pm 1]}$
as a $\Z[q^{\pm1}]$-algebra.

We omit the subscript $_\tw$ if there is no confusion about the
choice of $\tw$.

\begin{remark}
In \cite{KKKO18}, we use $D(s, t^+ )$ instead of $D[t,s]$,
where $k^+\seteq\min\Bigl( \{t \in [k+1,l]\  \mid   \im_t=\im_k \} \sqcup \{ l+1 \} \Bigr)$.
We hope that there is no confusion by the differences of the notations.

\end{remark}

\begin{theorem} [{\cite{KL09,R08,R11,VV09}}] \label{thm: KLR iso}
There exists a $\Z[q^{\pm 1}]$-algebra isomorphism:
$$ {\rm ch}_q \colon K(R\gmod) \isoto \An_{\Z[q^{\pm 1}]}.$$
Furthermore, under the isomorphism ${\rm ch}_q$, the upper global
basis $\mathbf B^\up$ of $\An_{\Z[q^{\pm 1}]}$ corresponds to the
set of the isomorphism classes of self-dual simple $R$-modules.
\end{theorem}

Since $D(\mu,\zeta)$ is a member of the upper global basis, there exists a unique real simple module $\DC(\mu,\zeta)$ in $R(\zeta-\mu)\gmod$
such that ${\rm ch}_q(\DC(\mu,\zeta))=D(\mu,\zeta)$. We call it \emph{the determinantial module}. For a reduced expression
$\tw=s_{\im_1}s_{\im_2}\cdots s_{\im_l}$ of $w \in \W_\fing$ and $0\le t\le s\le l$ such that $t=0$ or $\im_s=\im_t$, let
$\DC_{\tw}[t,s]$ be a simple module such that ${\rm
ch}_q(\DC_{\tw}[t,s]) = D_{\tw}[t,s]$. In particular,  we call
$\DC_{\tw}[k]\seteq\DC_{\tw}[k,k]$
 \emph{the cuspidal module} for $1 \le k \le l$.

\begin{proposition} [{\cite[Proposition 10.2.5]{KKKO18}, \cite[Proposition 4.6]{KKOP18}}]
For a reduced expression $\tw=s_{\im_1}s_{\im_2}\cdots s_{\im_l}$ of
$w \in \W_\fing$, we have the followings$\colon$ \bnum
\item For $1 \le a  < b \le l$ with $\im_a=\im_b$, there exists an exact sequence in $R\gmod:$
\begin{align} \label{eq: T-system in QH}
0  \to  \conv_{\substack{\jm\in \Iff; \; d(\im_a,\jm)=1} }  \DC[a(\jm)^+,b(\jm)^-] &\to  \DC[a^+ ,b]
\conv  \DC[a,b^-] \to  \DC[a,b]   \conv  \DC[a^+ ,b^-]   \to   0.
\end{align}
Here we understand $\DC[t,s] \simeq \one$ for $t >s$.
For the notations $k(\jm)^{\pm}$, see \eqref{eq: nota +,-} below.
\item For $1 \le a \le b  <  c \le l$ with $\im_a=\im_b=\im_c$, we have
\begin{align*}
\DC[b^+,c]   \hconv \DC[a ,b]  \simeq \DC[a,c].
\end{align*}
\end{enumerate}
\end{proposition}

We call an exact sequence of the form in~\eqref{eq: T-system in QH}
a {\em T-system}.

\bigskip
For  an $R(\beta)$-module $M$,  define
\begin{align*}
&\mathrm W (M)\seteq \set{\gamma \in  \rl^+\cap (\beta-\rl^+)}{e(\gamma, \beta-\gamma) M \neq 0}, \\
&\mathrm W^* (M)\seteq \set{\gamma \in  \rl^+\cap
(\beta-\rl^+)}{e(\beta-\gamma,\gamma) M \neq 0}.
\end{align*}
For $w \in \W_\fing$, we denote by $R_w\gmod$ the full subcategory
of $R\gmod$ whose objects $M$ satisfy $\mathrm W(M) \subset
\rl^+\cap w \rl^-$. Here $R$ denote the quiver Hecke algebra
associated to $\fing$.

\begin{theorem}[\cite{KKOP18}] \label{thm: KKOP18main} \hfill
\bna
\item $R_w\gmod$ is  the  smallest subcategory of $R\gmod$ which is stable under taking subquotients, extensions, convolutions, grading shifts, and  contains  the determinantial modules
$\{  \DC_{\tw}[k] \mid 1\le k\le l \}$, where
$\tw=s_{\im_1}s_{\im_2}\cdots s_{\im_l}$ is a reduced expression  of
$w$.
\item There exists a $\Z[q^{\pm1}]$-algebra isomorphism between $K(R_w\gmod)$ and  $A_{\Z[q^{\pm1}]}(\finn(w))$. Here, $A_{\Z[q^{\pm1}]}(\finn(w))$ is the
$\Z[q^{\pm1}]$-subalgebra of $ \An_{\Z[q^{\pm 1}]}$ generated by
$D[t,s]$ with $1\le t\le s\le l$.
\item The set of self-dual  simple modules in $R_w\gmod$  coincides with the upper global basis of $A_{\Z[q^{\pm1}]}(\finn(w))$ under the isomorphism ${\rm ch}_q$.
\end{enumerate}
In particular, for any $\La \in \wl^+$, the determinantial module
$\DC(w\La,\La)$ strongly commutes with all the simple modules in $R_w\gmod$.
\end{theorem}

\vskip 2em
\section{Quantum affine \SW duality and T-systems} \label{sec: SW dailty and T-system}
In this section, we briefly recall the quantum affine
\SW duality functor from the category $R\gmod$ of finite-dimensional
modules to $\Ca_\g$, constructed by Kang-Kashiwara-Kim
in \cite{KKK18A} (cf.\ \cite{KKOP20C}), which preserves several
invariants of those categories.
We also review the affine cuspidal
modules, \pbws and reflections, defined in
\cite{KKOP20C}. Then we will construct affine determinantial
modules, which are real simple modules associated to  $i$-boxes and
\pbws. In the last part, we will investigate the commuting condition
and  T-systems among affine determinantial modules, which can
be understood as a vast generalization of the T-systems in $R\gmod$ and
among KR-modules.

\subsection{Quantum affine \SW duality functor} \label{subsec: QAWS}

Let $\fing$ be a simple Lie algebra with a  simply-laced Cartan
matrix  $\cmf=(c_{\im,\jm})_{\im,\jm\in\Iff}$. Let $\Dd \seteq \{
\Rt_\im \}_{\im\in \Iff} \subset \Ca_\g$ be a family of real simple
modules of $\Ca_\g$. The family $\Dd $ is called a \emph{strong
duality datum} associated with $\fing$ if it satisfies the following
conditions (\cite[Definition 4.7]{KKOP20C}): \ben
\item
$\mathsf{L}_\im$ is a root module for each $\im \in \Iff$,
\item for any $\im,\jm\in \Iff$ with $\im\ne \jm$, $\de(\mathsf{L}_\im, \D^k \mathsf{L}_\jm) = -\delta(k=0)c_{\im,\jm}$.
\end{enumerate}

Let us denote by $\Rc$ the symmetric quiver Hecke algebra associated
with $\cmf$.

\begin{theorem}[{\cite{KKK18A,KKOP20C}}] \label{thm:gQASW duality} For a given strong duality datum $\Dd$ associated with $\fing$, there exists a functor
$$\widehat{\F}_{\Dd} \col  R_\cmf\gMod\rightarrow  U_q'(\g)\Md,$$
where $U_q'(\g)\Md$ denotes the category of integrable $U_q'(\g)$-modules.
Moreover, $\widehat{\F}_{\Dd}$
satisfies the following properties$\colon$ \bna
\item $\widehat{\F}_\Dd(L(\im))\simeq \mathsf{L}_\im$.
\item The functor $\widehat{\F}_\Dd$ induces an exact monoidal functor on $R\gmod$ which preserves simplicity$\colon$
\begin{align}\label{eq: F fd}
\F_\Dd \col \Rc\gmod \rightarrow \Ca_\g.
\end{align}
Namely, $\F_\Dd$ sends finite-dimensional graded $\Rc$-modules to
modules in $\Ca_\g$ and,   for any $M_1, M_2 \in  \Rc\gmod$, we have
isomorphisms
$$\F_\Dd(\Rc(0)) \simeq \ko, \quad \F_\Dd(M_1 \circ M_2) \simeq \Fd(M_1) \tens \Fd(M_2),$$
and $\F_\Dd$ sends simple modules to simple modules.
\item For $M \in \Rc\gmod$, we have
 $$\widehat{\F}_\Dd(M_z) \simeq \F_\Dd(M)_{\e^z} $$
 \rmo see \cite{KKOP20C} for more details\/\rmf.
\end{enumerate}
\end{theorem}

We call the functor $\Fd$ the \emph{quantum affine \SW duality functor}.

\begin{definition}[{\cite[\S 5.2]{KKOP20C}}]\label{def: CD}
The category $\Ca_\Dd$ is defined to be  the smallest full
subcategory of $\Ca_\g$ such that \ben
\item   it contains $\mathsf{L}_\im$ for any $\im\in I_\fing$
and the trivial module $\one$,
\item  it is stable by taking subquotients, extensions, and tensor products.
\ee
\end{definition}
Hence $\Ca_\Dd$ contains $\F_\Dd( M )$ for any module
$M\in R_{\cmf}\gmod$.

\smallskip
The Grothendieck ring $K(R_{\cmf}\gmod)$ has a
$\Z[q^{\pm1}]$-algebra structure. We set
$$  K(R_{\cmf}\gmod)|_{q=1} \seteq K(R_{\cmf}\gmod)/(q-1)K(\Rc\gmod).$$

\begin{theorem}[\cite{KKOP20C}]
\label{thm: isomorphisms} Let $\Dd$ be a strong duality datum. Then
there exist a $\Z$-algebra isomorphism
\begin{align*}
{\rm ch}\col K(R_{\cmf}\gmod)|_{q=1} \isoto
K(\Ca_\Dd),
\end{align*}
and a $\C$-algebra isomorphism
$$  \C[N] \simeq  \C\tens\bl K(R_{\cmf}\gmod)|_{q=1}\br,$$
where $K(\Ca_\Dd)$ denotes the Grothendieck ring of $\Ca_\Qd$ and
$\C[N]$ denotes the coordinate ring of the unipotent group $N$
associated with the nilpotent Lie subalgebra $\finn$ of $\fing$.
\end{theorem}

\subsection{ $\La$, $\tLa$ and $\de$ under the duality functors}
Keep the notations $\fing$ and
$\cmf=(c_{\im,\jm})_{\im,\jm\in\Iff}$. Let $\Dd$ be a strong duality
datum associated with $\fing$. For each $k \in \Z$, we define a
subcategory $\D^k \Ca_\Dd$ as follows: $\D^k \Ca_\Dd$ is the
smallest full subcategory of $\Ca^0_\g$ such that (a) it contains
$\D^k(\F_\Dd(L))$ for all simple modules $L$ in $R_\cmf\gmod$, and
(b) it is stable by taking subquotients, extensions and tensor
products. By~\cite[Proposition 5.4]{KKOP20C}, there is no
non-trivial module which is contained in distinct $\D^k \Ca_\Dd$ and
$\D^l \Ca_\Dd$ simultaneously.
Note that we have
$$  \de(M,N)=0   $$
for a simple module $M$ in $\D^k \Ca_\Dd$ and
a simple module $N$ in $\D^l \Ca_\Dd$
if $|k-l|>1$.

We recall one of the main results of \cite{KKOP20C}, which tells
that the invariants $\La,\de$ are preserved under the functor
$\F_\Dd$:

\begin{proposition}[{\cite[Theorem 4.12]{KKOP20C}}]\label{prop: de ve}\hfill
\bnum
\item For  simple modules $M,N \in R_\cmf\gmod$, we have
\eqn \Lambda(M,N)&& = \La\bl\F_\Dd(M),\F_\Dd(N)\br,\\
\tLa(M,N)&&  =  \de\bl\D \F_\Dd(M),\F_\Dd(N)\br. \eneqn
\item For a simple module $M$ in $R_\cmf\gmod$ and $\im \in \Iff$, we have
$$\de\bl\D \F_\Dd(L(\im)),\F_\Dd(M)\br  = \ve_\im(M) \qtq \de\bl\D^{-1} \F_\Dd(L(\im)),\F_\Dd(N)\br  =  \ve_{\im}^*(N).$$
\end{enumerate}
\end{proposition}

\subsection{Affine cuspidal modules,  \pbw and reflections}
\label{subsec:pbw}
Let $\fing$ and $\cmf=(c_{\im,\jm})_{\im,\jm\in\Iff}$ be as in the
last subsection. Let $\Dd \seteq \st{ \Rt_\im }_{\im\in \Iff} $ be a
strong duality datum associated to $\fing$, and let $\ell$ be the length of the longest
element $w_0$ of the Weyl group of $\fing$, and  let
$\tw_0=s_{\im_1}\ldots s_{\im_\ell}$ be a reduced expression of
$w_0$. We extend the definition of $\seq{\im_k}_{1\le k\le \ell}$ to
$\seq{\im_k}_{k\in\Z}$ by
$$\im_{k+\ell} = (\im_k)^* \qt{for any $k\in \Z$.}$$
We set $\hhw=\seq{\im_k}_{k \in \Z}$, and we call $(\Dd,\hhw)$ a {\em
\pbw.}

Note that $s_{\im_k}s_{\im_{k+1}} \cdots s_{\im_{k+\ell-1}}$ is a
reduced expression of $w_0$ for any $k \in \Z$.

\Def\label{def:pbw}
 We define a sequence of simple $U_q'(\g)$-modules $\seq{ \bS{k}^{\Dd,\hhw} }_{ k\in \Z }$ in $\Ca_\g $ as follows:
\ben
\item For $1 \le k \le \ell$, $\bS{k}^{\Dd,\hhw} \seteq \F_\Dd( \DC_{\tw_0}[k])$ and, for all $k\in\Z$, it is defined by
\item $\bS{k+\ell}^{\Dd,\hhw} = \D( \bS{k}^{\Dd,\hhw})$ for any $k\in \Z$.
\ee The modules $\bS{k}^{\Dd,\hhw} $  are called the \emph{affine
cuspidal modules} corresponding to the \pbw $(\Dd,\hhw)$.  For
simplicity of notation, we drop the superscript ${}^{\Dd,\hhw}$ if
there is no risk of confusion.
\edf

{\em Throughout this section, we fix a \pbw $(\Dd,\hhw)$} and
frequently use the notations:
\begin{equation} \label{eq: nota +,-}
\begin{aligned}
&s^+\seteq\min\{t \mid s<t ,\; \im_t=\im_s \}, \ &&s^-\seteq\max\{t
\mid
t<s, \; \im_t=\im_s \},\\
&s(\jm)^+\seteq\min\{t \mid s \le  t,\; \im_t=\jm \}, \ &&
s(\jm)^-\seteq\max\{t \mid t \le s,\; \im_t=\jm \}
\end{aligned}
\end{equation}
for $s\in\Z$ and $\jm\in I_\fing$.

\begin{proposition} [{\cite[Proposition 5.7]{KKOP20C}}] \label{prop: lemma of KKOP20C} \hfill
\bnum
\item For every $k \in \Z$, $\bS{k}$ is a root module.\label{it:root}
\item \label{it: far com}For any $a,b \in \Z$ such that $|a-b|\ge \N$, we have  $ \de(\bS{a},\bS{b})=\delta(|a-b|=\N)$.
\item \label{it: unmixed} For any $a<b$, the ordered pair $(\bS{b},\bS{a})$ is strongly unmixed.
\end{enumerate}
\end{proposition}

For any $\jm \in I_\fing$, we set
\begin{align} \label{Def: refl}
 \Refl_\jm (\Dd) \seteq  \{ \Refl_\jm (\Rt_\im)  \}_{\im\in \Iff}\qtq
 \Refl_\jm^{-1} (\Dd) \seteq  \{ \Refl^{-1}_\jm (\Rt_\im)  \}_{\im\in \Iff},
\end{align}
where
$$
\Refl_\jm (\Rt_\im) \seteq
\begin{cases}
 \D \Rt_\im  & \text{ if } \im=\jm, \\
\Rt_\jm \hconv \Rt_\im &  \text{ if }  c_{\im,\jm}=-1, \\
\Rt_\im &  \text{ if } c_{\im,\jm}=0,
\end{cases}
\qtq \Refl^{-1}_\jm (\Rt_\im) \seteq
\begin{cases}
 \D^{-1} \Rt_\im  & \text{ if } \im=\jm, \\
\Rt_\im \hconv \Rt_\jm &  \text{ if }  c_{\im,\jm}=-1, \\
\Rt_\im &  \text{ if } c_{\im,\jm}=0.
\end{cases}
$$
It is easy to see that $\Refl_\jm \circ \Refl^{-1}_\jm (\Dd) = \Dd$
and $ \Refl^{-1}_\jm \circ \Refl_\jm (\Dd) = \Dd $ for any $\jm\in
\Iff$.

\Prop[{\cite[Proposition 5.8, Proposition
5.9]{KKOP20C}}]\label{prop:refl} Let
$\bl\Dd,\hhw=\seq{\im_k}_{k\in\Z} \br$ be a \pbw and $\jm\in\Iff$. \bnum
\item
The data $\Refl_\jm (\Dd)$ and $\Refl^{-1}_\jm (\Dd)$ are also
strong duality data associated with $\cmf$.
\item Set $\Refl^{\pm1}(\hhw)= \seq{\im_k'}_{k\in\Z} $ with $\im'_k=\im_{k\pm1}$.
Then $\Refl(\Dd,\hhw)\seteq\bl\Refl_{\im_1}(\Dd),\Refl(\hhw)\br$ and
$\Refl^{-1}(\Dd,\hhw)\seteq\bl\Refl^{-1}_{\im_1}(\Dd),\Refl^{-1}(\hhw)\br$
are \pbws.

\item
We have
$$\cusp_k^{\Refl^{\pm1}(\Dd,\hhw)}=\cusp^{\Dd,\hhw}_{k\pm1}\qt{for any $k\in\Z$.}$$
\ee \enprop

For each interval $[a,b]$, we denote by $\Ca^{[a,b],\Dd,\hhw}_{\g}$
the smallest full subcategory of $\Ca_\g$  satisfying the following
conditions (see \cite[\S 6.3]{KKOP20C})$\colon$ \ben
\item it is stable under taking subquotients, extensions, tensor products and
\item  it contains $\bS{k}$ for all $a \le k \le b$ and the trivial module $\mathbf{1}$.
\end{enumerate}

Hence we have
$$\Ca^{[a,b],\,\Refl(\Dd,\hhw)}_{\g}=\Ca^{[a+1,b+1],\,\Dd,\hhw}_{\g}.$$

\begin{definition} [{\cite[Definition 6.1]{KKOP20C}}]
\label{def:complete dd}
A strong duality datum $\Dd$ in $\cato$ is called \emph{complete}
if, for any simple module $M \in \Ca_\g^0$, there exist simple
modules $M_k \in \Ca_\Dd$ $(k\in \Z)$ such that \ben
\item $M_k \simeq \one $ for all but finitely many $k$,
\item $M \simeq \hd ( \cdots \tens \D^2 M_2 \tens \D M_1 \tens \ M_0 \tens \D^{-1} M_{-1} \tens \cdots  ).$
\ee We say that a \pbw $(\Dd,\hhw)$ is complete if $\Dd$ is a
complete duality datum.
\end{definition}

Note that if $\Dd$ is complete, then $\cmf$ is the Cartan matrix of the simply-laced root system associated with $\g$  in \S \ref{subsec: E(M)} (\cite[Proposition 6.2]{KKOP20C}).

\begin{proposition}[{\cite[Theorem 6.9]{KKOP20C}}] \label{prop:Ca -infty,infty}
For any complete \pbw $(\Dd,\hhw)$, we have
$$  \Ca^{[-\infty,\infty],\Dd,\hhw}_{\g} = \Ca_\g^0.$$
\end{proposition}

\Lemma[{\cite[Lemma 6.8]{KKOP20C}}]
 Let $(\Dd,\hhw)$ be a complete \pbw,
$[a,b]$ an interval, and $M$ a simple module in $\catCO$. Then, $M$
belongs to $\Ca^{[a,b],\Dd,\hhw}_\g$ if and only if
$$\text{$\de(\D \bS{k},M)=0$ for $k>b$ and
$\de(\D^{-1}\bS{k},M)=0$ for $k<a$.}
$$
\enlemma

\Cor\label{cor:MX} Let $(\Dd,\hhw)$ be a complete \pbw, $[a,b]$
an interval, $M$ a real simple module in
$\Ca^{[a,b],\Dd,\hhw}_\g$ and $X$ a simple module in $\catCO$. If $
M\hconv X$ and $X\hconv M$ belong to $\Ca^{[a,b],\Dd,\hhw}_\g$, then
so does $X$. \encor \Proof For $k>b$, we have $\de(\D
\bS{k},M)=\de(\D^2 \bS{k},M)=0$. Hence Lemma~\ref{lem: normal} below
implies
$$0=\de(\D \bS{k}, X\hconv M)=\de(\D \bS{k}, X).$$
Similarly, for $k<a$, we have $\de(\D^{-1} \bS{k},M)=\de(\D^{-2}
\bS{k},M)=0$. Hence we have
$$0=\de(\D^{-1} \bS{k}, X\hconv M)=\de(\D^{-1} \bS{k}, X).$$
Therefore, $X$ belongs to $\Ca^{[a,b],\Dd,\hhw}_\g$. \QED

\Lemma\label{lem: normal} Let $L$, $M$ be real simple modules and
$X$ a simple module. \bnum
\item \label{it: LX}
If $\de(L,M)=\de(\D L,M)=0$, then we have
$$\de(L, X\hconv M)=\de(L,X).$$
\item \label{it: LX2}
If $\de(L,M)=\de(\D^{-1} L,M)=0$, then we have
$$\de(L, M\hconv X)=\de(L,X).$$
\ee \enlemma \Proof \eqref{it: LX}\ By the assumption, $(L,X,M)$ and $(X,M,L)$
is normal sequences. Hence by \cite[Lemma 4.16]{KKOP19C}, we have
$$\de(L, X\hconv M)=\de(L, X)+\de(L,M) =\de(L,X).$$

\snoi \eqref{it: LX2} can be similarly proved by using the normality of
$(L,M,X)$ and $(M,X,L)$. \QED

\subsection{\AD modules associated with $i$-boxes}

In this subsection, we introduce the notion of  $i$-boxes, assign a
simple module to each $i$-box, and then study its properties.

\begin{lemma} \label{lem: simply}
For any $a \in \Z$, we have $ \de(\bS{a^+},\bS{a})=1$.
\end{lemma}
\Proof Set $\im=\im_a$. By Proposition~\ref{prop:refl}, we may
assume that $a=1$. If $\ell=1$, then the assertion easily follows
from $\cusp_{a^+}=\D \cusp_a$ and
Proposition~\ref{prop: lemma of KKOP20C}\;\eqref{it:root}.
Assume that $\ell>1$. Then we have
$a=1<a^+\le\ell$. Hence we have $\bS{a}\simeq\F_\Dd(L(\im))$ and
$\bS{a^+}\simeq \F_\Dd(M)$, where $M\seteq\DC(\tw_{\le
a^+}\La_{\im},s_{\im}\La_{\im})$. Since $s_{\im}\tw_{\le a^+}\La_\im
\succeq\tw_{\le a^+}\La_\im$ and $s_\im\tw_{\le
a^+}\La_\im=s_{\im_2}\cdots s_{\im_{a^+}}\La_\im\preceq
s_{\im}\La_{\im}$, \cite[Proposition 10.2.4]{KKKO18} implies that
$\eps_\im(M)=-(\al_\im,\tw_{\le a^+}\La_{\im})$ and
$\eps_\im^*(M)=0$.

On the other hand,
 \cite[Corollary 3.8]{KKOP18} implies that
\eqn
\de(\cusp_{a^+},\cusp_a)=\de(L(\im),M)&&=\eps_\im(M)+\eps^*_\im(M)+\bl\al_\im,\wt(M)\br\\
&&=-(\al_\im,\tw_{\le a^+}\La_{\im})+ (\al_\im,\tw_{\le
a^+}\La_{\im}-s_{\im}\La_\im)=1.\qedhere \eneqn \QED

For $a,b \in \Z \sqcup \{ \pm \infty \}$,  an \emph{interval}
$[a,b]$ is the set of integers between $a$ and $b$:
\begin{align*}
 [a,b] \seteq \{ k \in \Z \ | \  a \le k \le b \}.
\end{align*}
If $a>b$, we understand $[a,b]=\emptyset$. We call $\max(b-a+1,0)$
the length of $[a,b]$ and denote it by $|\,[a,b]\,|$.

\begin{definition} Let $(\Dd,\hhw)$ be a \pbw.
\ben
\item We say that an interval $\ci=[a,b]$ is an \emph{$i$-box} if
$-\infty<a\le b<+\infty$ and $\im_a=\im_b$.
\item
For an $i$-box $[a,b]$, we set
$$ |[a,b]|_\phi \seteq  \big| \st{ s \mid \text{$s\in[a,b]$
and $\im_a=\im_s=\im_b$}}\big|.$$
\item For an $i$-box $[a,b]$, we define
\begin{align*}
M^{\Dd,\hhw}[a,b] \seteq \hd(\bS{b} \tens \bS{b^-} \tens \cdots \tens \bS{a^+} \tens \bS{a} ).
\end{align*}
We call $M^{\Dd,\hhw}[a,b]$ \emph{an affine determinantial module}.
\end{enumerate}
\end{definition}

Sometimes $M^{\Dd,\hhw}[a,b]$ with $a>b$ appears, and we understand
$M^{\Dd,\hhw}[a,b]=\one$ in such a case.

 For simplicity of notation, we sometimes drop the superscripts ${}^{\Dd,\hhw}$
if there is no risk of confusion.

\begin{proposition} \label{prop:unm}\hfill
\bnum
\item
For any $i$-box $[a,b]$, the module $M[a,b]$ is simple.
\item For $a\le b<a'\le b'$,
the pair $\bl M[a',b'], M[a,b] \br$ is strongly unmixed. \ee
\end{proposition}
\begin{proof}
(i) follows from Proposition~\ref{prop: s un and normal} and (ii)
follows from~\eqref{it: unmixed} in Proposition~\ref{prop: lemma of
KKOP20C}.
\end{proof}
Note that we shall see that $M[a,b]$ is real in Theorem~\ref{thm:
commuting ab} below.

\Lemma \label{lem: image} For any $i$-box $[a,b]\subset[1,\N]$, the module
$M[a,b]$ is isomorphic to the image of a determinantial module under
the functor $\F_\Dd$. More precisely, we have
$$   M^{\Dd,\hhw}[a,b] \simeq \F_\Dd(\DC_{\tw_0}[a,b]). $$
\enlemma

\begin{proof}
This follows from \cite[Proposition 4.1, Proposition 4.6]{KKOP18}
and the fact that $\F_\Dd$ is an exact monoidal functor sending simple
modules to simple modules.
\end{proof}

\begin{corollary}\label{cor:Tsystem}  For an $i$-box $[a,b]$ with $ b-a +1 \le \N$,
we have an exact sequence in $\Ca_\g$ as follows$\col$
\begin{align*}
0 \to  \dtens_{ \substack{\jm \in  \Iff; \;
d(\im_a,\jm)=1}} M[a(\jm)^+,b(\jm)^-]  \to   M[a^+,b] \tens M[a,b^-]
\to   M[a,b] \tens M[a^+,b^-] \to 0.
\end{align*}
\end{corollary}

\begin{proof}
The above exact sequence is the image of~\eqref{eq: T-system in QH}
under the exact functor $\F_\Dd$.
\end{proof}

\subsection{Commuting condition between \Ad  modules}
In this subsection, we will give a sufficient condition that
simple modules $M[a_1,b_1]$ and $M[a_2,b_2]$ commute. We start with
the following definition on the pair of $i$-boxes $[a_1,b_1]$ and
$[a_2,b_2]$.

\begin{definition}
We say that $i$-boxes $[a_1,b_1]$ and $[a_2,b_2]$ \emph{commute} if
we have either
$$a_1^- < a_2 \le b_2 <  b_1^+ \quad  \text{ or } \quad a_2^- < a_1 \le b_1 <  b_2^+.$$
\end{definition}

The main goal of this subsection is to prove the theorem
which tells us that \Ad modules $M[a_1,b_1]$ and  $M[a_2,b_2]$
commute if $[a_1,b_1] $ and $[a_2,b_2]$ commute (see
Theorem~\ref{thm: commuting ab} below).

\Prop \label{prop: commuting main} For any $i$-box $[a,b]$ and $s \in
\Z$ such that $a^-<s < b^+$, we have
$$\de(M[a,b],\bS{s})=0.$$
\enprop

\begin{proof}
We shall argue by induction on $|[a,b]|_{\phi}$.

\mnoi {\rm (i)} Assume first that $a<s<b$ and $\im_a=\im_b=\im_s$.
Note that
$$M[a,b] = \hd( M[s^+,b] \tens \bS{s} \tens M[a,s^-])$$
and hence
$$ \de(\bS{s},M[a,b]) = \de\left(\bS{s},\hd( M[s^+,b] \tens \bS{s} \tens M[a,s^-]) \right).$$
Since the pairs $(M[s^+,b],\bS{s})$ and $(\bS{s}, M[a,s^-])$ are
unmixed by Proposition~\ref{prop:unm}, we have \eqn
&&\de\left(\bS{s},\hd( M[s^+,b] \tens \bS{s} \tens M[a,s^-] ) \right)\\
&&\hs{10ex}=\de\left(\bS{s},\hd( M[a,s^+] \tens \bS{s}) \right) +\de\left(\bS{s},\hd( \bS{s} \tens M[s^-,b]) \right) \\
&&\hs{10ex} = \de\left(\bS{s}, M[a,s] \right)+
\de\left(\bS{s},  M[s,b]  \right)  =0 \eneqn by
Proposition~\ref{prop: three}. Here, the last equality follows from
the induction hypothesis on $|[a,b]|_{\phi}$.

\mnoi {\rm (ii)} Assume that $a<s<b$ and $\im_a=\im_b \ne \im_s$.
Set $c=\min\st{k\in[a,b]\mid\im_k=\im_b, s<k}$. Then, we have $a\le
c^- < s < c\le b$, and
$$ M[a,b]=\hd(  M[c,b] \tens  M[a,c^-]).$$
Now our assertion follows from the induction hypothesis.

\mnoi (iii) It remains to prove the cases
$$
{\rm (a)} \ \ a^- < s \le a \qtq {\rm (b)} \ \ b \le s < b^+.
$$

We shall prove only the case  {\rm (b)} since the other case is
similarly proved.

\snoi {\rm (iii-1)} Assume first that $s>a+\N$. Then we have
$$ M[a,b] =  M[a^+,b] \hconv \bS{a}.  $$
On the other hand we have $\de(\bS{s},\bS{a})=0$ and
$\de(\bS{s},M[a^+,b])=0$ by Proposition~\ref{prop: lemma of KKOP20C}
and the induction hypothesis, respectively. Thus our assertion
follows from Proposition~\ref{pro:subquotient}.

\mnoi {\rm (iii-2)} Assume that $s \le a+\N$. We may assume that
$a=1$ by Proposition~\ref{prop:refl}.  Hence it is enough
to show that
$$\text{$M[1,b]$ commutes with $\bS{s}$ if
$\im_1=\im_b$, $b \le s < b^+$ and $s \le 1+\N$.}$$

Assume first that $s \le \N$. Then we have  $\bS{s} =
\F_\Dd(\DC[s] )$ and $\F_\Dd\bl\DC(\tw_{\le
s}\La_{\im_1},\La_{\im_1})\br\simeq M[1,b]$ by Lemma~\ref{lem:
image}. Then the assertion follows from the fact that
 $\DC[s]\in R_{\tw_{\le s}}\gmod $ and
$\DC(\tw_{\le s}\La_{\im},\La_{\im})$ commutes with every simple
module in $R_{\tw_{\le s}}\gmod$  by Theorem~\ref{thm: KKOP18main}.

\medskip
Now we assume that $s=\N+1$.
In this case, we have $\cusp_s=\D\cusp_1$.

We divide this case into two
sub-cases:

\noindent (1) Assume that $b=\N+1=s$. Then we have
\begin{align*}
M[a,b] &= \bS{b} \hconv M[1,b^-] = \bS{s} \hconv
\F_\Dd(\DC(w_0\La_{\im_1},\La_{\im_1})).
\end{align*}

By Proposition~\ref{prop: de ve} and Lemma~\ref{lem: image}, we have
$$ \de(\bS{s},M[1,b^-])=\de(\D \bS{1},M[1,b^-]) =\ve_{\im_1}(\DC(w_0\La_{\im_1},\La_{\im_1})).$$
Then \cite[Lemma 9.15]{KKKO18} tells that
$$\ve_{\im_1}(\DC(w_0\La_{\im_1},\La_{\im_1} )) =
\max\bl 0,-(\al_{\im_1},w_0\La_{\im_1})\br\le1.$$ Now
\cite[Proposition 2.17]{KKOP20C} implies that
$$ \de(\bS{s},M[a,b])=0,$$
yielding our assertion in this case.

\snoi (2) Assume that  $b\le\N$. Then we have $b < s=\N+1 < b^+$ and
hence $\im_1^*=\im_{\N+1} \ne \im_b=\im_1$.
Since $w_0\La_{\im_1}=s_{\im_1}\cdots s_{\im_\ell}\La_{\im_1}=s_{\im_1}\cdots s_{\im_b}\La_{\im_1}$, we have $M[1,b]=\DC(w_0\La_{\im_1},\La_{\im_1})$.
Hence we have
\eqn
 \de(\bS{s},M[1,b])=\de(\D \bS{1},M[1,b])&&=
\ve_{\im_1}\bl\DC(w_0\La_{\im_1},\La_{\im_1})\br\\* &&= \max\bl
0,-(\al_{\im_1},w_0\La_{\im_1})\br = \max\bl
0,(\al_{\im^*_1},\La_{\im_1})\br=0, \eneqn
where the last equality follows from $\im_1^*
\ne \im_1$.
Thus we obtain the desired result.
\end{proof}

\begin{proposition} \label{prop: de1}
For any $i$-box $[a,b]$, we have
$$ \de(\bS{b^+},M[a,b]) = \de(\bS{a^-},M[a,b])=1.$$
\end{proposition}

\begin{proof}
We shall only prove that $\de(\bS{b^+},M[a,b]) =1$ since the other
assertion can be proved similarly. Since $\bS{b^+}\hconv
M[a,b]\simeq M[a,b^+]$, we have $M[a,b]\simeq
M[a,b^+]\hconv\D\bS{b^+}$.
Noting that $\bS{b^+}$ is a root module (see Definition~\ref{def: root module}),
 we have
$$0\le\de(\bS{b^+}, M[a,b])\le
\de(\bS{b^+}, M[a,b^+])+\de(\bS{b^+}, \D \bS{b^+})=1.$$

Assume that $\de(\bS{b^+}, M[a,b])=0$. Then we have $M[a,b^+]\simeq
\bS{b^+}\tens M[a,b]$, and
$$0=\de\bl \bS{b},M[a,b^+]\br=\de\bl \bS{b},\bS{b^+}\br+\de\bl \bS{b},M[a,b]\br
=1$$ by Lemma~\ref{lem: simply}, which is a contradiction. Thus we
conclude that $\de(\bS{b^+}, M[a,b])=1$.
\end{proof}

\begin{theorem} \label{thm: commuting ab} \hfill
\bna
\item \label{it: Mab real} For any $i$-box $[a,b]$, $M[a,b]$ is a real simple module in $\Ca_\g$.
\item \label{it: Mab com}  If two $i$-boxes  $[a_1,b_1]$ and $[a_2,b_2]$ commute, then
$M[a_1,b_1]$ and $M[a_2,b_2]$ commute.
\end{enumerate}
\end{theorem}

\begin{proof}
By induction on $|\,[a,b]\,|_\phi$, we may assume that
 $a\le b^-$ and   $M[a,b^-]$ is real simple.
Hence ~\eqref{it: Mab real} follows from $M[a,b]=\bS{b}\hconv M[a,b^-]$,
Lemma~\ref{cor: real 1} and Proposition~\ref{prop: de1}.

\snoi \eqref{it: Mab com} is an immediate consequence of Proposition~\ref{prop:
commuting main}.
\end{proof}

\Lemma\label{lem: at most simply linked} For any $i$-box $[a,b]$, we
have
$$\de(M[a,b],M[a^-,b^-])\le1.$$
\enlemma

\begin{proof}
Note that $M[a^-,b^-] = M[a,b^-] \hconv \bS{a^-}$.
Thus Theorem~\ref{thm: commuting ab} and Proposition~\ref{prop: de1}
tell that
\[\de(M[a,b],M[a^-,b^-]) \le \de(M[a,b],M[a,b^-])+\de(M[a,b], \bS{a^-}) =1.\qedhere\]
\end{proof}
Later, we see that $\de(M[a,b],M[a^-,b^-])=1$ in the course of the
proof of Theorem~\ref{th:Tsystem}.
\Lemma\label{lem:unique} For  any
$i$-box $[a,b]$, we have
$$\de(\D \bS{b},M[a,b])=1\qtq \de(\D ^{-1}\bS{a},M[a,b])=1.$$
\enlemma

\Proof Since $M[a,b]\simeq\bS{b}\hconv M[a,b^-]$ and $\de(\D^k
\bS{b},M[a,b^-])=0$ for $k\ge1$, the first equality  $\de(\D
\bS{b},M[a,b])=1$ follows from \cite[Lemma 6.7]{KKOP20C}. The second
equality can be proved similarly. \QED

For an interval $[a,b]$ with $a\le b$, we define $i$-boxes
\begin{align}\label{def:int}
[a,b\}\seteq[a,b(\im_a)^-] \quad \text{ and } \quad \{a,b]\seteq[a(\im_b)^+,b].
\end{align}
 For notations, see \eqref{eq: nota +,-}.

\Lemma\label{lem:s=b}
If $a<b$,
then $M[a,b\}$ and $\cusp_{b^-}$ commute.
\enlemma

\Proof
Set $c=b(\im_a)^-$ so that $a\le c\le b<c^+$ and $M[a,b\}=M[a,c]$.

\snoi
(i) If $b^->a^-$, then  $M[a,c]$ and $\cusp_{b^-}$ commute
since $a^-<b^-<b<c^+$.

\snoi
(ii) If $b^-\le a^-$, then  $\im_a\not=\im_b$ and hence $c<b$.
Then
$(b^-)^-<b^-\le a^-<a\le c<b=(b^-)^+$, which implies that
$M[a,c]$ and $\cusp_{b^-}=M[b^-,b^-]$ commute.
\QED

\subsection{T-systems among \Ad modules}  A T-system was introduced in \cite{KNS} as a system of differential equations associated with solvable lattice models. It was conjectured
in \cite{FR99} that the {\it q-characters} of Kirillov-Reshetikhin
modules solve the T-system.
The T-system can
be written as a short exact sequence in terms of KR-modules.  In
\cite[Appendix]{OT19B}, T-systems are  presented in terms of
notations $W^{(i)}_{k,x}$ and also $V(i^k)_x$ for $i\in I_0$, $k \in
\Z_{\ge1}$ and $x \in \ko^\times$.
For instance,  for  a
simply-laced untwisted $\g$ and $i\in I_0$, T-system is given as follows:
\begin{align*}
0 \to \dtens_{ \substack{  j \in I_0 ;\;  d(i,j)=1}}
W^{(j)}_{k,x(-q)} \to  W^{(i)}_{k,x(-q)^2} \tens W^{(i)}_{k,x} \to
W^{(i)}_{k-1,x(-q)^2} \tens W^{(i)}_{k+1,x} \to 0,
\end{align*}
for any $i \in I_{0}$  and $x\in\cor^\times$.

In our context, it is paraphrased
as follows (see Theorem~\ref{Thm: KR as AD} below).

\begin{theorem}\label{th:Tsystem}
For an arbitrary quantum affine algebra $U_q'(\g)$ and an arbitrary
$i$-box $[a,b]$ such that $a<b$, we have an exact sequence
\begin{align} \label{eq: T-system in terms of M[a,b] again}
0 \to \hs{-2ex} \dtens_{ \substack{ \jm \in\Iff; \; d(\im_a,\jm)=1}} \hs{-2ex}
M[a(\jm)^+,b(\jm)^-]  \to   M[a^+,b]  \tens M[a,b^-]    \to   M[a,b]
\tens M[a^+,b^-] \to 0.
\end{align}
{\em We call it also a T-system.}
\end{theorem}

Note that the left term and the right term in~\eqref{eq: T-system in
terms of M[a,b] again} are simple.

By this theorem, we have also an exact sequence
\begin{align}
0 \to M[a,b]\tens M[a^+,b^-] \to   M[a,b^-]\tens  M[a^+,b]   \To
\hs{-2ex} \dtens_{ \substack{ \jm \in \Iff; \; d(\im_a,\jm)=1}} \hs{-2ex}
M[a(\jm)^+,b(\jm)^-]
 \to 0.
\end{align}

\medskip

The rest of the subsection is devoted to the proof of this theorem.
We begin with the following lemma. \Lemma\label{lem:xyz} Let $X,Y,Z$
be real simple modules in $\catg$, and we assume the following
conditions: \bnum
\item \label{it:cond1}
$(X,Y,Z)$ is a normal sequence,
\item \label{it:cond2}
$\La(Y, X)+\La(Y,Z)-\La(Y,X\hconv Z)=2\de(X,Y)$. \ee Then, the
composition $X\tens Y\tens Z\To[{\rmat{X,Y}}] Y\tens X\tens Z\epito
Y\hconv(X\hconv Z)$ is an epimorphism and it induces an isomorphism
$\hd(X\tens Y\tens Z)\simeq Y\hconv(X\hconv Z)$. \enlemma
\begin{proof}
Set $\la=\de(X,Y)$. Then we have a commutative diagram:
$$\xymatrix@C=10ex{
X\tens Y_{\e^z}\tens Z\ar[d]^{\Rren_{X,Y_{\e^z}}}
\ar[drr]^-{z^\la\id_{X\tens Y_{\e^z}\tens
Z}}\ar@/^1.5pc/[drrrr]|-{z^\la\Rren_{Y_{\e^z},Z}}
&&\ar@{}[dll]|(.7){\textcircled{\scriptsize A}}\\
 Y_{\e^z}\tens X\tens Z\ar[rr]_{\Rren_{Y_{\e^z},X}}\ar@{->>}[d]\ar@{}[rrrrd]_-{
{\textcircled{\scriptsize B}}} &&X\tens Y_{\e^z}\tens
Z\ar[rr]_{\Rren_{Y_{\e^z},Z}}
&&X\tens Z\tens Y_{\e^z}\ar@{->>}[d]\\
 Y_{\e^z}\tens (X\hconv Z)\ar[rrrr]_{z^\la\Rren_{Y_{\e^z},X\nabla Z}}&&&&(X\hconv Z)\tens Y_{\e^z}
}
$$
up to constant in $\cor[[z]]^\times$. Here the commutativity of
$\textcircled{\scriptsize A}$ follows from Proposition~\ref{prop:
de(M,N)} and the one of $\textcircled{\scriptsize B}$ from
Lemma~\ref{lem:Rmat} and the hypothesis~\eqref{it:cond2}. Hence, the
following diagram is commutative:
$$\xymatrix@C=10ex{
X\tens Y_{\e^z}\tens Z\ar[d]^{\Rren_{X,Y_{\e^z}}}
\ar[drrrr]^-{\Rren_{Y_{\e^z},Z}}\\
 Y_{\e^z}\tens X\tens Z\ar@{->>}[d]&&&&X\tens Z\tens Y_{\e^z}\ar@{->>}[d]\\
 Y_{\e^z}\tens (X\hconv Z)\ar[rrrr]_{\Rren_{Y_{\e^z},X}}&&&&(X\hconv Z)\tens Y_{\e^z}
}
$$
Setting $z=0$, we obtain a commutative diagram:
$$\xymatrix@C=10ex{
X\tens Y\tens Z\ar[d]^{\rmat{X,Y}}
\ar[dr]^-{\rmat{Y,Z}}\ar@{->>}[rrr]&&&\hd(X\tens Y\tens Z)\akete[-1.2ex]\ar@{>->}[dd]\\
 Y\tens X\tens Z\ar@{->>}[d]&X\tens Z\tens Y\ar@{->>}[d]\ar[dr]^{\rmat{X,Z}}\\
 Y\tens (X\hconv Z)\ar[r]_{\rmat{Y,X\nabla Z}}&(X\hconv Z)\tens Y\akew[1ex]
\ar@{>->}[r]&Z\tens X\tens Y\ar[r]_{\rmat{X,Y}}&Z\tens Y\tens X }
$$
Let $S$ be the image of the composition
$$ Y\tens X\tens Z\epito Y\tens (X\hconv Z)\To[{\rmat{Y,X\nabla Z}}]
(X\hconv Z)\tens Y\monoto Z\tens X\tens Y.$$ Then $S\simeq Y\hconv
(X\hconv Z)$ is a simple module. By the above commutative diagram,
the image $K$ of the composition of $X\tens Y\tens
Z\To[{\rmat{Y,Z}}]X\tens Z\tens Y\To[{\rmat{X,Z}}] Z\tens X\tens Y$
is contained in $S$. On the other hand, the image of $K\to Z\tens
Y\tens X$ is isomorphic to the simple module $\hd(X\tens Y\tens Z)$.
Therefore $K$ is non-zero and hence $K=S$, and it is isomorphic to
 $\hd(X\tens Y\tens Z)$.
\QED

\Proof[{Proof of {\rm Theorem~\ref{th:Tsystem}}}] We argue by
induction on $\disp{[a,b]}\ge2$.

Assume first that $\disp{[a,b]}=2$, i.e., $b=a^+$. If $\ell>1$ then
we have $b-a+1\le\ell$ and hence we have the desired result by
Corollary~\ref{cor:Tsystem}. If $\ell=1$, then we have $\bS{b}=\D
\bS{a}$ and  we have an exact sequence
$$0\To\one\To\D \bS{a} \tens \bS{a} \To \D \bS{a} \hconv\bS{a}\To0$$
since $\cusp_a$ is a root module.
Hence, for any $\ell$, we have obtained an exact sequence:
\begin{align} \label{eq: sT-system}
0 \to \hs{-2ex} \dtens_{   \jm \in  \Iff; \; d(\im_a,\jm)=1} \hs{-2ex}
M[a(\jm)^+,(a^+)(\jm)^-]  \to\cusp_{a^+}\tens\cusp_{a}\to\cusp_{a^+}\htens\cusp_{a}\to 0.
\end{align}

\mnoi Now we can assume $\disp{[a,b]}>2$.

In order to prove the theorem, it is enough to show \ben[(I)]
\item $\de(M[a,b^-],M[a^+,b])\le 1$, \label{item a}
\item $M[a^+,b]\hconv M[a,b^-]\simeq M[a,b]\tens M[a^+,b^-]$, \label{item b}
\item $M[a,b^-]\hconv M[a^+,b]\simeq \dtens_{d(\im_a,\jm)=1}M[a(\jm)^+,b(\jm)^-]$,
 \label{item c}
\item $M[a^+,b]\hconv M[a,b^-]\not\simeq M[a,b^-]\hconv M[a^+,b]$.\label{item d}
\ee Indeed, \eqref{item a} and \eqref{item d} imply
$\de(M[a,b^-],M[a^+,b])=1$, and \cite[Proposition~4.7]{KKOP19C}
implies the existence of a short exact sequence
$$0\To  M[a,b^-]\hconv M[a^+,b]\To
M[a^+,b]\tens M[a,b^-]\To M[a^+,b]\hconv M[a,b^-]\To0.$$

\mnoi Now let us show (I)--(IV).

\snoi \eqref{item a} follows from Lemma~\ref{lem: at most simply
linked}.

\snoi Let us show \eqref{item b}. The composition
$$M[a^+,b]\tens M[a,b^-]\monoto
M[a^+,b]\tens\bS{a} \tens M[a^+,b^-]\epito M[a,b]\tens M[a^+,b^-]$$
does not vanish by Proposition~\ref{prop:rcomp}. Since $M[a,b]$ and
$M[a^+,b^-]$ commute by Theorem~\ref{thm: commuting ab}, $
M[a,b]\tens M[a^+,b^-]$ is a simple module. Hence we obtain an
epimorphism $M[a^+,b]\tens M[a,b^-]\epito  M[a,b]\tens M[a^+,b^-]$.
Thus we obtain \eqref{item b}.

\mnoi Let us show \eqref{item c}.

\noi We shall first prove \eq \hd\bl
M[a^{++},b]\tens\bS{a}\tens\bS{a^+}\br\simeq \bS{a}\hconv
M[a^{+},b].\label{claim1} \eneq

In order to prove \eqref{claim1}, we shall apply Lemma~\ref{lem:xyz}
with $X\seteq M[a^{++},b]$, $Y\seteq \bS{a}$, and $Z\seteq\bS{a^+}$.
Since $(X,Z)$ is unmixed, $(X,Y,Z)$ is a normal sequence.

We shall verify~\eqref{it:cond2} in Lemma~\ref{lem:xyz}. Note
that~\eqref{it:cond2} is equivalent to
$$\La(Y,Z)-\La(Y,X\hconv Z)=\La(X,Y).$$
We have $\La(M[a^{++},b]\hconv \bS{a^+},\bS{a})=
\La(M[a^{++},b],\bS{a})+\La(\bS{a^+},\bS{a})$ because
$(M[a^{++},b],\bS{a})$ is unmixed. On the other hand, we have
$$2=2\de(M[a^+,b],\bS{a})=\La(M[a^{++},b]\hconv \bS{a^+},\bS{a})
+\La(\bS{a},M[a^{++},b]\hconv \bS{a^+})$$ and
$$2=2\de(\bS{a^+},\bS{a})=\La(\bS{a^{+}},\bS{a})
+\La(\bS{a},\bS{a^+}).$$ Hence we obtain
$$2-\La(\bS{a},M[a^{++},b]\hconv \bS{a^+})=
\La(M[a^{++},b],\bS{a})+\bl 2-\La(\bS{a},\bS{a^+})\br,$$ which
implies~\eqref{it:cond2}. Hence, Lemma~\ref{lem:xyz} implies
\eqref{claim1}.

\snoi Thus, we have a chain of morphisms \eqn
M[a^+,b^-]\tens\bS{a}\tens M[a^{+},b]
\epito&& M[a^+,b^-]\tens\bl\bS{a}\hconv M[a^{+},b]\br\\
\simeq&& M[a^+,b^-]\tens \hd\bl M[a^{++},b]\tens\bS{a}\tens\bS{a^+}\br\\
\simeq&& M[a^+,b^-]\tens \bl M[a^{++},b]\hconv
(\bS{a}\hconv\bS{a^+})\br. \eneqn

On the other hand, the pair  $(M[a^{+},b^-], \bS{a}\hconv\bS{a^+})$
is unmixed since
\begin{align*}
\bS{a}\hconv\bS{a^+}\simeq\dtens\limits
_{ \substack{ \jm \in \Iff; \; d(\im_a,\jm)=1}} M[a(\jm)^+,(a^+)(\jm)^-]\qt{by \eqref{eq: sT-system}
\quad  and} (a^+)(\jm)^-<a^+.
\end{align*}
Therefore, $(M[a^+,b^-],  M[a^{++},b],\; \bS{a}\hconv\bS{a^+})$
is a normal sequence
and
 we have an epimorphism
$$M[a^{+},b^-]\tens\bS{a}\tens  M[a^{+},b]
\epito \bl M[a^{+},b^-]\hconv
M[a^{++},b]\br\hconv(\bS{a}\hconv\bS{a^+}).
$$

By the induction hypothesis, we have
$$M[a^{+},b^-]\hconv M[a^{++},b]\simeq
\dtens\limits _{d(\jm,\im_a)=1}M[(a^+)(\jm)^+,b(\jm)^-].$$ Hence
\cite[Lemma~3.2.22]{KKKO18} implies that \eqn
&&\bl M[a^{+},b^-]\hconv M[a^{++},b]\br\hconv(\bS{a}\hconv\bS{a^+})\\
&&\hs{10ex}\simeq \dtens\limits _{ \substack{ \jm \in \Iff; \; d(\im_a,\jm)=1}}
\bl M[(a^+)(\jm)^+,b(\jm)^-] \hconv M[a(\jm)^+,(a^+)(\jm)^-]\br\\
&&\hs{10ex}\simeq A, \eneqn where $A\seteq\dtens\limits
_{  \jm \in \Iff; \; d(\im_a,\jm)=1} M[a(\jm)^+,b(\jm)^-]$.

Thus we obtain an epimorphism
\begin{align*}
M[a^+,b^-]\tens\bS{a}\tens M[a^{+},b]\epito A.
\end{align*}

 \mnoi Since $\de(M[a^{+},b^-],\bS{a})=1$, we have an exact sequence
$$0\To \bS{a}\hconv M[a^{+},b^-]\To
M[a^{+},b^-]\tens \bS{a}\To M[a,b^-]\To0.$$

\vs{1ex} Now, we shall show

\noi \eq\text{ $(\bS{a}\hconv M[a^{+},b^-])$ and $M[a^{+},b]$
commute.}\label{claim2} \eneq We have \eqn \de\bl\bS{a} \hconv
M[a^{+},b^-],\ms{3mu}M[a^{+},b]\br \le \de\bl\bS{a},M[a^{+},b]\br+
\de\bl M[a^{+},b^-],M[a^{+},b]\br=1+0=1. \eneqn Hence it is enough
to show
$$\de\bl\bS{a}\hconv M[a^{+},b^-],\ms{3mu}M[a^{+},b]\br
\not=\de\bl\bS{a},M[a^{+},b]\br+ \de\bl M[a^{+},b^-],M[a^{+},b]\br.
$$
If the equality holds, then \cite[Lemma 2.23]{KKOP20C} implies that
$\bl M[a^{+},b],\bS{a},M[a^{+},b^-]\br$ is a normal sequence. Hence
we have \eqn \La\bl  M[a^{+},b],M[a^{+},b^-])+ \La\bl\bS{a}
,M[a^{+},b^-]\br-\La\bl  M[a^{+},b]\hconv\bS{a},M[a^{+},b^-]\br =0.
\eneqn On the other hand, we have
$$\La\bl  M[a^{+},b],M[a^{+},b^-]\br=-\La\bl M[a^{+},b^-],M[a^{+},b]\br,$$
and \eqn \La\bl  M[a^{+},b]\hconv\bS{a},M[a^{+},b^-]\br
&&=-\La\bl M[a^{+},b^-], M[a^{+},b]\hconv\bS{a}\br\\
&&=-\La\bl M[a^{+},b^-], M[a^+,b]\br-\La\bl M[a^{+},b^-], \bS{a}\br,
\eneqn where the last equality follows from the normality of $\bl
M[a^{+},b^-],M[a^+,b],\bS{a}\br$.

Thus, we obtain \eqn 0&&=\La\bl  M[a^{+},b],M[a^{+},b^-])+
\La\bl\bS{a} ,M[a^{+},b^-]\br-\La\bl  M[a^{+},b]\hconv\bS{a},M[a^{+},b^-]\br\\
&&=-\La(M[a^{+},b^-],M[a^{+},b])+\La\bl\bS{a}  ,M[a^{+},b^-]\br\\
&&\hs{28ex}-\bl-\La (M[a^{+},b^-], M[a^+,b])-\La( M[a^{+},b^-], \bS{a})\br\\
&&=2\de\bl \bS{a}  ,M[a^{+},b^-]\br, \eneqn which contradicts
$\de\bl \bS{a}  ,M[a^{+},b^-]\br=1$. Thus we obtain \eqref{claim2}.
In particular, $\bl\bS{a}\hconv M[a^{+},b^-]\br\tens  M[a^{+},b]$ is
a simple module.

\mnoi
\eqn
\bullet&&\text{The simple modules $\bl\bS{a}\hconv M[a^{+},b^-]\br\tens
M[a^{+},b]$
and $A$ are not isomorphic.}
\eneqn

\noi
Indeed, we have $\de\bl\D \bS{b},\ms{3mu}(\bS{a}\hconv
M[a^{+},b^-])\tens  M[a^{+},b]\br \ge\de\bl\D \bS{b},
M[a^{+},b]\br=1$ by Lemma~\ref{lem:unique} and $\de\bl\D
\bS{b},A\br=0$.

\vs{1ex} Thus the composition
$$ \bl\bS{a}\hconv M[a^{+},b^-]\br\tens  M[a^{+},b]
\monoto M[a^{+},b^-]\tens \bS{a} \tens  M[a^{+},b] \epito A$$
vanishes, and hence $ M[a^{+},b^-]\tens \bS{a}\tens  M[a^{+},b]
\epito A$ factors through $ M[a,b^-]\tens  M[a^{+},b]$. Thus we
obtain an epimorphism
$$ M[a,b^-]\tens  M[a^{+},b]\epito A,$$
which completes the proof of \eqref{item c}.

\mnoi Finally, \eqref{item d} follows from $\de\bl \D \bS{b},
M[a,b^-]\hconv M[a^{+},b] \br= \de\bl \D \bS{b}, M[a,b]\tens
M[a^+,b^-] \br=1$ and $\de\bl \D \bS{b}, M[a^+,b]\hconv M[a,b^-]
\br= \de\bl \D
\bS{b},\dtens_{\jm \in \Iff; \; d(\im_a,\jm)=1}M[a(\jm)^+,b(\jm)^-]\br=0$.
\end{proof}

\vskip 2em
\section{Admissible chains of $i$-boxes} \label{sec: Admissible chain}
In this section, we construct commuting families
of real simple modules consisting of \Ad modules and
then investigate relations among the families, called box moves.
Throughout this section,
we fix a \pbw $(\Dd,\hhw)$
and keep the data and notations in the previous section.

\subsection{Chains of $i$-boxes}
Recall $[a,b\}$ and $\{a,b]$ defined in \eqref{def:int}.
\begin{definition}  \label{def: ad ch of i box}\hfill
\ben
\item A chain $\frakC$ of $i$-boxes
$$ (\ci_k = [a_k,b_k] )_{1 \le k \le l} \qt{for $l \in \Z_{\ge 1} \sqcup \{ \infty \}$}$$
is called \emph{admissible} if
$$\text{$\tc_k = [\ta_k,\tb_k] \seteq \bigcup_{1 \le j \le k} [a_j,b_j]$ is an interval with $|\tc_k| = k$ for $k=1,\ldots,l$}$$
and
$$ \text{either   $[a_k,b_k] = [\ta_k,\tb_k\}$ or $[a_k,b_k] = \{\ta_k,\tb_k]$  for $k=1,\ldots,l$.}
$$

\item The interval $\tc_k$  is called
the {\em envelope} of $\ci_k$, and $\tc_l$ is called the {\em range}
of $\frakC$.
\end{enumerate}
\end{definition}
Note that the chain is uniquely determined by its envelopes:
\begin{align}\label{eq: T_k}
\ci_k=[a_k,b_k] \seteq T_{k-1}[\ta_{k},\tb_{k}] =
\bc  \ [\ta_{k},\tb_{k}\} &\text{ {\rm (i)} if $\ta_k=\ta_{k-1}-1$,}\\
\ \{\ta_{k},\tb_{k}]&\text{ {\rm (ii)} if $\tb_k=\tb_{k-1}+1$,} \ec
\end{align}
for $1<k\le l$. In case  {\rm (i)} in~\eqref{eq: T_k}, we write
$T_{k-1}=\LL$, while $T_{k-1}=\RR$ in case (ii).

\medskip

Thus, to an admissible chain of $i$-boxes of length $l$, we can
associate a pair $(a, \mathfrak{T})$ consisting of an integer $a$
and a sequence $\mathfrak{T} = ( T_1,T_2,\ldots, T_{l-1})$ such that
$T_i \in \{ \LL,\RR\}$ $(1 \le i < l)$,
$$ a= a_1=b_1=\ta_1=\tb_1 ,  \quad  [\ta_k,\tb_k] =
\begin{cases}
[\ta_{k-1}-1,\tb_{k-1}] & \text{ if } T_{k-1}=\LL, \\
[\ta_{k-1},\tb_{k-1}+1]  & \text{ if } T_{k-1}=\RR.
\end{cases}$$
Note that $[a_k,b_k]=[\ta_k,\tb_k\}$ or $\{\ta_k,\tb_k]$ according
that $T_{k-1}=\LL$ or $T_{k-1}=\RR$.

\Lemma\label{lem:chain} Let $\frakC=\seq{\ci_k}_{1\le k \le l}$ be an
admissible chain of $i$-boxes, and $\seq{\tc_k}_{1\le k \le l}$ its
envelopes. Then, we have \bnum
\item \label{it: a-b+} $(a_k)^-<\ta_k\le\tb_k<(b_k)^+$,
\item \label{it: a-b+2}
$(a_k)^-<a_j\le b_j<(b_k)^+$ if $1\le j\le k\le l$.\\
In particular, $\ci_j$ and $\ci_k$ commute. \ee \enlemma \Proof \eqref{it: a-b+}
Set $\im=\im_{a_k}=\im_{b_k}$. Then we have $a_k=\ta_k(\im)^+$.
Hence we have $(a_k)^-<\ta_k$. The other strict inequality on $b_k$
can be proved similarly.

\snoi \eqref{it: a-b+2}  This assertion follows from
\[ (a_k)^-<\ta_k\le a_j\le b_j\le\tb_k<(b_k)^+. \qedhere\]
\QED

For an admissible chain of $i$-boxes $\frakC=\seq{\ci_k}_{1\le k \le l}$
and $1 \le t \le l$, we set
$$  \frakC_{\le t} \seteq \seq{\ci_k}_{1\le k \le t}.$$

\begin{proposition} \label{prop: commuting and cotained}
Let $\frakC=\seq{\ci_k}_{1\le k \le l}$ be an admissible chain of
$i$-boxes and let $\ci=[a,b]$ be an $i$-box such that $\ci\subset
\tc_l$ and $\ci$ commutes with all $\ci_k$ $(1 \le k \le l)$. Then
there exists $s\in[1,l]$ such that $\ci = \ci_s$.
\end{proposition}

\begin{proof}
We can assume that $\ci \subset \tc_l$ and $\ci \not\subset
\tc_{l-1}$ and $l \ge 2$, which implies
$$\text{(i)} \ \  a=\ta_l<\ta_{l-1} \quad \text{ or } \quad \text{(ii)} \ \ b=\tb_l >\tb_{l-1}.$$
Since the case (ii) can be proved similarly, we
only treat the case (i).

Then $\ci_l=[a_l,b_l]$ satisfies $a =a_l=\ta_l$ and
 $b_l=\tb_l(\im_b)^-$. Hence, we have $b\le b_l$.
If $b=b_l$, we have $\ci=\ci_l$. Thus we may assume that $b < b_l$.
Hence $b^+\le b_l\le \tb_{l}=\tb_{l-1}$. Take the smallest $s\ge1$
such that $b^+ \le \tb_s$. Then $1\le s\le l-1$.

Assume first $s>1$. Then we have $\tb_{s-1}<b^+\le \tb_s$, and hence
$b^+=\tb_s=b_s$. Since $\ci$ and $\ci_s$ commute and $b^+ \le b_s$,
we have
$$  a_{s}^- < a \le b < b_{s}^+.$$
Hence we have $a_s^-<a<\ta_{l-1}\le a_s$, which contradicts
$\im_{a_s}= \im_{b_s}=\im_{b^+}=\im_a$.

\medskip
Hence we have $s=1$. Therefore, $\ta_l<b^+\le \ta_1$. Now take the
smallest $t\ge1$ such that $\ta_t\le b^+$. Since $\ta_l\le a\le
b<b^+$, we have $\ta_{l-1}=1+\ta_l\le b^+$. Hence we have $1\le t\le
l-1$.

If $t>1$, then we have $\ta_t\le b^+<\ta_{t-1},$ and hence
$b^+=\ta_{t}=a_t$. If $t=1$, then $b^+=\ta_1=a_1$. In any case, we
have $b^+=a_t$.

Since $\ci$ and $\ci_t$ commute and $b^+ \le b_t$, we have
$$  a_{t}^- < a \le b < b_{t}^+.$$
Hence, we obtain
 $a_{t}^- < a=\ta_{l}<\ta_{l-1}\le\ta_t\le a_t$, which contradicts
$\im_{a_t}=\im_{b^+}=\im_a$.
\end{proof}

\begin{corollary} \label{cor: membership}
Let $\frakC=\seq{\ci_k}_{1 \le k \le l}$ be an admissible chain of
$i$-boxes. Let $\ci=[a,b]$ be an $i$-box such that
$a^-<\ta_l\le a$ and $b \le \tb_l<b^+$, i.e. $a=(\ta_l)(\im)^+$ and
$b=(\tb_l)(\im)^-$ for some $\im\in\Iff $. Then $\ci$ is a member of
$\frakC$.
\end{corollary}

\begin{proof}
It follows from the fact that $\ci\subset\tc_l$ commutes with all
$\ci_k$ $(1\le k \le l)$.
\end{proof}

\medskip
For an admissible chain $\frakC$ of $i$-boxes, let us define a
family of real simple modules $\cM[](\frakC)$ as follows:
$$\cM[](\frakC) \seteq \st{ M(\ci_k)}_{1 \le k \le l}.$$

\begin{theorem} \label{thm: admissible chain commuting}
For any admissible chain $\frakC=\seq{\ci_k}_{1\le k\le l}$, $\cM[](\frakC)$
forms a commuting family of real simple modules.
\end{theorem}

\begin{proof}
By Lemma~\ref{lem:chain}, $\ci_{k_1}$ and $\ci_{k_2}$ commute for
all $k_1, k_2 \in[1,l]$.  Thus our assertion follows from
Theorem~\ref{thm: commuting ab}.
\end{proof}

\subsection{Box moves}\label{sec:box moves}
For an admissible chain $\frakC=\seq{\ci_k}_{1\le k\le l}$ with the
associated pair $(a, \mathfrak{T})$ and for $1 \le s <l$, we say
that an $i$-box $\ci_{s}$ is \emph{movable} if $s=1$ or $T_{s-1} \ne
T_s$ $(s\ge2)$. Note that the last condition (for $s\ge2$) is
equivalent to the condition $\ta_{s+1}=\ta_{s-1}-1$ and
$\tb_{s+1}=\tb_{s-1}+1$.

For a movable $\ci_s$ in $\frakC$, we define a new  admissible chain
$B_s(\frakC)$ whose associated pair $(a',\mathfrak{T}')$ is given
\begin{align*}
& {\rm (i)} \ \begin{cases}
a' = a \pm  1  &\text{if $s=1$ and $T_1 = \RR$ (resp.\ $\LL$),}  \\
a'=a & \text{if $s>1$,}
\end{cases} \allowdisplaybreaks \\
& {\rm (ii)} \ T'_{k} = T_k \ \  \text{ for } k \not\in \{ s-1,s\}
\quad \text{ and } \quad  \ T'_{k}\ne T_k \ \ \text{ for } k \in \{
s-1,s\}.
\end{align*}
That is, $B_s(\frakC)$ is the admissible chain obtained from
$\frakC$ by moving $\tc_s$ by $1$ to the right or to the left inside
$\tc_{s+1}$. We call $B_s(\frakC)$ the {\em box move} of $\frakC$ at
$s$.

\begin{proposition} \label{prop: changing order}
Let $\frakC=\seq{\ci_k}_{1 \le k \le l}$ be an admissible chain of
$i$-boxes and let $k_0$ be a movable $i$-box $(1\le k_0<l)$. Set
$B_{k_0}(\frakC) = \seq{\ci_k'}_{1 \le k \le l}$. Assume that
$\im_{\ta_{k_0+1}} \ne \im_{\tb_{k_0+1}}$, i.e., $\tc_{k_0+1}$ is
not an $i$-box. Then we have
$$\ci_k' = \ci_{\mathfrak{s}_{k_0}(k)}, $$
where $\mathfrak{s}_{k_0} \in \sym_l$ is the transposition of $k_0$
and $k_0+1$.
\end{proposition}

\begin{proof}
Set $p=k_0+1$.By the assumption, we have
$$
[\ta_p,\tb_p\} = [\ta_p,\tb_p-1\} \qtq\{\ta_p,\tb_p] =
\{\ta_p+1,\tb_p],
$$
which implies the desired result.
\end{proof}

\begin{proposition} \label{prop: real T-move}
Let $\frakC=\seq{\ci_k}_{1 \le k \le l}=(c,\mathfrak{T}) $ be an
admissible chain of $i$-boxes and let $\ci_{k_0}$ be a movable
$i$-box. Assume that $\im_{\ta_{k_0+1}} = \im_{\tb_{k_0+1}}$, i.e.,
$\tc_{k_0+1}$ is an $i$-box. Set $\ci_{k_0+1}=[a,b]$ with
$a=\ta_{k_0+1}$ and $b=\tb_{k_0+1}$, and
set $B_{k_0}(\frakC) = \seq{\ci'_k}_{1\le k\le l}$.
Then we have \bnum
\item $\ci_{k_0}=[a^+,b]$ and $\ci'_{k_0}=[a,b^-]$ if $T_{k_0-1}=\RR$,\label{it:1}
\item $\ci_{k_0}=[a,b^-]$ and $\ci'_{k_0}=[a^+,b]$ if $T_{k_0-1}=\LL$,\label{it:2}
\end{enumerate}
In particular, we have an exact sequence
\eq &&0 \to  \tens_{
d(\im_a,\jm)=1} M[a(\jm)^+,b(\jm)^-]  \to  X\tens Y
  \to   M(\ci_{k_0+1}) \tens M[a^+,b^-]   \to 0.
\label{eq:Tsystem2} \eneq where $(X,Y)=\bl
M(\ci_{k_0}),M(\ci'_{k_0})\br$ in case \eqref{it:1} and $(X,Y)=\bl
M(\ci'_{k_0}),M(\ci_{k_0})\br$ in case \eqref{it:2}.
\end{proposition}

\begin{proof}
By the assumption that $\im_{\ta_{k_0+1}} = \im_{\tb_{k_0+1}}$,
$[a,b] \seteq [a_{k_0+1},b_{k_0+1}]= [\ta_{k_0+1},\tb_{k_0+1}]$.

\noindent In case \eqref{it:1}, we have
$$[a_{k_0},b_{k_0}]= \{ a_{k_0+1}+1,\tb_{k_0+1}] =[a^+,b] \quad \text{ and }
\quad \ci'_{k_0}=[\ta_{k_0+1},\tb_{k_0+1}-1\}=[a,b^-].$$ The proof
in case \eqref{item:2} is similar.

\snoi The last statement follows from Theorem~\ref{th:Tsystem}.
\end{proof}

\begin{remark} \label{rmk: combinatorial T-system}
When $\ci_{k_0}$ is movable and $\tc_{k_0+1}=\ci_{k_0+1}=[a,b]$ is
an $i$-box, Corollary~\ref{cor: membership}
tells that the $i$-boxes $\{
[a(\jm)^+,b(\jm)^-] \ | \ d(\im_a,\jm)=1 \}$ appearing
in~\eqref{eq:Tsystem2} are all contained in $\frakC_{\le k_0+1}$ and
hence in $\frakC$.
 Note that $[a^+,b^-]$  commutes with $i$-boxes in  $\{[a^+,b],[a,b^-],[a,b]\}=\{ \ci_{k_0-1}, \ci_{k_0}, \ci_{k_0+1} \}$, and
it also commutes with all i-boxes in $\frakC_{\le k_0-2}$
since $(a^+)^- =\ta_{k_0+1} < \ta_k\le a_k\le b_k \le \tb_k < \tb_{k_0+1}=(b^{-})^+$ for $1\le k \le k_0-2$.
It follows that $[a^+,b^-]$ is contained in $\frakC$.
Thus, the operation $B_{k_0}$ under this
situation can be understood as a combinatorial analogue of
T-system (see Lemma~\ref{lem: b is mu} below,
 also \cite[Section 13]{GLS11} and \cite[Section 12]{GLS13S}).
\end{remark}

\begin{definition}
For admissible chains $\frakC^{(1)}$ and $\frakC^{(2)}$ of the same
length $l \in \Z_{\ge 1}$, we say that they are
\emph{$T$-equivalent}, denoted by $\frakC^{(1)} \overset{T}{\sim}
\frakC^{(2)}$, if there exists a sequence $(p_1,p_2,\ldots,p_r) \in
\{ 1,2,\ldots,  l-1 \}^r$ $(r \in \Z_{\ge 1})$ such that
$$   B_{p_r}(\cdots(B_{p_2}(B_{p_1}(\frakC^{(1)}))\cdots)  = \frakC^{(2)}.$$
\end{definition}

The following lemma is almost obvious. \Lemma\label{lem: finite
sequence are T-equi} \hfill \bnum
\item The binary relation $\overset{T}{\sim}$ for admissible chains of finite length is an equivalence relation.
\item If $\frakC^{(1)}$ and $\frakC^{(2)}$ are $T$-equivalent, then they have the same range.
\item Two admissible chains $\frakC^{(1)}$ and $\frakC^{(2)}$
with the same range are $T$-equivalent.
\end{enumerate}
\enlemma

\begin{example}
Let us consider a \pbw $(\Dd,\hhw)$ with  $\fing=A_3$ and
$$
\ba{rcccccccl}  \hhw = &\cdots &s_{\im_{-2}}\;&s_{\im_{-1}}\;&\,s_{\im_0}\;&\;s_{\im_1}&\;s_{\im_2}\;&\;s_{\im_3}&\cdots\\
 =&  \cdots &s_{3}&s_{2}&s_{3}&s_{1}&s_{2}&s_{3}&\cdots \ea
 $$
and an admissible chain $\frakC^{(1)}=\seq{\ci^{(1)}_k}_{1\le k\le
3}=(0,\mathfrak{T}=(\LL,\LL))=([0],[-1,0\},[-2,0\})=([0],[-1],[-2,0])$
associated to $\hhw$.

Note there exists only one movable box $\ci^{(1)}_1$ in
$\frakC^{(1)}$. Thus we have \eqn \frakC^{(2)}=\seq{\ci^{(2)}_k}_{1\le
k\le 3} \seteq B_1 (\frakC^{(1)}) &=& \bl-1,
\mathfrak{T}=(\RR,\LL)\br\\*
&=&\bl[-1],\{-1,0],[-2,0\}\br=\bl[-1],[0],[-2,0]\br,\eneqn since
$[-1,0]$ is not an $i$-box (see Proposition~\ref{prop: changing
order}).

For $\frakC^{(2)}$,  the second $i$ box $\ci^{(2)}_2$ is movable and
hence we have \eqn \frakC^{(3)}=\seq{\ci^{(3)}_k}_{1\le k\le 3} \seteq
B_1 (\frakC^{(2)}) &&= \bl-1, \mathfrak{T}=(\LL,\RR)\br\\*
&&=\bl[-1],[-2,-1\},\{-2,0] \br=\bl[-1],[-2],[-2,0]\br,\eneqn
since
$[-2,0]$ is an $i$-box (see Proposition~\ref{prop: real T-move}).

Finally, $\ci^{(3)}_1$ is movable hence we have \eqn && \frakC^{(4)}
= B_1 (\frakC^{(3)}) = \bl-2, \mathfrak{T}=(\RR,\RR)\br
=\bl[-2],\{-2,-1],\{-2,0]\br=\bl[-2],[-1],[-2,0]\br.\eneqn
\end{example}

\vskip 2em
\section{$\rmQ$-data and associated \pbws} \label{sec: subcategories}

In this section, we recall the notion of $\rmQ$-data
introduced in \cite{FO20}.  A $\rmQ$-datum
is a generalization of a Dynkin quiver with a height
function. Then we attach a complete \pbw to each $\rmQ$-datum.

\subsection{{\rm Q}-data}\label{subsec: Q-datum}
For each untwisted quantum affine algebra $U_q'(\g)$, we associate a
finite simple Lie algebra $\gf$ of simply-laced type as follows (see
\S\,\ref{subsec: E(M)} and  \cite{KKOP20A}):
\renewcommand{\arraystretch}{1.5}
\begin{align} \label{Table: root system}
\small
\begin{array}{|c||c|c|c|c|c|c|c|}
\hline
 \g  & A_n^{(1)} \ (n\ge1)  & B_n^{(1)} \ (n\ge2) & C_n^{(1)} \ (n\ge3)  & D_n^{(1)} \ (n\ge4) & E_{6,7,8}^{(1)} & F_{4}^{(1)} & G_{2}^{(1)}  \\ \hline
 \gf  & A_n & A_{2n-1}    & D_{n+1}   &  D_n & E_{6,\,7,\,8} & E_{6} & D_{4}  \\
\hline
\mathsf{h}^\vee&n+1&2n-1&n+1&2n-2&12,\,18,\,30&9&4\\
\hline
\end{array}
\end{align}

Here $\mathsf{h}^\vee\seteq\ang{c,\rho}=\sum_{i\in I}\mathsf{c}_i $ is
the dual Coxeter number of $\g$.

\vskip 0.7em Let $\Dynkin_\gf$ be the Dynkin diagram for $\gf$. Then
there exists a Dynkin diagram automorphism $\sigma$  of
$\Dynkin_\gf$ whose orbit set $\Dynkin_\gf^\sigma$ yields the Dynkin
diagram $\Dynkin_{\g_0}$: for $\g=(A DE)^{(1)}_n$-case, we have $\sigma={\rm id}$.
In the remaining cases,
$\sigma$ is given by $\vee$ or $\wvee$ in the following diagrams:
\begin{align*}
\hs{-1ex}\big( \Dynkin_{A_{2n-1}}: \xymatrix@R=0.5ex@C=3ex{ *{\circ}<3pt>
\ar@{-}[r]_<{1 \ \ }  &*{\circ}<3pt> \ar@{-}[r]_<{2 \ \ }  &   {}
\ar@{.}[r] & *{\circ}<3pt> \ar@{-}[r]_>{\,\,\,\  _{2n-2} }
&*{\circ}<3pt>\ar@{-}[r]_>{\,\,\,\,  _{2n-1} } &*{\circ}<3pt> }, \
k^{\vee} = 2n-k \big)\hs{-.5ex} &\Longrightarrow \Dynkin_{B_{n}}\hs{-.5ex}:\hs{-.5ex}
\xymatrix@R=0.5ex@C=3ex{ *{\circ}<3pt> \ar@{-}[r]_<{ 1 \ \ }
&*{\circ}<3pt> \ar@{-}[r]_<{2 \ \ }  &   {} \ar@{.}[r] &
*{\circ}<3pt> \ar@{}_>{ _{n-1}} &*{\circ}<3pt>\ar@{<=}[l]^>{ \qquad
\ \  \; \ _{n}}  },
\allowdisplaybreaks\\
  \hs{-.4ex} \left(\hs{-1ex}
\Dynkin_{D_{n+1}} \hspace{-.5ex} :
\hspace{-.5ex}\raisebox{1em}{\xymatrix@R=0.5ex@C=2.5ex{
& & &  *{\circ}<3pt>\ar@{-}[dl]^<{ \  _{n}} \\
*{\circ}<3pt> \ar@{-}[r]_<{1\ \ }  &*{\circ}<3pt>
\ar@{.}[r]_<{2 \ \ } & *{\circ}<3pt> \ar@{.}[l]^<{ \ \  _{n-1}}  \\
& & &   *{\circ}<3pt>\ar@{-}[ul]^<{\quad \ \  _{n+1}}   \\
}} \hspace{-1ex}   ,\hs{1ex}     k^{\vee} \hspace{-.5ex}  =
\hspace{-.5ex} \begin{cases} k & \text{ if } k \le n-1, \\ n+1 &
\text{ if } k = n, \\ n & \text{ if } k = n+1 \end{cases} \right)\hs{-.5ex}
&\Longrightarrow \Dynkin_{C_{n}}\hs{-.5ex}:\hs{-.5ex} \xymatrix@R=0.5ex@C=3ex{
*{\circ}<3pt> \ar@{-}[r]_<{1 \ \ }  &*{\circ}<3pt> \ar@{-}[r]_<{ 2 \
\ }  &   {}
\ar@{.}[r] & *{\circ}<3pt> \ar@{}_>{ _{n-1} } &*{\circ}<3pt>\ar@{=>}[l]^>{ \qquad \ \  \; \ _{ n } }  },
\allowdisplaybreaks\\
  \left(
\Dynkin_{E_{6}}: \raisebox{2em}{\xymatrix@R=3ex@C=3ex{
& & *{\circ}<3pt>\ar@{-}[d]_<{\quad \ \  2} \\
*{\circ}<3pt> \ar@{-}[r]_<{1 \ \ }  & *{\circ}<3pt> \ar@{-}[r]_<{3 \
\ }  & *{\circ}<3pt> \ar@{-}[r]_<{4 \ \ }  & *{\circ}<3pt>
\ar@{-}[r]_<{5 \ \ }  & *{\circ}<3pt> \ar@{-}[l]^<{ \ \ 6 }
}},\hs{.5ex} \begin{cases} 1^{\vee}=6, \ 6^{\vee}=1, \\ 3^{\vee}=5,
\ 5^{\vee}=3, \\ 4^{\vee}=4, \ 2^{\vee}=2 \end{cases}  \right)
&\Longrightarrow \Dynkin_{F_{4}}: \xymatrix@R=0.5ex@C=3ex{
*{\circ}<3pt> \ar@{-}[r]_<{ 1 \ \ }  &*{\circ}<3pt> \ar@{}[r]_<{2 \
\ }  &  *{\circ}<3pt> \ar@{<=}[l]^<{ \ \ 3 } &
*{\circ}<3pt>\ar@{-}[l]^>{ \qquad \ \  \; \ 4}  },
\allowdisplaybreaks\\
   \left(  \Dynkin_{D_{4}}:  \raisebox{1em}{
\xymatrix@R=0.5ex@C=3ex{
& &   *{\circ}<3pt>\ar@{-}[dl]^<{ \ 3} \\
*{\circ}<3pt> \ar@{-}[r]_<{1 \ \ }  &*{\circ}<3pt>
\ar@{-}[l]^<{2 \ \ }   \\
& &    *{\circ}<3pt>\ar@{-}[ul]^<{\quad \ \  4} \\
}}, \ \begin{cases} 1^{\wvee}=3, \ 3^{\wvee}=4, \ 4^{\wvee}=1, \\
2^{\wvee}=2 \end{cases} \right)  &\Longrightarrow\Dynkin_{G_{2}}:
\xymatrix@R=0.5ex@C=3ex{ *{\circ}<3pt> \ar@{-}[r]_<{ 1 \ \ }
&*{\circ}<3pt>
\ar@{<=}[l]^<{ \ \ 2 }  }.
\end{align*}
We denote by  $\If=\{1,2,\ldots,
|\If|\}$ the index set of $\Dynkin_\gf$. Let
$(\cmf,\wlf,\Pif,\wlf^\vee,\Pif^\vee)$ be the finite Cartan datum
associated to $\gf$. Let us denote by $\Wf$ the Weyl group, by $\{
\La_\im \ | \ \im \in \If \}$ the set of fundamental weights, by
$\Phi$ the set of roots, by $\wlf^+$ the set of dominant integral
weights, by $\al_\im$ the $\im$-th simple root, by $\rl$ the root
lattice of $\gf$, and by  $( \ , \ )$ the symmetric bilinear form on
$\rl$, all of which are associated to $\gf$. We also write $\Phi^+$
(resp.\ $\Phi^-$) for the set of positive (resp.\ negative) roots of
$\gf$
 and  $\rl^+$ (resp.\ $\rl^-$) for
the positive (resp.\ negative) root lattice of $\gf$. For $\be
=\sum_{\im \in \If}k_\im\al_\im \in \rl$, we set $|\be|\seteq
\sum_{\im \in \If}|k_\im|$.

Let $I_{\g_0}=\{ 1,2,\ldots, n \}$ be the index set of
$\Dynkin_{\g_0}$. Note that $\If =
I_{\g_0}$ when $\g=(ADE)_n^{(1)}$. To distinguish the index sets $\If$ and $I_{\g_0}$, we
use $\im$ for indices in $\If$ and $i$ for indices in $I_{\g_0}$.
We write $\pi\col \If\to I_{\g_0}$ for the
projection. Then the orbit of $\im$ is given as
$$\pi^{-1}\pi(\im)=\st{\sigma^m(\im)\mid m\in\Z}.$$
We denote by $\cmA=(a_{i,j})_{i,j\in I_{\g_0}}$ the Cartan matrix
associated to $\g_0$. Then we have
$$ a_{i,j}=-|\,\st{\jm\in \pi^{-1}(j)\mid d(\im,\jm)=1}\,|.$$
for $i\not=j\in I_{\g_0}$ and $\im\in\pi^{-1}(i)$. Set $\fd_i=|
\pi^{-1}(i)|$ and $\fD={\rm diag}(\fd_i \ | \ i \in I_{\g_0})$ the
diagonal matrix. Then, $\fD \cmA$ is symmetric.

For $\im\in \If$, we write also $\fd_\im$ for $\fd_{\pi(\im)}$, i.e.,
$$\fd_\im=|\pi^{-1}\pi(\im)|.$$
Note that $\fd_\im\in\st{1,\ord(\sigma)}$ and
$$(\ual_{\pi(\im)},\ual_{\pi(\im)})=2\fdi/\ord(\sigma).$$

\smallskip

In this subsection, we fix $\g$ and a pair $(\Dynkin_\gf, \sigma)$
of the Dynkin diagram $\Dynkin_\gf$ of finite type $\gf$ and an
automorphism $\sigma = {\rm id}$, $\vee$ or $\wvee$. If there is no
afraid of confusion, we simply write $\Dynkin$ for $\Dynkin_\gf$.

\begin{definition} [{\cite[Definition 2.5]{FO20}}]
\label{def: height} A function $\xi \colon \If \to \Z$ is called a
\emph{height function on $(\Dynkin, \sigma)$} if the following two
conditions are satisfied. \ben
\item \label{it:ht1}
For any $\im, \jm \in \If$ such that $d(\im,\jm)=1$ and $\fdii=
\fdjj$, we have $|\xi_{\im} - \xi_{\jm}| =\fdii$.
\item \label{it:ht2}
For any $\im,\jm \in \If$ such that $d(\im,\jm)=1$ and $1=\fdii <
\fdjj=\ord(\sigma)$, there exists a unique element $j^{\circ} \in
\pi^{-1}\pi(\jm)$ such that $|\xi_{\im} - \xi_{j^{\circ}}| = 1$ and
$\xi_{\sigma^{k}(j^{\circ})} = \xi_{j^{\circ}} + 2k$
for any $0 \le k < \ord(\sigma)$.
\end{enumerate}
We call the triple $\Qd = (\Dynkin, \sigma, \xi)$ a
\emph{$\rmQ$-datum} for $\g$.

\end{definition}

Note that, when $\sigma={\rm id}$, a $\rmQ$-datum coincides with a usual Dynkin quiver with a height function (\cite{KT10}).

\begin{remark}
The convention for a $\rmQ$-datum in this paper is associated to the
\emph {sink-adapted orientation} which is different from the
\emph{source-adapted orientation} in \cite{FO20,HL11}.
We take this convention in order to match
$M[a,b]$ as the image of a certain unipotent quantum minor under the
quantum affine \SW duality functor (see
\S\,\ref{subsec: QAWS} and \S\,\ref{subsec:subcategory} below).
\end{remark}

\begin{example} For an untwisted affine type $\g$, we give examples of
$\rmQ$-data $\Qd$:
\begin{enumerate}
\item $\Qd=\xymatrix@C=3ex@R=3ex{
  *{\circ}<3pt> \ar[r]^<{\ble{_{n-1}}}_<{1} & *{ \circ }<3pt> \ar[r]^<{\ble{_{n-2}}}_<{2}   & \cdots \ar[r]_<{ }   &*{\circ}<3pt> \ar[r]^<{ _{\ble{1}} }_<{ _{n-1} } &*{ \circ }<3pt>  \ar@{-}[l]_<{ \ \ _{\ble{0}} }^<{ _{n} }}$
\hs{5ex}for $\g=A_n^{(1)}$, $\ord(\sigma)=1$,
\item
$\Qd=\raisebox{2.3em}{\xymatrix@C=3ex@R=3ex{
 &&& *{ \circ }<3pt> & \\
  *{\circ}<3pt> \ar[r]^<{\ble{_{n-2}}  }_<{1} & *{ \circ }<3pt> \ar[r]^<{\ble{_{n-3}}}_<{2}   & \cdots \ar[r]^<{ }   &*{\circ}<3pt> \ar[r]^<{ \ \ble{1}}_{ \qquad _{n} } \ar[u]^>{ {\ble{0}} }_>{ _{n-1} } &*{ \circ }<3pt>
  \ar@{-}[l]_<{ \ble{0} }^>{ _{n-2}}   }}$  \hs{5ex} for $\g=D_n^{(1)}$, $\ord(\sigma)=1$,
\item $\Qd= \hspace{-2.5ex} \xymatrix@C=4ex@R=3ex{ *{ \circ }<3pt>
\ar@{->}[r]^{ _{\ble{2n-3}} \qquad \quad}_<{1}  &*{\circ}<3pt>
\ar@{->}[r]^{  _{\ble{2n-5}} \qquad \quad }_<{2}  &   {} \ar@{.}[r]&
*{\circ}<3pt> \ar@{->}[r]^{ \ble{1} \qquad \ \  }_{ _{n-1} \quad \ \
\ \ \ }&  *{\circ}<3pt> \ar@{->}[r]^{ \ble{0} \qquad \quad }_<{ _{n}
\ \ } & *{\circ}<3pt> \ar@{-}[l]_{  \qquad \ \ \ble{-1}  }^{ \ \ \
\quad\ \ _{n+1} }  & *{\circ}<3pt> \ar@{->}[l]_{  \qquad \ \ \ble{1}
}^<{_{n+2}}
  {} \ar@{.}[r] &  *{\circ}<3pt> \ar@{<-}[r]^<{ _{\ble{2n-7}} \ \ }_<{ _{2n-2} }&  *{\circ}<3pt> \ar@{-}[l]_<{ \ \ _{\ble{2n-5}} \ \ }^<{ _{2n-1} } } \hs{-1ex}$ for $\g=B_n^{(1)}, \; \ord(\sigma)=2$,
\item $\Qd =
\raisebox{2.3em}{\xymatrix@C=3ex@R=3ex{  &&&& *{ \circ }<3pt>  \ar@{->}[d]_<{\ble{0}}^<{ _{n} } \\
*{ \circ }<3pt> \ar@{->}[r]^{ _{\ble{n-3}} \qquad \quad}_<{1 \ \ }
&*{\circ}<3pt> \ar@{->}[r]^{  _{\ble{n-4}} \qquad \quad }_<{2  \ \
}  &   {} \ar@{.}[r]&  *{\circ}<3pt> \ar@{->}[r]^{ \ble{0} \qquad
\quad }_<{_{n-2}} &  *{\circ}<3pt> \ar@{<-}[r]^{ \ble{-1} \qquad
\quad \ }_<{ _{n-1} \ \ } & *{\circ}<3pt> \ar@{-}[l]_{  \qquad \ \
\ble{2} }^<{  \ \  _{n+1} }  }}$ for $\g=C_n^{(1)}$,
$\ord(\sigma)=2$,
\item $
\Qd = \raisebox{2.3em}{\xymatrix@C=4ex@R=3ex{ && *{\circ}<3pt>\ar@{->}[d]_<{\ble{-2}}^<{2} \\
*{ \circ }<3pt> \ar@{->}[r]^<{\ble{0}  \ \; }_<{1}  &*{\circ}<3pt>
\ar@{->}[r]^<{\ble{-2} \ \; }_<{3} &*{ \circ }<3pt> \ar@{->}[r]^<{
\; \ble{-3}   }_<{4} &*{\circ}<3pt> \ar@{-}[r]^<{\ble{-4}  \ \;
}_<{5} &*{\circ}<3pt> \ar@{->}[l]_<{\  \ble{-2}}^<{6} }} $\hs{2ex}  for
$\g=F_4^{(1)}$, $\ord(\sigma)=2$,
\item
$\Qd = \raisebox{2.3em}{\xymatrix@C=4ex@R=3ex{ & *{ \circ }<3pt> \ar@{->}[d]_<{\ble{1}}^<{3}  \\
*{ \circ }<3pt> \ar@{<-}[r]^<{\ble{-1} \ \ }_<{1}   &*{\circ}<3pt>
\ar@{<-}[r]^<{\ble{0}  \ \ \ \  }_<{2}   &*{ \circ }<3pt>
\ar@{-}^<{\ble{3} \ }_<{4} } }$  \hs{2ex}for $\g=G_2^{(1)}$,
$\ord(\sigma)=3$.
\end{enumerate}
Here, \ben[(1)]
\item an underline integer $\ble{*}$
is the value of $\xi_\im$ at each vertex $\im \in \Dynkin$,
\item an arrow $\im \to \jm$ means that $\xi_\im>\xi_\jm$ and $d(\im,\jm)=1$.
\ee
\end{example}

For a $\rmQ$-datum $\Qd$, we call a vertex $\im \in \If$
 a \emph{sink} of $\Qd$ if $  \xi_\im < \xi_\jm$ for all $\jm\in\If$ such that $d(\im,\jm)=1$.
We also call a vertex a \emph{source} of $\Qd$ if $  \xi_\im
-2\fdi > \xi_\jm -2\fdj$ for all $\jm\in\If$ such that
$d(\im,\jm)=1$.

For a $\rmQ$-datum  $\Qd=(\Dynkin,\sigma,\xi)$ and its sink $\im$,
we denote by $s_\im \Qd$ the $\rmQ$-datum
$(\Dynkin,\sigma,s_\im \xi)$ where $s_\im \xi$ is  the height
function defined as follows (\cite[Lemma 2.11]{FO20}):
\begin{align*}
(s_\im\xi)_\jm =  \xi_\jm+ \delta_{\im\jm} \times 2\fdi,
\end{align*}
Then $s_\im \Qd$ becomes a $\rmQ$-datum associated to the same $\g$
of $\Qd$.  Similarly, we can define the $\rmQ$-datum
$s_\im^{-1}\Qd$ for a source $\im \in \If$ of $\Qd$.  We call
these operations as (combinatorial) \emph{reflection functors} on
$\rmQ$-data associated to $\g$.

\medskip

For a reduced expression $\tw=s_{\im_1}s_{\im_2}\cdots s_{\im_l}$ of
$w \in  \Wf$  and $1 \le k \le l$, we set \eq &&\hs{2ex}  \tw_{\le
k} \seteq s_{\im_1}s_{\im_2}\cdots s_{\im_k}, \quad\tw_{< k} \seteq
s_{\im_1}s_{\im_2}\cdots s_{\im_{k-1}} \quad \text{ and } \quad
\be^{\tw}_k  \seteq s_{\im_1}\cdots s_{\im_{k-1}} \al_{\im_k}. \eneq

Let $\Qd$ be a $\rmQ$-datum associated to $\g$ and let $\Wf$ be the
Weyl group of type $\Dynkin_\gf$. For $w\in \Wf$ and its reduced
expression $\tw=s_{\im_1}s_{\im_2}\cdots s_{\im_l}$, we say that
$\tw$ is \emph{adapted to $\Qd$} (or simply {\em $\Qd$-adapted}\/)
if
$$  \text{ $\im_k$ is a sink of $s_{\im_{k-1}}s_{\im_{k-2}}\cdots s_{\im_{1}}\Qd$ for all $1 \le k \le l$}.$$

\begin{proposition} [{\cite[Corollary 2.21]{FO20}}] \label{prop: adapted characterization}
Let $w_0=s_{\im_1}\cdots s_{\im_\ell}$ be a reduced expression of
$w_0$. Then the following  conditions~\eqref{it: c1} and~\eqref{it:
c2} are equivalent$\col$ \bnum
\item \label{it: c1} The following two conditions hold$\colon$
\bna
\item For any $\im$, $\jm\in \If$ such that $d(\im,\jm)=1$
     and any $s$ such that $1\le s < s^+\le \ell$ and $\im=\im_s$,
     we have
$$
-a_{\pi(\jm),\pi(\im)} = \begin{cases}
|\{t\mid s<t<s^+,   \pi(\jm)=\pi(\im_t) \}| & \text{ if $\fdi < \fdj$,} \\[1ex]
|\{t\mid s<t<s^+,   \jm=\im_t \}| & \text{ otherwise}.
\end{cases}
$$
Here $s^+\seteq\min\bl\{p \mid s<p\le\ell ,\; \im_p=\im_s
\}\cup\st{+\infty}\br$ and $s^-\seteq\max\bl\{p \mid s>p\ge1 ,\;
\im_p=\im_s \}\cup\st{-\infty}\br$.
\item
 If $s^-<t<s<s^+\le\ell$ and $d(\im_s,\im_t)=1$,
$\fd_{\im_s}<\fd_{\im_t}$, then there exists $t'$ such that
$s<t'<s^+$ and $\im_{t'}=\sigma(\im_t)$. \ee
\item \label{it: c2} The expression $s_{\im_1}\cdots s_{\im_\ell}$ is adapted to some $\rmQ$-datum $\Qd$ associated to $\g$.
\end{enumerate}
\end{proposition}

In \cite{FO20}, for each $\rmQ$-datum $\Qd$, a unique element
$\tau_\Qd = s_{\im_1}\cdots s_{\im_{r}}\sigma \in\Wf \sigma$ is
defined, which can be understood as a generalization of a Coxeter
element associated to a Dynkin quiver  in the  $\sigma={\rm id}$
case.
Here, we regard $\Wf \sigma$ as the subset of the automorphism group
of the root lattice of $\Dynkin$. We call $\tau_\Qd$ the  {\em
$\Qd$-adapted Coxeter element}.
We refer the reader to \cite{FO20} for its definition, but give some of its properties
instead.

\Rem[{ \cite[\S 2]{FO20}}] \label{rem:qdata}
\hfill \bna
\item The element $(\tau_\Qd)^{\ord(\sigma)}$ is contained in $\Wf$, of length $\mathsf{t} \seteq \ord(\sigma) \times \vert I_{\g_0}\vert$ and
has a $\Qd$-adapted reduced expression $s_{\im_1} s_{\im_2} \cdots
s_{\im_{\mathsf{t}}}$.
\item
Any  $\Qd$-adapted reduced expression $s_{\im_1} s_{\im_2} \cdots
s_{\im_{\mathsf{t}}}$ of $(\tau_\Qd)^{\ord(\sigma)}$ satisfies the
following property. Set $\Qd' \seteq  s_{\im_{\mathsf{t}}} \cdots
s_{\im_1}\Qd$, and let $\xi$ and $\xi'$ be the height function of
$\Qd$ and $\Qd'$. Then we have
$$\text{$\xi_\im'=\xi_\im+2 \times \ord(\sigma)$ for any $\im\in \If$.}$$
\item Let $\dC$ be the dual Coxeter number of $\g_0$. Then one can check (see \cite[Table 1]{FO20}) that
$$ \mathsf{h}^\vee =\dfrac{2\, |\Phi^+|}{\ord(\sigma)\, |I_{\g_0}|}.$$
\item \label{item:ell}
Let $\tw_0=s_{\im_1}\cdots s_{\im_\N}$ be a $\Qd$-adapted reduced expression of the longest element $w_0$.
Set $\Qd'=\tw_0^{-1}\Qd\seteq s_{\im_\N}\cdots s_{\im_1}\Qd$ and let
$\xi$ and $\xi'$ be the height function of $\Qd$ and $\Qd'$. Then we
have
$$\text{$\xi_\im'= \xi_{\im^*}+\ord(\sigma)\dC$ for any $\im\in \If$.}$$
In particular, setting $\im_{\ell+1}=(\im_1)^*$,
$s_{\im_2}\cdots s_{\im_{\ell}}s_{\im_{\ell+1}}$ is an $(s_{\im_1}\Qd)$-adapted reduced expression of $w_0$.

\ee \enrem

For a $\rmQ$-datum $\Qd=(\Dynkin,\sigma,\xi)$ associated to $\g$, we
define
\begin{align}\label{eq: I_Q}
\hI_\Qd \seteq \{ (\im,p) \in \If \times \Z \ | \ p -\xi_\im \in
2\fdi\Z\}.
\end{align}

We define the quiver $\Psi_{\Qd}$ whose set of vertices is
$\widehat{I}_\Qd$ and arrows are assigned in the following way: for
$(\im,p),(\jm,q) \in \widehat{I}_{\Qd}$, we have
\begin{align}\label{eq: Psi_g}
(\im,p) \to (\jm,q)  \quad \text{ if } \ \  d(\im,\jm)=1 \ \ \text{
and } \ \ q-p=\min\{ \fdi,\fdj \}.
\end{align}

Let $\widehat{\Phi} \seteq \Phi^+ \times \Z$. For each $\im \in
\If$, we define
$$ \ga^\Qd_\im \seteq (1-\tau_\Qd^{\fdi})\La_\im \in \Phi^+.$$
In \cite{FO20,HL11}, it is shown that there exists a unique
bijection $\phi_\Qd\col\hI_\Qd  \to \widehat{\Phi}$ defined
inductively as follows:
\begin{eqnarray}&&
\parbox{70ex}{
\ben
\item $\phi_\Qd(\im,\xi_\im)=(\ga^\Qd_\im,0)$
\item if $\phi_\Qd(\im,p)=(\be,m)$, then we define
\bna
\item $\phi_\Qd(\im,p+2\fdi) = (\tau_\Qd^{\fdi} (\be),m) \qquad\qquad\quad$ if $ \tau_\Qd^{\fdi}(\be) \in \Phi^+$,
\item $\phi_\Qd(\im,p+2\fdi) = (-\tau_\Qd^{\fdi}(\be),m+1)\qquad \ \ $ if $\tau_\Qd^{\fdi}(\be) \in \Phi^-$,
\item $\phi_\Qd(\im,p-2\fdi) = (\tau_\Qd^{-\fdi}(\be),m)\qquad\qquad \ \ $ if $\tau_\Qd^{-\fdi}(\be) \in \Phi^+$,
\item $\phi_\Qd(\im,p-2\fdi) = (-\tau_\Qd^{-\fdi}(\be),m-1)\qquad    $ if $\tau_\Qd^{-\fdi}(\be) \in \Phi^-$.
\end{enumerate}
\end{enumerate}
}\label{eq: twisted bijection}
\end{eqnarray}

\begin{figure}
\begin{gather*}
\begin{split}
{\tiny   A^{(1)}_{3}}&=
\raisebox{2.3em}{\scalebox{0.55}{\xymatrix@C=3ex@R=3ex{
\cdots& (1,-4) \ar[dr] && (1,-2)  \ar[dr] && (1,0)  \ar[dr] && (1,2)  \ar[dr] && \cdots\\
\cdots &&  \ar[ur] \ar[dr] (2,-3) &&  \ar[ur] \ar[dr] (2,-1) &&  \ar[ur] \ar[dr] (2,1) &&  \ar[ur] \ar[dr] (2,3)\\
\cdots& (3,-4) \ar[ur] && (3,-2)  \ar[ur] && (3,0)  \ar[ur] && (3,2)
\ar[ur] &&  \cdots
}}} \allowdisplaybreaks\\
{\tiny B^{(1)}_{2}} &=
\raisebox{2.3em}{\scalebox{0.55}{\xymatrix@C=1.5ex@R=3ex{
\cdots&& \ar[dr] (1,-7) &&&& \ar[dr] (1,-3) &&&& \ar[dr] (1,1) & \cdots \\
\cdots& \ar[ur] (2,-8) && \ar[dr] (2,-6) && \ar[ur] (2,-4) && \ar[dr] (2,-2) &&  (2,0)\ar[ur] &&   \cdots \\
\cdots &&&& \ar[ur] (3,-5) &&&& \ar[ur] (3,-1)&&& \cdots
}}}\allowdisplaybreaks\\
{\tiny C^{(1)}_{3}} &=
\raisebox{2.3em}{\scalebox{0.55}{\xymatrix@C=1.5ex@R=3ex{
\cdots&(1,-6)\ar[dr] &&  \ar[dr] (1,-4) && \ar[dr] (1,-2) && \ar[dr] (1,0) && \ar[dr] (1,2) &&    \cdots \\
\cdots&&  \ar[ur]\ar[dr] (2,-5) && \ar[ur]\ar[ddr] (2,-3) && \ar[ur]\ar[dr] (2,-1) &&\ar[ur]\ar[ddr] (2,1) &&\ar[ur]\ar[dr] (2,3) &   \cdots \\
\cdots&&& (3,-4)\ar[ur] &&&& (3,0)\ar[ur] &&&& \cdots \ar[ul]\\
\cdots&(4,-6) \ar[uur]   &&&& (4,-2)\ar[uur] &&&& (4,2)\ar[uur] &&
\cdots
}}} \allowdisplaybreaks \\
{\tiny D^{(1)}_4} &=
\raisebox{2.3em}{\scalebox{0.55}{\xymatrix@C=3ex@R=3ex{
\cdots& (1,-6) \ar[dr] && (1,-4)  \ar[dr] && (1,-2)  \ar[dr] && (1,0)  \ar[dr] &&  \cdots\\
\cdots &&  \ar[ur] \ar[dr]\ar[ddr] (2,-5) &&  \ar[ur] \ar[ddr]\ar[dr] (2,-3) &&  \ar[ur] \ar[ddr]\ar[dr] (2,-1) &&  \ar[ur] \ar[ddr]\ar[dr] (2,1)\\
\cdots& (3,-6) \ar[ur] && (3,-4)  \ar[ur] && (3,-2)  \ar[ur] && (3,0)  \ar[ur] &&   \cdots \\
\cdots& (4,-6) \ar[uur] && (4,-4)  \ar[uur] && (4,-2)  \ar[uur] && (4,0)  \ar[uur] &&  \cdots \\
}}}
\end{split}
\end{gather*}
\caption{Some examples of the quivers $\Psi_\g$.} \label{fig:quiver
Psi_examples}
\end{figure}
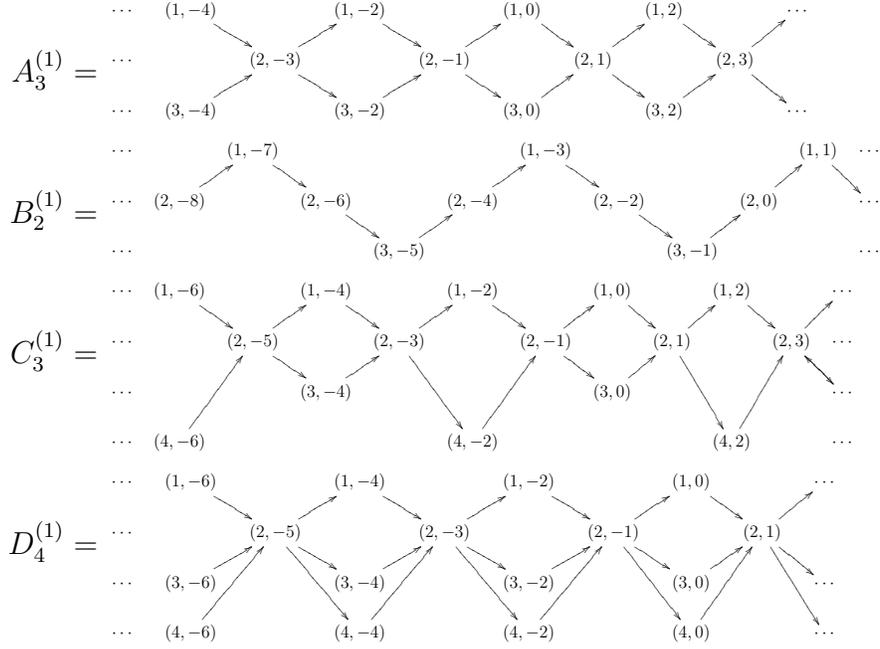

Note that the full subquiver $\Gamma_\Qd$ of $\hI_\Qd$ with
$$I_\Qd\seteq\phi_\Qd^{-1}(\Phi^+\times\{0\})\subset\If\times\Z$$
as the set of vertices
is isomorphic to the Auslander-Reiten (AR) quiver if $\Qd$ is associated to $\g=(ADE)_n^{(1)}$,
to the twisted AR-quiver if $\Qd$ is associated to
$\g=(BCF)_n^{(1)}$, and to the triply-twisted AR-quiver if $\Qd$ is
associated to $\g=G_2^{(1)}$ (\cite{OhSuh19}). Throughout this
paper, we call $\Gamma_\Qd$ the AR-quiver in a uniform way.

\begin{proposition}[{\cite[Corollary 2.39]{FO20}}] \label{prop: dual}
Let $\mathsf{h}^\vee$ be the dual Coxeter number of $\g_0$. Then we
have
$$  \phi_\Qd(\im^*,\xi_\im+\ord(\sigma)\mathsf{h}^\vee) =(\ga^{\Qd}_{\im },1). $$
Hence, for $k \in \Z$ and $\phi_\Qd(\im,p)=(\be,0) \text{ for } \be
\in \Phi^+$, we have
$$  \phi_\Qd(\im^{k*},p+k\;\ord(\sigma)\mathsf{h}^\vee) =
\bl(-w_0)^k\beta,k\br
$$
Here $\im^{k*}\seteq
\im^{\overbrace{*\cdots*}^{\text{$|k|$-times}}}$. In particular, we
have
$$I_\Qd=\st{(\im,p)\in\hI_\Qd\mid
\xi_\im\le p<\xi_{\im^*}+\ord(\sigma)\mathsf{h}^\vee}.$$
\end{proposition}

\Prop [ \cite{B99,FO20,OhSuh19}] \label{prop: reading redez}
For a $\rmQ$-datum $\Qd$, we have the followings$\col$ \bnum
\item \label{it: sequnce Q} If a sequence $\seq{(i_k,p_k)}_{1\le k\le \ell}$ of elements of $I_\Qd$
satisfies \bna
\item \label{it: IQ} $I_\Qd=\st{(\im_k,p_k)\mid 1\le k\le \ell}$, and
\item \label{it: arrow} $(\im_k,p_k)\to(\im_{k'},p_{k'})$ in $\Gamma_\Qd$ implies $k<k'$,
\ee then $\tw_0=s_{\im_1}\cdots s_{\im_\ell}$ is a $\Qd$-adapted reduced
expression of $w_0$, and
$$
p_k \seteq\xi_{\im_k}( s_{\im_{k-1}} \cdots s_{\im_1} \Qd).
$$
Here, we set $ \xi_{\jm} (\Qd)\seteq \xi_\jm$
for a $\rmQ$-datum $\Qd = (\Dynkin, \sigma, \xi) $.
\item
Conversely, if $\tw_0=s_{\im_1}\cdots s_{\im_\ell}$ is a $\Qd$-adapted reduced
expression  of $w_0$,  then $\seq{(i_k,p_k)}_{1\le k\le \ell}$ with
$p_k \seteq\xi_{\im_k}( s_{\im_{k-1}} \cdots s_{\im_1} \Qd)$
satisfies {\rm(a)} and {\rm(b)}
in~\eqref{it: sequnce Q}.

\item For any sequence $\seq{(i_k,p_k)}_{1\le k\le \ell}$ satisfying {\rm(a)} and {\rm(b)}
in~\eqref{it: sequnce Q}, we have $\phi_\Qd(\im_k,p_k) = (\be_k,0)$
$(1 \le k \le \ell)$, where $\be_k=s_{\im_1}\cdots
s_{\im_{k-1}}(\al_{\im_k})$.
\ee \enprop

\begin{remark} \label{rmk: Up to constant}
The statements below follow from \cite{HL11} and \cite{FO20}.

\noindent
(a) For $\rmQ$-data $\Qd$ and $\Qd'$ associated to the
same $\g$,
$\hI_\Qd$ and $\hI_{\Qd'}$ differ by constant integer
(see \cite{HL11} and \cite{FO20}); that is,
there exists a unique $k \in \{0, 1, \ldots, 2\ord(\sigma)-1\}$ such
that
$$ \hI_{\Qd'} = \{ (\im,p+k) \mid (\im,p) \in \hI_\Qd \}.$$

Note also that there exists a unique pair $(\epsilon,m) \in
\st{0,1}\times[0,\ord(\sigma)-1]$ such that
$$ \hI_{\Qd'} = \{ (\im,p+\epsilon)\mid(\sigma^m(\im),p) \in \hI_\Qd \}.$$

Thus we use the notation $\hI_\g$ and $\Psi_\g$ instead of $\hI_\Qd$
and $\Psi_\Qd$  when we can neglect integer shift.

\mnoi (b)  A  subset $\K$ of $\If\times\Z$ is equal to
$\hI_\Qd$ for some $\rmQ$-datum $\Qd$ if and only if $\K$ satisfies
the following condition:
\eq&&\hs{2.8ex}\left\{
\parbox{77ex}{
\bnum
\item For any $\im\in \If$, we have $\st{m\in\Z\mid (\im,m)\in \K}=2\fdii\,\Z+a$ for some $a\in\Z$,
\item
if $d(\im,\jm)=1$, $\fd_\im\le\fd_\jm$ and $(\jm,p)\in \K$, then we
have $(\im, p+\fdi)\in \K$,
\item \label{it: sig} for any $(\im,m)\in \K$, we have $\bl\sigma(\im),m+2\br\in \K$.
\ee
}\right.
\label{eq: hIQd} \eneq
\end{remark}

\begin{definition} \label{def: admissible sequence}
We say that an infinite sequence
$\hw\seteq\seq{(\im_k,p_k)}_{k\in\Z}$ in $\If\times\Z$ is {\em
admissible} if the following conditions are satisfied:
\ben 
\item \label{it:ssp} $p_{s^+}=p_s+2\fd_{\oi_s}$,
\item \label{it:tsts}
$p_t=p_s+\min(\fd_{\im_s},\fd_{\im_t})$ if $d(\im_s,\im_t)=1$
and  $t^-<s<t<s^+$,
\item \label{it:sigma+2}
 $p_t-(p_s+2)\in2\fd_{\im_s}\Z$ if $s,t\in\Z$ satisfy
$\im_t=\sigma(\im_s)$,
\item \label{it:ii} $s_{\im_k}\cdots s_{\im_{k+\ell-1}} = w_0$
for every $k \in \Z$, where $\ell$ denotes the length of the longest
element $w_0 \in \Wf$.
\end{enumerate}
\end{definition}

Note that \eqref{it:ii} implies that $\im_{k+\ell}=(\im_k)^*$ for any $k\in\Z$.
Note also that, when
$s<t$ and $\im_s\not=\im_t$,
the condition $t^-<s<t<s^+$ is equivalent to
$\im_s,\im_t\not\in\st{\im_k\mid s<k<t}$.

\begin{lemma} \label{lem:tendency}
Let $\hw=\seq{(\im_s,p_s)}_{s\in\Z}$ be an admissible sequence,
and assume that  $s,t\in\Z$ satisfy $d(\im_s, \im_t)=1$.
Then we have the followings$\col$
\bnum
\item \label{it: condi}
$p_s<p_t$ if and only if $s<t$.
\item \label{it: condii}
$p_s-p_t-a\in2a\Z$ where $a=\min(\fd_{\im_s}.\fd_{\im_t})$.

\ee
\end{lemma}

\begin{proof}
It is enough to show that
$p_s<p_t$ and $p_s-p_t-a\in2a\Z$
under the condition $s<t$.

Take $t'\in\Z$ such that $\im_{t'}=\im_t$ and
$(t')^-<s<t'$ and then take $s'\in\Z$ such that
$\im_{s'}=\im_s$ and $s'<t'<(s')^+$.
Then we have
$(t')^-<s'<t'<(s')^+$, and hence
$$p_{t'}=p_{s'}+a$$
by \eqref{it:tsts}.
Since $s<s'$ and $t'<t$, \eqref{it:ssp} implies $p_s\le p_{s'}$ and $p_{t'}\le p_t$.
Thus we obtain
$$p_s\le p_{s'}<p_{t'}\le p_t$$
and
$$p_s+a\equiv p_{s'}+a\equiv p_{t'}\equiv p_t \ \mod2a.$$
Note that $ 2\fd_{\im_s},2\fd_{\im_t} \in 2a\Z$.
\QED

\Prop\label{prop:admseq}
Let $A$ be the set of pairs $(\Qd,\tw_0)$ of a $\rmQ$-datum $\Qd$ and a $\Qd$-adapted reduced expression $\tw_0$ of $w_0$,
and let $B $ be the set of admissible sequences $\hw$.
Then there exists a bijective map $\vrho\col A\To B$ defined
as follows.
\bnum
\item \label{item:it}
For a pair $(\Qd,\tw_0)$ of a $\rmQ$-datum $\Qd = (\Dynkin, \sigma, \xi) $ and a $\Qd$-adapted reduced expression $\tw_0=s_{\jm_1}\cdots s_{\jm_\ell}$ of $w_0$,
we define $\seq{(\im_s,p_s)}_{s\in\Z}$ by:
\eqn
\im_{s+m\ell}&&=\bc\jm_s&\text{if $m$ is even,}\\
(\jm_s)^*&\text{if $m$ is odd,}\ec
\qt{for $s,m\in\Z$ such that $1\le s \le \ell$,}\\
p_k&& =\bc
\xi_{\im_k}\bl s_{\im_{k-1}} \cdots s_{\im_1} \Qd\br&\text{if $k\ge1$,}\\
\xi_{\im_k}\bl (s_{\im_{k}})^{-1}(s_{\im_{k+1}})^{-1}
\cdots (s_{\im_0})^{-1} \Qd\br
&\text{if $k\le 0$.}
\ec
\eneqn
Here, we set $ \xi_{\jm} (\Qd)\seteq \xi_\jm$
for a $\rmQ$-datum $\Qd = (\Dynkin, \sigma, \xi) $.\\
Then $\seq{(\im_s,p_s)}_{s\in\Z}$ is an admissible sequence and we set
$\vrho(\Qd,\tw_0)=\seq{(\im_s,p_s)}_{s\in\Z}$.

\item \label{it: converse}
Conversely, for an admissible sequence $\hw=\seq{(\im_s,p_s)}_{s\in\Z}$,
$(\Qd,\tw_0)=\vrho^{-1}(\hw)$
is given by
\eqn
\tw_0=s_{\im_1}\cdots s_{\im_\ell}\quad \text{ and } \quad  \xi_\im =p_{\uxi_\im},
\eneqn
where $\uxi_\im \seteq\min \st{  k \in \Z_{\ge 1} \mid  \im_k=\im }$.
\ee

Moreover, if $(\Qd,\tw_0)$ and  $\hw$ correspond by $\vrho$,
the following properties hold.
\bna
\item The map $s\mapsto (\im_s,p_s)$ gives a bijection $\Z\isoto\hI_\Qd$.
\label{item:1}
\item For $1\le k\le \ell$, $\phi_\Qd(\im_k,p_k)=(\beta_k,0)$,
where $\beta_k=s_{\im_1}\cdots s_{\im_{k-1}}\al_{\im_k}$.
\label{item:2}
\item $\im_{s+\ell}=(\im_s)^*$ and $p_{s+\ell}=p_s+\ord(\sigma)\;\mathsf{h}^\vee$,
where $\mathsf{h}^\vee$ is the dual Coxeter number of $\g_0$.
\label{item:3}

\item Set
$\tw_0{}'=s_{\im_2}s_{\im_3}\cdots s_{\im_\ell}s_{\im_{\ell+1}}$
and $\hw{}'\seteq\seq{(\im_k',p_k')}_{k\in \Z}$
where $(\im_k',p_k')=(\im_{k+1},p_{k+1})$.
Then $\im_1$ is a sink of $\Qd$, and $\hw{}'=\vrho\bl s_{\im_1}\Qd, \tw_0{}'\br$.
\label{item:5}

\item If $s^-<t<s$,  $d(\im_s,\im_t)=1$ and $\fd_{\im_s}<\fd_{\im_t}$,
then there exists $t'$ such that $s<t'<s^+$ and
$\im_{t'}=\sigma(\im_t)$.
\label{item:6}
\ee

\enprop

\Proof
\eqref{item:it} By the construction of $\seq{\im_k}_{k\in\Z}$,
it is evident that we have $w_0=s_{\im_{k+1}}\cdots s_{\im_{k+\ell}}$
for any $k$.
By Remark~\ref{rem:qdata},
$\im_k$ is a sink of $s_{\im_{k-1}} \cdots s_{\im_1} \Qd$ for any $k\ge1$.
Similarly $\im_k$ is a source of
$(s_{\im_{k}})^{-1}(s_{\im_{k+1}})^{-1}\cdots (s_{\im_0})^{-1} \Qd$.
Hence, there exists a sequence of Q-data $\seq{\Qd_k}_{k\in\Z}$ such that
\eqn
&&\Qd_0=\Qd,\\
&&\text{$\im_{k+1}$ is a sink of $\Qd_k$ and $\Qd_{k+1}=s_{\im_{k+1}}\Qd_k$,}\\
&&\text{$\im_{k}$ is a source of $\Qd_k$ and $\Qd_{k-1}=(s_{\im_{k}})^{-1}\Qd_k$,}
\eneqn
Hence, $p_k=\xi_{\im_k}(\Qd_{k-1})$ and  $(\im_k,p_k)\in\hI_\Qd$ for any $k$.

Now, the conditions \eqref{it:ssp}, \eqref{it:sigma+2} and
\eqref{it:ii}
in Definition~\ref{def: admissible sequence}
are almost obvious by the construction (see also \eqref{eq: hIQd}).

Let us show \eqref{it:tsts}.
The following properties are almost obvious by the construction:
\eq
&&\hs{-5ex}\ba{rcl}
{\rm(\alpha)}&\akew[.5ex]&\text{$p_s<p_t$ if $d(\im_s,\im_t)=1$ and $s<t$.
}\\[.5ex]
{\rm(\beta)}&&\text{$p_t-p_s-\min(\fd_{\im_s},\fd_{\im_t})\in 2\min(\fd_{\im_s},\fd_{\im_t})\Z$
 if $d(\im_s,\im_t)=1$.\akew[3ex]}
\ea\label{eq:alb}
\eneq
Now assume $d(\im_s,\im_t)=1$ and $t^-<s<t<s^+$.
Then \eqref{it:ssp} and ($\alpha$) implies
$p_{t^-}=p_t-2\fd_{\im_t}<p_s<p_t<p_{s^+}=p_s+2\fd_{\im_s}$.
Hence we obtain
\eqn
 p_s<p_t<p_s+2\min(\fd_{\im_s},\fd_{\im_t}).
\eneqn
Then ($\beta$) implies that $p_t=p_s+\min(\fd_{\im_s},\fd_{\im_t})$.
Thus $\hw$ is an admissible sequence.

\mnoi
\eqref{it: converse} Conversely, let  $\hw=\seq{(\im_s,p_s)}_{s\in\Z}$ be an admissible sequence.
Let us show that $\seq{\xi_\im}_{\im\in\If}$ is a height function.

\snoi
(ii-a)\ Let us show \eqref{it:ht1} in Definition~\ref{def: height}.
Assume that $\im,\jm\in\If$ satisfy $d(\im,\jm)=1$ and $\fd_\im=\fd_\jm$.
In order to see $|\xi_\im-\xi_\jm|=\fd_\im$, we may assume that $\uxi_\im<\uxi_\jm$.
Then we have
$(\uxi_{\jm})^-\le 0<\uxi_{\im}<\uxi_{\jm}$.
Hence $p_{(\uxi_{\jm})^-}=\xi_\jm-2\fd_\jm<p_{\uxi_{\im}}=\xi_\im<p_{\uxi_{\jm}}=\xi_{\jm}$ by
Lemma~\ref{lem:tendency}\;\eqref{it: condi}.
Since $\xi_\im-(\xi_\jm-\fd_\im)\in 2\fd_{\im}\Z$
by Lemma~\ref{lem:tendency}\;\eqref{it: condii}.
we obtain $\xi_\im=\xi_\jm-\fd_\im$.

\snoi
(ii-b)\ Let us show (ii) in Definition~\ref{def: height}.
Assume that $\im,\jm\in\If$ satisfy $d(\im,\jm)=1$ and $\fd_\im=1<\fd_\jm$.
We may assume that $\uxi_{\jm'}\ge\uxi_{\jm}$ for any $\jm'\in\pi^{-1}\pi(\jm)$.
Set $s_0=\uxi_{\jm}$. Then by \eqref{it:ssp} and \eqref{it:sigma+2}
in Definition~\ref{def: admissible sequence},
 for any $k\in\Z$, there exists
$s_k\in\Z$ such that $\im_{s_k}=\sigma^k(\jm)$ and $p_{s_k}=p_{s_0}+2k$.
Since $p_{\uxi_i}\equiv p_{s_0}+1 \mod2$, there exists
$t_k$ such that $\im_{t_k}=\im$ and $p_{t_k}=p_{s_0}+2k+1$.
Then we have
$p_{s_k}<p_{t_k}<p_{s_{k+1}}$, which implies
$s_k<t_k<s_{k+1}$ by
Lemma~\ref{lem:tendency}\;\eqref{it: condi}.
Hence $\seq{s_k}_{k\in\Z}$ is a strictly increasing sequence.
Since $\im_{s_{-1}}=\jm$ and $s_{-1}<s_0=\uxi_\jm$, we have $s_{-1}<0$.
Hence
$s_{-1}\le 0<\uxi_\im$ and $(\uxi_{\im})^-\le 0<s_0$.
Then Lemma~\ref{lem:tendency}\;\eqref{it: condi} implies
$\xi_{\jm}-2<\xi_{\im}$ and $\xi_{\im}-2<\xi_\jm$.
Hence we have
$|\xi_\im-\xi_\jm|\le1$.
Since $\xi_\jm\not=\xi_\im$ by Lemma~\ref{lem:tendency}\;\eqref{it: condii}, we obtain
$|\xi_\im-\xi_\jm|=1$.
On the other hand, if $1\le k\le \ord(\sigma)-1$, then we have
$(s_k)^-=s_{k-\ord(\sigma)}\le s_{-1}<0$, which implies
$s_k=\uxi_{\sigma^k(\jm)}$, and we obtain
$\xi_{\sigma^k(\jm)}=p_{s_k}=p_{s_0}+2k=\xi_\jm+2k$.

\mnoi
\eqref{item:1} is obvious.

\snoi
\eqref{item:2} follows from Lemma~\ref{prop: reading redez}.

\snoi
\eqref{item:3} follows from the definition and
Remark~\ref{rem:qdata}\;\eqref{item:ell}.

\snoi
\eqref{item:5} follows from \eqref{item:it}.

\snoi
\eqref{item:6}
 By \eqref{it:ssp} and \eqref{it:sigma+2}
in Definition~\ref{def: admissible sequence},
there exists $t'$ such that $\im_{t'}=\sigma(\im_t)$ and $p_{t'}=p_t+2$.
Then we have
$p_s\le p_t+1<p_{t'}=p_t+2<p_s+2=p_{s^+}$, where
the first inequality follows from
$p_s-2=p_{s^-}<p_t$.
Hence $s<t'<s^+$ by  Lemma~\ref{lem:tendency}\;\eqref{it: condi}.
\QED

\subsection{Associated fundamental modules, and twisted case}
For each untwisted affine $\g$ and $(\im,p) \in \If\times\Z$, we
define the fundamental module $V_\g(\im,p)\in \Ca_\g$ in the
following way: set $\qm \seteq q^{1/\ord(\sigma)}$ and
\begin{align} \label{eq: V(i;p)}
V_\g(\im,p) \seteq \begin{cases}
V(\varpi_{\pi(\im)})_{(-\qm)^p} & \text{ if } \g=A^{(1)}_n, \; C_n^{(1)}, \; D^{(1)}_n, \; E_{6,7,8}^{(1)}, \; G_2^{(1)}, \\
V(\varpi_{\pi(\im)})_{(-1)^{d(\im,n)}(\qm)^p} & \text{ if } \g=B_n^{(1)}, \\
V(\varpi_{\pi(\im)})_{(-1)^{d(\im ,2)}(\qm)^p} & \text{ if }
\g=F_4^{(1)}.
\end{cases}
\end{align}
Note that $q_{\pi(\im)}=(\qm)^{\fdi}$. By \cite{HL11,OhSuh19},
$V_\g(\im,p)$'s are distinct. The fundamental modules
$\seq{V_\g(\im,p)}_{(\im,p)\in\If\times\Z}$ do not depend on the
choice of $\rmQ$-data.

\medskip
Now let us consider the twisted affine types
$\fing^{(2)}=A_{2n}^{(2)}$ $(n \ge 1)$, $A_{2n-1}^{(2)}$ $(n \ge
2)$, $D_{n+1}^{(2)}$ $(n \ge 3)$, $E_6^{(2)}$ and
$\fing^{(3)}=D_4^{(3)}$. By \cite[Theorem 4.15]{Her10}, there exists
a ring isomorphism
\begin{align}\label{eq: chi}
\chi_t\col K(\Ca^0_{\fing^{(1)}})\isoto  K(\Ca^0_{\fing^{(t)}})
\qquad (t=2,3)
\end{align}
sending fundamental modules to fundamental modules, and KR-modules
to KR-modules (see also \cite{KKKO15III}). Here,
$K(\Ca^0_{\fing^{(1)}})$ and $K(\Ca^0_{\fing^{(t)}})$ denote the
Grothendieck rings of $\Ca^0_{\fing^{(1)}}$ and
$\Ca^0_{\fing^{(t)}}$, which are called the Hernandez-Leclerc
subcategory of $U'_q(\fing^{(1)})$ and  $U'_q(\fing^{(t)})$,
respectively (see Subsection~\ref{subsec:subcategory} below for
definitions). Note that $\chi_t$ commutes with $\D$.

The image of fundamental modules by $\chi_t$ can be described as
follows.

For $(\im,p) \in \If\times\Z$, we set
\begin{align*}
V_{\fing^{(t)}}(\im,p)\seteq V(\varpi_{\pi(\im)})_{A(\im,p)},
\end{align*}
where
\begin{align*}%
&\pi(\im)\seteq \begin{cases} \im & \text{ if } {\rm (i)} \  \fing^{(2)}=A^{(2)}_N \text{ and } \im \le \lceil N/2 \rceil  \\ & \quad \text{ or } {\rm (ii)} \  \fing^{(2)}=D^{(2)}_{n+1}  \text{ and } \im <n,  \\[1ex]
N+1- \im  & \text{ if } \fing^{(2)}=A^{(2)}_N \text{ and } \im > \lceil N/2 \rceil, \\
  n  & \text{ if } \fing^{(2)}=D^{(2)}_{n+1} \text{ and }  \im=n,n+1, \\
\end{cases} \allowdisplaybreaks \\[1ex]
&\bc
 \pi(1)=\pi(6)=1,  \quad \pi(3)=\pi(5)= 2,  \quad \pi(4)=3,  \quad
\pi(2)=4 &\text{if $\fing^{(2)}=E_6^{(2)}$},\\
\pi(1)=\pi(3)=\pi(4)= 1, \quad \pi(2) = 2 &\text{if
$\fing^{(3)}=D_4^{(3)}$,} \ec \end{align*} and
\begin{align*}
& A(\im,p)\seteq
\begin{cases}
 (-q)^{p} & \text{if } {\rm (i)} \  \fing^{(2)}=A^{(2)}_N \text{ and } \im \le  \lceil N/2 \rceil  \\
& \quad  \text{ or } {\rm (ii)} \  \fing^{(2)}=E^{(2)}_6 \text{ and } \im=1,3, \\
(-1)^N (-q)^{p} & \text{if } \fing^{(2)}=A^{(2)}_N \text{ and } \im >  \lceil N/2 \rceil,  \\
 (\sqrt{-1})^{n+1-\im} (-q)^{p} & \text{if } \fing^{(2)}=D^{(2)}_{n+1} \text{ and } \im <n, \\
(-1)^\im (-q)^{p} & \text{if } \fing^{(2)}=D^{(2)}_{n+1} \text{ and } \im =n,n+1, \\
- (-q)^p & \text{if }  \fing^{(2)}=E^{(2)}_6 \text{ and } \im=5,6, \\
 (\sqrt{-1}) (-q)^p & \text{if }  \fing^{(2)}=E^{(2)}_6 \text{ and } \im=2, 4,\\
\bigl(\delta_{\im,1}-\delta_{\im,2}+\delta_{\im,3} \omega +
\delta_{\im,4} \omega^2\bigr)(-q)^{p}&\text{if
$\fing^{(3)}=D^{(3)}_4$.}
\end{cases}
\end{align*}
Here $\omega$ is the third root of unity.  Recall that we follow the
enumeration of vertices of the Dynkin diagram of $\g$ as \cite{Kac}
except the $A^{(2)}_{2n}$ case given in \eqref{Eq: DD}.

\medskip
Then we have \eq
&&\chi_t\bl[V_{\fing^{(1)}}(\im,p)]\br=[V_{\fing^{(t)}}(\im,p)]
\qt{for any $(\im,p)\in\hI_{\fing^{(1)}}$.}\eneq

For each twisted quantum affine algebra $U_q'(\g)$, we associate a
finite simple Lie algebra $\gf$ and $\Psi_\g$ as follows:
\renewcommand{\arraystretch}{1.5}
\begin{align} \label{Table: root system for twisted}
\small
\begin{array}{|c||c|c|c|c|c|c|c|}
\hline
 \g    & A_{n}^{(2)} \ (n\ge 3)  & D_{n+1}^{(2)} \ (n\ge 3) & D_4^{(3)}   & E_6^{(2)}   \\ \hline
 (\gf,\Psi_\g) & \left(A_{n},\Psi_{A_n^{(1)}} \right) & \left(D_{n+1},\Psi_{D_{n+1}^{(1)}}\right) &\left(D_4,\Psi_{D_{4}^{(1)}} \right)          & \left(E_{6},\Psi_{E_{6}^{(1)}}\right)       \\  [1ex]  \hline
\end{array}
\end{align}
\noindent Also, the set of $\rmQ$-data associated to
$U_q'(\fing^{(t)})$ $(t \ge 2)$ are the same as the one for
$U_q'(\fing^{(1)})$. In particular, $\sigma=\id$ in the twisted
case.  Note also that the dual Coxeter number of $\fing^{(1)}$ and
the one of $\fing^{(t)}$ coincide.

In the sequel, we write simply $V(\im,p)$ for $V_\g(\im,p)$ if there
is no afraid of confusion.

\subsection{\pbw associated with a $\rmQ$-datum}\label{subsec:subcategory}

For each $\rmQ$-datum $\Qd$, the subcategory $\Ca_\Qd $ of $\Ca_\g$,
introduced in \cite{HL11}, is defined as the smallest subcategory of
$\Ca_\g$ containing $V(\im,p)$ for all $(\im,p) \in I_\Qd $ and the
trivial module $\mathbf{1}$ and is stable under taking subquotients,
extensions and tensor products.

Theorem~\ref{thm:gQASW duality 2} and Proposition~\ref{prop: Det
fun} below are proved for untwisted (resp.\ twisted) affine $A$ and
$D$ types in \cite[Theorem 4.3.1, Theorem 4.3.4]{KKK15B} (resp.\
\cite[Theorem 5.1]{KKKO16D}), for untwisted affine $B$ and $C$ types
in \cite[Theorem 6.3, Theorem 6.5]{KO18}, and for the remaining
exceptional affine types \cite[Theorem 6.3, Theorem 6.7, Theorem
6.13, Theorem 6.15]{OT19} (see also \cite[Proposition
6.5]{KKOP20C}):

\begin{theorem}[{\cite{KKK15B,KKKO16D,KO18,OT19}}] \label{thm:gQASW duality 2}
Let $U_q'(\g)$ be a quantum affine algebra  and let $\Qd$ be a
$\rmQ$-datum  associated to $\g$. Set
\begin{align}\label{eq: Vj}
V_\Qd(\al_\jm) \seteq V(\im,p) \qt{for $\jm\in\If$, where
$\phi_\Qd^{-1}(\al_\jm,0)=(\im,p)$.}
\end{align}
Then, we have \bna
\item The family $\Dd_\Qd\seteq\{ V_\Qd(\al_\jm) \}_{\jm \in \If}$ is a {\em complete}
 duality datum associated with the Cartan matrix $\cmf$ of type $\gf$. Hence the functor
$$\F_\Qd\seteq\F_{\Dd_\Qd}\col R_{\cmf}\gmod \rightarrow \Ca_\Qd   \text{ in~\eqref{eq: F fd} is exact.}$$
\item The functor $\F_\Qd$ sends simple modules to simple modules.
\ee
\end{theorem}

For an admissible sequence $\hw= \seq{(\im_k, p_k)}_{k\in\Z}$ in
$\If\times\Z$, we say that $(\Dd_\Qd,\hhw)$ is the associated \pbw,
where $\Qd$ is the corresponding $\rmQ$-datum and
$\hhw=(\im_k)_{k\in\Z}$.

\begin{proposition}[{\cite{KKK15B,KKKO16D,KO18,OT19}}]\label{prop: Det fun}
Let $\hw= \seq{(\im_k, p_k)}_{k\in\Z}$ be an admissible sequence in
$\If\times\Z$, and let $(\Dd_\Qd,\hhw)$ be the associated
\pbw.
Then we have
$$\cusp_k^{\hw}\seteq\cuspS{k}^{\Dd_\Qd,\hhw}\simeq V(\im_k,p_k).$$
\end{proposition}

We can say more.
By~\eqref{eq: KR1},
we have the following theorem.
\begin{theorem}  \label{Thm: KR as AD}
For an admissible sequence $\hw$ in $\hI_\g$ and any  $i$-box
$[a,b]$,
$$\text{$M^\hw[a,b]\seteq M^{\Dd_\Qd,\hhw}[a,b]$ is a KR-module over $U_q'(\g)$.}$$
\end{theorem}

Hence we can understand KR-modules as a special case of \Ad modules.

\begin{remark} \label{rem: KR}
Let $\hw$ be an admissible sequence. Then, for any $i$-box $[a,b]$
with $\im_a=\im_b=\im$ with $|[a,b]|_\phi = k$, we have
\begin{align*}
(\im_b,p_b) = (\im,p_a+2\fdi \times(k-1)) \quad \text{ and } \quad
M^{\hw}[a,b] = W^{(\pi(\im))}_{k,\epsilon\,(\qm)^{p_a}},
\end{align*}
where $\qm=q^{1/\ord(\sigma)}$ and
$\epsilon \in \C^\times$ is determined by~\eqref{eq: V(i;p)}
and~\eqref{eq: chi}.
\end{remark}

\begin{definition} \label{def: Category [a,b]}
For each interval $[a,b]$ and an admissible sequence $\hw$, we set
$$\Ca^{[a,b],\hw}_{\g}\seteq\Ca^{[a,b],\Dd_\Qd,\hhw}_{\g}.$$
Namely, $\Ca^{[a,b],\hw}_{\g}$ is the smallest full subcategory of
$\Ca_\g$  satisfying the following conditions$\colon$ \ben
\item it is stable under taking subquotients, extensions, tensor products and
\item  it contains $\cuspS{s}^{\hw}\simeq V_\Qd(\im_s,p_s)$ for all $a \le s \le b$ and the trivial module $\mathbf{1}$.
\end{enumerate}
\end{definition}

\smallskip
Now let us compare $\Ca_\g^{[a,b],\hw}$ with the subcategories
introduced in \cite{HL10,HL11,HL16}, which are stable under taking
subquotients, extensions, tensor products. Those categories depends
on the choice of $\hI_\g$, and we choose one. \bna
\item The subcategory $\Ca_\g^0$ of $\Ca_\g$, introduced in \cite{HL10,KKO18}, is defined as the smallest subcategory of $\Ca_\g$
containing $V(\im,p)$ for all $(\im,p) \in \hI_\g$. Thus $\Ca_\g^0$
can be identified with  $\Ca^{[-\infty,\infty],\hw}_{\g}$ for any
admissible sequence $\hw$.

\item Let $\hI_\g^- \seteq \hI_\g \cap (\If \times \Z_{\le 0})$. The subcategory $\Ca_\g^-$ of $\Ca_\g$, introduced in \cite{HL16}, is defined as the smallest subcategory of $\Ca_\g$
containing $V(\im,p)$ for all $(\im,p) \in \hI_\g^-$ and stable
under taking subquotients, extensions, tensor products. Take
a unique $\rmQ$-datum $\Qd$ with  $\xi_\im\in[1,2\fdi] $
 and $\hI_\Qd=\hI_\g$,
and let $\tw_0$ be a $\Qd$-adapted reduced expression of $w_0$. Now
let $\hw$ be the admissible sequence corresponding to $(\Qd,\tw_0)$
(see Proposition~\ref{prop:admseq}). Then, $\Ca_\g^-$ coincides with
the category
 $\Ca^{[-\infty,0],\hw}_{\g}$.
In a similar way, we can define $\Ca_\g^{< t}$ for $t \in \Z$
generated by $V(\im,p)$ with $(\im,p)\in\hI_\g^{< t} \seteq \hI_\g
\cap (\If \times \Z_{< t})$. Then it is equal to
$\Ca^{[-\infty,0],\hw}_{\g}$ taking $\hw$ such that $\hI_\g^{< t}
=\st{(\im_k,p_k)}_{k\le 0}$.

\item For each untwisted affine $\g$ of simply-laced type and $N \in \Z_{\ge 1}$, the subcategory $\Ca^N_\g$ of $\Ca_\g$, introduced in \cite{HL10}, is defined as the smallest subcategory of $\Ca_\g$
containing $V(\im,p)$ for all $(\im,p) \in \hI_\g \cap (\If \times
{[-2N+1,0]})$. By taking a $\rmQ$-datum  $\Qd$ with $[1,2]  \ni
\xi_\im $ and $(\im,\xi_\im) \in \hI_\g$, and a $\Qd$-adapted
reduced expression $\tw_0$, $\Ca_\g^N$  coincides with
$\Ca^{[a,0],\hw}_{\g}$, where $\hw$ is the admissible sequence
corresponding to $(\Qd,\tw_0)$ and $a = 1 - (N \times |\If|)$.

\item For each $\rmQ$-datum $\Qd$, the subcategory $\Ca_\Qd$ of $\Ca_\g$
coincides with  $\Ca^{[1,\ell],\hw}_{\g}$ for any corresponding
admissible sequence $\hw$.

\ee

Sometimes, we write $\Ca^{< \xi}_{\g}$ for the monoidal category
$\Ca_\g^{[-\infty,0],\hw}$, since it depends only on the choice of a
$\rmQ$-datum $\Qd= (\Dynkin, \sigma, \xi)$. Indeed, by setting
\begin{align} \label{eq: hI<xi}
\hI^{< \xi}_\g \seteq \{ (\im,p) \ | \  (\im,p) \in \hI_\Qd, p <
\xi_\im \},
\end{align}
we have $$\hI^{< \xi}_\g=\st{(\im_k,p_k)\mid k\le0}.$$

\vskip 2em
\section{Cluster algebra structure and monoidal categorification} \label{sec: monoidal categorification}
In this section, we briefly recall the definition of a cluster
algebra with  small  modifications as in \cite{KKOP19C}.
We also briefly review the main result of \cite{KKOP19C} on
$\Uplambda$-monoidal categorification of cluster algebras, which is
an application of the invariants $\La,\Li$ and $\de$ on $\Ca_\g$,
and can be understood as a quantum affine analogue of  the result
for quiver Hecke algebras in \cite{KKKO18}. After reviewing
several properties of monoidal seeds of  various kinds developed in
\cite{KK18} and \cite{KKOP20A}, we will construct \Lad monoidal seeds
associated to \pbws $(\Dd,\hhw)$ and the quivers $\GLS(\hhw)$
introduced in~\cite{GLS11}.
In the last part, we will review the result
in~\cite{HL16} which gives a cluster algebra structure on
$K(\Ca_{\g}^-)$ associated to the initial quiver $\HL$, and prove
that the initial seed of $K(\Ca_{\g}^-)$ in~\cite{HL16}  lifts to a
$\Uplambda$-admissible monoidal seed.
For more details on cluster algebras and monoidal categorification,
we refer the reader to \cite{BZ05,FZ02} and \cite{KKOP19C}.

\smallskip

{\em From now on, $\shc$ is  a full subcategory of $\uqm$
containing the trivial module $\one$ and stable
under taking tensor products, subquotients and extensions.} Note
that its Grothendieck group $K(\shc)$ has a ring structure with the
$\Z$-basis consisting of the isomorphism classes of simple modules.

\subsection{Cluster algebras} Fix a countable index set $\K=\K^\ex \sqcup \K^\fr$ which decomposes into a subset $K^\ex$
of {\em exchangeable} indices  and  a  subset $K^\fr$ of {\em
frozen} indices.

Let $\widetilde{B}=(b_{ij})_{(i,j)\in \K \times \Kex}$ be an
integer-valued matrix such that \eq &&
\parbox{77ex}{
\ben
\item  for each $j \in \Kex$, there exist finitely many $i \in \K$ such that $b_{ij} \ne 0$,
\item  the {\it principal part} $B \seteq (b_{ij})_{i,j \in \Kex}$ is skew-symmetric.
\end{enumerate}
}\label{eq: condition B} \eneq We call $\tB$ an {\em exchange
matrix.} We extend the definition of $b_{ij}$ for $(i,j)\in \K\times
\K$ by:
$$\text{
$b_{ij}=-b_{ji}$ if $i\in \Kex$ and $j\in \Kfr$,  and $b_{ij}=0$ for
$i,j\in \Kfr$, }$$ so that $(b_{ij})_{i,j\in \K}$ is skew-symmetric.

To the matrix $\tB$, we associate the quiver $\mathcal{Q}_{\tB}$
such that the set of vertices is $\K$ and the number of arrows from
$i\in\K$ to $ j\in\K$ is $\max(0,b_{ij})$. Then, $\mathcal{Q}_{\tB}$
satisfies the following conditions:
\begin{eqnarray} &&
\hs{2ex}\left\{\parbox{70ex}{
\begin{enumerate}[{\rm (a)}]
\item the set of vertices of $\mathcal{Q}_{\tB}$ is labeled by $\K$,
\item $\mathcal{Q}_{\tB}$ does not have any loop, any $2$-cycle, nor arrow between frozen vertices,
\item each  exchangeable vertex $v$ of $\mathcal{Q}_{\tB}$ has {\it finite degree}; that is, the number of arrows incident with $v$ is finite.
\end{enumerate}
}\right.\label{eq: Quiver condition}
\end{eqnarray}

Conversely, for a quiver satisfying~\eqref{eq: Quiver condition}, we
can associate a matrix $\tB$ satisfying~\eqref{eq: condition B} by
taking
\begin{align*}
b_{ij} \seteq \text{(the number of arrows from $i$  to $j$)}
\hspace{-.2ex}   -  \hspace{-.2ex}  \text{(the number of arrows from
$j$  to $i$)}.
\end{align*}

\smallskip
We say that a $\Z$-valued skew-symmetric $\K \times \K$-matrix
$L=(\la_{ij})_{i,j\in \K}$ is \emph{compatible} with
$\widetilde{B}$ (or $(L,\tB)$ is a compatible pair), if
$$   \sum_{k \in \K} \lambda_{ik}b_{kj} =2\delta_{i,j} \qquad \text{ for each $i \in \K$ and $j \in \Kex$.}$$

Let $\st{ X_i}_{i\in\K}$ be the set of mutually commuting
indeterminates.

\begin{definition}
For a commutative ring $\scrA$, we say that a triple $\Seed=(\st{ x_i
}_{i \in \K}, L,\tB) $ is a \emph{$\Uplambda$-seed} in $\scrA$ if
\ben
\item $\st{ x_i}_{i \in \K}$ is a family of elements of $\scrA$ and
there exists an injective algebra homomorphism $\Z[X_i;i \in \K]$ to $\scrA$ such that $X_i \mapsto x_i$,
\item $(L,\tB)$ is a compatible pair.
\end{enumerate}
\end{definition}

For a $\Uplambda$-seed $\Seed=(\st{ x_i }_{i \in \K}, L,\tB)$, we
call the set $\st{ x_i }_{i \in \K}$  the \emph{cluster} of $\Seed$
and its elements the \emph{cluster variables}. An element of the
form   $x^{{\bf a}}$ $\bigl({\bf a} \in \Z_{\ge 0}^{\oplus
\K}\bigr)$ is called a \emph{cluster monomial}, where
$$ x^{{\bf c}}  \seteq  \prod_{i \in \K} x_{i}^{c_i} \quad \text{ for }  \ {\bf c}=(c_{i})_{i\in\K} \in \Z^{\oplus \K}.$$

Let $\Seed=(\{ x_i \}_{i \in \K}, L,\tB)$ be a $\Uplambda$-seed in a
field $\mathfrak{K}$ of characteristic $0$. For $k \in \Kex$, we
define
\begin{enumerate}
\item[{\rm (a)}]\hs{2ex}$\mu_k(L)_{ij} =
\begin{cases}
  -\la_{kj}+\displaystyle\sum _{t\in\K} \max(0, -b_{tk}) \la_{tj} \quad \  & \text{if} \ i=k, \ j\neq k, \\
  -\la_{ik}+\displaystyle\sum _{t\in\K} \max(0, -b_{tk}) \la_{it} & \text{if} \ i \neq k, \ j= k, \\
   \la_{ij} & \text{otherwise,}
\end{cases}$

\vs{2ex}
\item[{\rm (b)}]\hs{2ex}$\mu_k(\tB)_{ij} =
\begin{cases}
  -b_{ij} & \text{if}  \ i=k \ \text{or} \ j=k, \\
  b_{ij} + (-1)^{\true(b_{ik} < 0)} \max(b_{ik} b_{kj}, 0) & \text{otherwise,}
\end{cases}
$

\vs{1.5ex}
\item[{\rm (c)}]\hs{2ex}$   \mu_k(x)_i  =\begin{cases}
x^{{\bf a}'}  +   x^{{\bf a}''}, & \text{if} \ i=k, \\
x_i & \text{if} \ i\neq k,
\end{cases}
$\\
\hs{1ex}where ${\bf a}'\seteq(a_i')_{i\in\K} \in \Z^{\oplus \K}$ and ${\bf
a}''\seteq(a_i'')_{i\in\K} \in \Z^{\oplus \K}$ are defined as
follows:
\begin{align*}
a_i'= \begin{cases}
  -1 & \text{if} \ i=k, \\
 \max(0,b_{ik}) & \text{if} \ i\neq k,
\end{cases} \qquad
a_i''= \begin{cases}
  -1 & \text{if} \ i=k, \\
 \max(0,-b_{ik}) & \text{if} \ i\neq k.
\end{cases}
\end{align*}
\end{enumerate}

Then the triple
$$\mu_k(\Seed) \seteq (   \{ \mu_k(x)_i\}_{i \in \K}, \mu_k(L),\mu_k(\tB) )$$
becomes a new $\Uplambda$-seed in $\mathfrak{K}$ and we call it the
\emph{mutation} of $\Seed$ at $k$.

The \emph{cluster algebra $\A(\Seed)$ associated to the
$\Uplambda$-seed} $\Seed$
 is the $\Z$-subalgebra of the field $\mathfrak{K}$ generated by all the cluster variables in the $\Uplambda$-seeds obtained from
$\Seed$ by all possible successive mutations.

A \emph{cluster algebra structure associated to a $\Uplambda$-seed
$\Seed$} on a $\Z$-algebra $\A$ is a family $\mathscr{F}$ of
$\Uplambda$-seeds in $\A$ such that
\ben
\item for any $\Uplambda$-seed $\Seed$ in $\mathscr{F}$, the cluster algebra $\A(\Seed)$ is isomorphic to $\A$,
\item any mutation of a $\Uplambda$-seed in $\mathscr{F}$ is in $\mathscr{F}$, \item for any pair $\Seed$, $\Seed'$ of $\Uplambda$-seeds in $\mathscr{F}$,
$\Seed'$ can be obtained from $\Seed$ by a finite sequence of
mutations.
\end{enumerate}

\subsection{Monoidal seeds and their mutations}

Let $\shc$ be a full subcategory of $\Ca_\g$
containing the trivial module $\one$ and stable under taking
tensor products, subquotients and extensions.

\begin{definition} \hfill
\ben
\item A \emph{monoidal seed in $\shc$} is
a  quadruple $\seed = (\{ M_i\}_{i\in \K },\widetilde B\KK)$
consisting of a  commuting family $\{ M_i\}_{i\in\K}$ of
real simple modules in $\shc$,  an  integer-valued
$\K\times\Kex$-matrix $\widetilde B =
(b_{ij})_{(i,j)\in\K\times\Kex}$ satisfying the conditions
in~\eqref{eq: condition B}, an index set $\K$ and an index set
$\Kex\subset \K$ of exchangeable vertices.
\item For $i\in\K$, we call $M_i$ the $i$-th {\em cluster variable module} of $\seed$.
\end{enumerate}
\end{definition}

For a monoidal seed  $\seed=(\{M_i\}_{i\in\K}, \widetilde B\KK)$,
let $\La^\seed=(\La^\seed_{ij})_{i,j\in\K}$ be the skew-symmetric
matrix given by $\La^\seed_{ij}=\Lambda(M_i,M_j)$.

\begin{definition} \label{def:admissible}
We say that a monoidal seed $\seed=(\{M_i\}_{i\in\K}, \widetilde
B\KK)$  in $\shc$ is  \emph{admissible} if it satisfies the following conditions:
\ben
\item for each
$k\in\Kex$, there exists a simple object $M'_k$ of $\shc$  such that
 there is an exact sequence in $\shc$
\begin{align*}
0 \to  \dtens_{b_{ik} >0} M_i^{\tens  b_{ik}} \to M_k \otimes M_k'
\to \dtens_{b_{ik} <0} M_i^{\tens  (-b_{ik})} \to 0,
\end{align*}

\item
$M_k'$
commutes with $M_i$  for each $k\in \Kex$ and any $i\in
\K\setminus\st{k}$.
\ee
\end{definition}

Note that we have also an
exact sequence \eqn
 &&0 \to  \dtens_{b_{ik} <0} M_i^{\tens(-b_{ik})} \to M_k' \tens M_k \to
  \dtens_{b_{ik} >0} M_i^{\tens b_{ik}} \to 0.\label{eq:ses_mutation2}
\eneqn
Note also that $M'_k$ is unique up to an isomorphism if it exists, since
$M_k\hconv M'_k\simeq \dtens_{b_{ik} <0} M_i^{\tens  (-b_{ik})}$.

\Lemma If  a monoidal seed $\seed=(\{M_i\}_{i\in\K}, \widetilde B\KK)$ is {\em admissible}, $M'_k$
in Definition~\ref{def:admissible} is real for any $k\in\Kex$.
\enlemma
\Proof
It follows from Proposition~\ref{prop:real}~\eqref{it:real2}.
\QED

By this lemma, if $\seed$ is admissible,
then the quadruple $$\mu_k(\seed)
\seteq  \bl\{M_i\}_{i\neq k}\cup\{M_k'\},\mu_k(\widetilde B)\KK\br$$
is a monoidal seed in $\shc$. We call $ \mu_k(\seed)$ the {\em
mutation} of $\seed$ in direction $k$.

\begin{proposition}[{\cite[Proposition 6.4]{KKOP19C}}] \label{prop:condition simplified}
Let $\seed=(\{M_i\}_{i\in\K},\widetilde B\KK)$ be an admissible
monoidal seed in $\shc$. Let  $k\in\Kex$ and let $M'_k$ be as in
{\rm Definition~\ref{def:admissible}}. Then we have the following
properties. \bnum
\item \label{it: d2}  For any $j \in \K$, we have $(\Lambda^\seed\;\widetilde{B})_{jk}=-2 \de(M_j,M_k')$.
\item For any $j\in\K$, we have
\begin{align*}
&\La(M_j,M'_k)=-\La(M_j,M_k)-\sum_{b_{ik}<0}\La(M_j,M_i)b_{ik},\\
&\La(M'_k,M_j)=-\La(M_k,M_j)+\sum_{b_{ik}>0}\La(M_i,M_j)b_{ik}.
\end{align*}
\end{enumerate}
\end{proposition}

\Def[{\cite[Definition 6.5]{KKOP19C}}] Let  $\seed=(\{M_i\}_{i\in
\K},\widetilde B; \K,\Kex)$ be an admissible monoidal seed.
\ben
\item
We say that a monoidal seed $\seed$ is \emph{$\Uplambda$-admissible}
if  $M'_k$ in Definition~\ref{def:admissible} satisfies
$\de(M_k,M'_k)=1$.
\item
If $\seed$ is $\Uplambda$-admissible, we say that the mutation
$\mu_k(\seed)$ of $\seed$ at $k\in \Kex$ is a {\em
$\Uplambda$-mutation},
\item
We say that a monoidal seed $\seed$ is \emph{\ca} if $\seed$ admits
successive $\Uplambda$-mutations in all possible directions. \ee
\edf

For a monoidal seed $\seed=(\{M_i\}_{i\in \K},\widetilde B)$ in
$\shc$, we define the triple $[\seed]$ in $K(\shc)$ by
$$[\seed]\seteq\bl\{ [M_i] \}_{i \in \K},-\La^\seed,\tB\br.$$

If $\seed$ is a $\Uplambda$-admissible monoidal seed, then $[\seed]$
is a $\Uplambda$-seed.

\begin{definition}[{\cite[Definition 6.7]{KKOP19C}}]
A category $\shc$ is called a \emph{$\Uplambda$-monoidal
categorification} of a cluster algebra $\A$ if \ben
\item the Grothendieck ring $K(\shc)$ is isomorphic to $\A$,
\item there exists a \ca monoidal seed $\seed=(\{ M_i \}_{i \in \K},\tB\KK)$ in $\shc$ such that
$$[\seed]\seteq(\{ [M_i] \}_{i \in \K},-\La^\seed,\tB)$$
is  an initial $\Uplambda$-seed of $\A$.
\end{enumerate}
\end{definition}

Now we present the main result of \cite{KKOP19C}:

\Th[{\cite[Theorem 6.10]{KKOP19C}}] \label{th:main KKOP19C} Let
$\seed=(\{M_i\}_{i\in \K},\widetilde B\KK)$ be a
$\Uplambda$-admissible monoidal seed  in $\shc$, and set
$$[\seed]\seteq(\{ [M_i]\}_{i\in\K}, -\Lambda^\seed,\widetilde B).$$
We assume that the algebra $K(\shc)$ is isomorphic to the cluster
algebra $\A([\seed])$. Then, we have

\smallskip
\hs{5ex}\parbox{55ex}{ $\bullet$ $\seed$ is \ca, and
\\
$\bullet$ $\shc$ gives a $\Uplambda$-monoidal categorification of
$\A([\seed])$.} \enth

\begin{definition} \label{def: max com}
A family of real simple modules $\{M_i\}_{i\in \K}$ in $\shc$ is
called a {\em real commuting family in $\shc$} if it
satisfies:
\ben
\item $\{M_i\}_{i\in \K}$ is mutually  commuting.
\ee
It is called  a {\em maximal real commuting family in $\shc$} if it
satisfies further : \ben
\item[{\rm(2)}]
if a simple  module $X$ commutes with all the $M_i$'s, then
$X$ is isomorphic to $\bigotimes_{i \in \K}M_i^{\tens a_i}$ for some
$\mathbf{a}=\st{a_i}_{i\in \K} \in \Z_{\ge 0}^{\oplus \K}$. \ee
\end{definition}

\begin{corollary}[{\cite[Corollary 6.11]{KKOP19C}}]\label{cor:main KKOP19C}
Let  $\seed=(\{M_i\}_{i\in\K},\widetilde B\KK)$ be a
$\Uplambda$-admissible  monoidal
seed in $\shc$ and assume  that the algebra $K(\shc)$ is isomorphic
to $\A([\seed])$. Then  the following statements hold{\rm:} \bnum
  \item Any cluster monomial in $K(\shc)$  is the isomorphism class of a  real simple object in $\shc$.
  \item Any cluster monomial in $K(\shc)$ is a Laurent polynomial of the initial cluster variables with coefficient in $\Z_{\ge0}$.
\item For $k\in\Kex$ and the $k$-th cluster variable module $\widetilde{M}_k$ of a monoidal seed $\widetilde{\seed}$ obtained by successive $\Uplambda$-mutations from the initial monoidal seed
$\seed$,  we have
$$\de(\widetilde{M}_k,\widetilde{M}_k')=1.$$
Here $\widetilde{M}_k'$ is the $k$-th cluster variable module of
$\mu_k(\widetilde{\seed})$.
\item  Any monoidal cluster $\{\widetilde{M}_i\}_{i\in\K}$ is a maximal real commuting family.
\end{enumerate}
\end{corollary}

\subsection{Properties of $\Uplambda$-admissible monoidal seeds} \label{subsec: Property of ad mseed}

Recall the definitions of $\rW_0$, $\Delta_0$, $\rE(M)$ ($M \in \Ca_\g$) given in \S~\ref{subsec: E(M)}.

\Lemma\label{lem:EB} Let  $\seed=(\{M_i\}_{i\in\K},\widetilde B\KK)$
be  an admissible monoidal seed in
$\Ca_\g$. Then we have
$$\sum_{i\in \K}\rE(M_i) b_{ik} =0\qt{for any $k\in\Kex$,}$$
where $\tB=( b_{ij})_{(i,j)\in \K\times \Kex}$.
\enlemma
\Proof By
the definition, there is a short exact sequence with a simple $M'_k$:
\begin{align*}
0 \to  \dtens_{b_{ik} >0} M_i^{\tens  b_{ik}} \to M_k \otimes M_k'
\to \dtens_{b_{ik} <0} M_i^{\tens  (-b_{ik})} \to 0,
\end{align*}
Then we have \eqn \sum_{b_{ik}
>0}\rE(M_i)b_{ik}&&=\rE\bl\dtens_{b_{ik} >0} M_i^{\tens
b_{ik}}\br
\\&&=\rE(M_k)+\rE(M_k')=
\rE\bl\dtens_{b_{ik} <0} M_i^{\tens  (-b_{ik})}\br =\sum_{b_{ik}
<0} \rE(M_i)(-b_{ik}) \eneqn by \cite[Lemma 3.1]{KKOP20A}. Hence we
have the desired result. \QED

\Lemma[{cf.\ \cite[Lemma 3.2]{KK18}}]\label{lem:uniqueB} Let
$\seed=(\{M_i\}_{i\in\K},\widetilde B \KK)$ be a
 $\Uplambda$-admissible in $\catCO$,  and
$\tB=(b_{ij})_{(i,j)\in \K\times \Kex}$. Assume that $\K$ is a
finite set. \bnum
\item Then we have
$\dim\bl\sum_{i\in\K}\Q\rE(M_i)\br\le |\Kfr|$.
\item Assume further that $\dim\bl\sum_{i\in\K}\Q\rE(M_i)\br=|\Kfr|$.
Then, for any $k\in\Kex$, $(b_{ik})_{i\in \K}$ is a unique element
$(v_i)_{i\in \K}$ of $\Q^\K$ such that \eqn \sum_{i\in
\K}\rE(M_i)v_i=0\qtq \sum_{i\in
\K}(\La^\seed)_{ji}v_i=-2\delta_{j,k} \qt{for any $j\in\Kex$.}
\eneqn \ee \enlemma \Proof Let $f\col\Q^{\oplus \K}\To \Q^{\oplus
\K\ex}\soplus \rW_0$ be the linear map given by $(\La^\seed)_{ki }$
and $\rE(M_i)$. Then, $\Im(f)$ contains $\Q^{\oplus \K\ex}\soplus 0$
by Proposition~\ref{prop:condition simplified} and
Lemma~\ref{lem:EB}. Moreover, the image of the composition
$f\col\Q^{\oplus \K}\To \Q^{\oplus \K\ex}\soplus \rW_0\to\rW_0$  is
$\sum_{i\in \K}\Q\rE(M_i)$, which implies that $\Im(f)=\Q^{\oplus
\K\ex}\soplus\bl\sum_{i\in \K}\Q\rE(M_i)\br$. Hence the dimension of
$\Im(f)$ is equal to $|\Kex|+\dim\bl\sum_{i\in \K}\Q\rE(M_i)\br$.
Then, $|\K|\ge\dim\bl\Im(f)\br$ implies (i).

\snoi (ii) If $\dim\bl\sum_{i\in\K}\Q\rE(M_i)\br=|\Kfr|$, then $f$
is injective. \QED

\Prop\label{prop:mutation must be}
Let  $\seed=(\{M_i\}_{i\in\K},\widetilde B\KK)$ be
a  $\Uplambda$-admissible monoidal seed  in
$\catCO$ with $\tB=(b_{ij})_{(i,j)\in \K\times \Kex}$, and let $k
\in \Kex$. Assume that \bnum
\item
$\K$ is a finite set and $\dim\bl\sum_{i\in\K}\Q\rE(M_i)\br\ge
|\Kfr|$,
\item there exist a real simple module $X \in \shc$
and an exact sequence
\begin{align*}
0 \to A \to M_k \tens X \to B \to 0,
\end{align*}
such that \bna
\item $X$ commutes with $M_j$ for all $j \in \K \setminus \{k\}$,
\item $\de(M_k,X)=1$,
\item $A=\dtens_{i \in \K} M_i^{\tens m_i}$,
$B=\dtens_{i \in \K} M_i^{\tens n_i}$ for some
$m_i,n_i\in\Z_{\ge0}$. \ee \ee Then we have $b_{ik}=m_i-n_i$.

If we have furthermore $m_in_i=0$ for all $i \in \K$, then  we have
$$X \simeq M'_k,$$
where $M'_k$ is given in {\rm Definition~\ref{def:admissible}}. \enprop
\Proof We shall apply Lemma~\ref{lem:uniqueB}. Set
$\la_{ij}=\La(M_i,M_j)$.
 Then we have
\eqn \sum_{i\in \K}m_i \rE(M_i)=
\rE(A)=\rE(M_k)+\rE(X)=\rE(B)=\sum_{i\in \K}n_i\rE(M_i).
\eneqn For any $j\in\K$, we have \eqn\sum_{i\in
\K}  \la_{ij} m_i&&=\La(A,M_j)=\La(X\hconv M_k, M_j)
=\La(X, M_j)+\La(M_k, M_j)\\
\sum_{i\in \K} \la_{ji}  n_i&&=\La(M_j,B)=\La(M_j,M_k\hconv
X)=\La(M_j, M_k)+\La(M_j,X) \eneqn Hence we have \eqn \sum_{i\in
\K}\la_{ji} (n_i-m_i)&&=
\La(X, M_j)+\La(M_k, M_j)+\La(M_j, M_k)+\La(M_j,X)\\*
&&=2\bl\de(M_j,X)+\de(M_j, M_k)\br=2\delta_{j,k}. \eneqn Thus,
Lemma~\ref{lem:uniqueB} implies that $b_{ik} =m_i-n_i$.

\mnoi Now assume that $m_in_i=0$ for all $i \in \K $. Then we have
$n_i=\max\st{0,- b_{ik} }$ and
$$B\simeq\dtens_{ b_{ik} <0}M_i^{\otimes (-b_{ik} )}.$$
Hence, we obtain
$$M_k\hconv X\simeq B\simeq M_k\hconv M'_k,$$
which implies that $X\simeq M'_k$. \QED

As an immediate application of
Proposition~\ref{prop:mutation must be},
we can show that the exchange matrix is uniquely determined for
a   $\Uplambda$-admissible  monoidal seed.

\Prop\label{prop:uniqueexchange} Let
$\seed_0=(\{M_i\}_{i\in\K_0},\tB_0; \K_0,\Kex_0)$ and
$\seed=(\{M_i\}_{i\in\K},\tB \KK)$ be two
$\Uplambda$-admissible  monoidal seeds in $\catCO$ such
that $\K_0\subset \K$ and $\Kex_0\subset\Kex$. Assume that $\K$ is a
finite set and $\dim\bl\sum_{i\in\K}\Q\rE(M_i)\br\ge |\Kfr|$. Then
$$\tB\vert_{\K_0\times\Kex_0}=\tB_0 \quad \text{ and }\quad
\tB\vert_{(\K\setminus\K_0)\times\Kex_0}=0.$$ \enprop

The following lemma is almost obvious by the definition.
\Lemma\label{lem:restrseed} Let $\seed=(\{M_i\}_{i\in \K},\tB\KK)$
be a monoidal seed in $\shc$. Let $\K_0$ be a subset of $\K$ with a
decomposition $\K_0=\Kex_0\sqcup\Kfr_0$ such that $\Kex_0\subset
\Kex$. Set
$$\seed\vert_{(\K_0,\,\Kex_0)}\seteq
\bl\st{M_i}_{i\in\K_0},\,\tB\vert_{(\K_0)\times\Kex_0};\K_0,\,\Kex_0\br.$$
Assume that \eq\text{$b_{ij} =0$ if $i\in\K\setminus\K_0$ and
$j\in\Kex_0$.} \label{eq:inout} \eneq Then, we have \bnum
\item
$\bl\mu_s(\tB)\br_{ij }=0$ if $s\in\Kex_0$, $i\in\K\setminus \K_0$
and $j\in\Kex_0$,
\item if $\seed$ is \Lad, then we have
\eqn \bl\mu_s\seed\br\vert_{(\K_0,\Kex_0)} =\bc
\mu_s(\seed\vert_{(\K_0,\,\Kex_0)})&\text{if $s\in\Kex_0$,}\\
\seed\vert_{(\K_0,\,\Kex_0)}&\text{if $s\in\K\setminus \K_0$.}
\ec\eneqn \ee In particular, if $\seed$ is a \ca monoidal seed in
$\shc$, then so is $\seed\vert_{(\K_0,\,\Kex_0)}$. \enlemma

\subsection{Monoidal seeds and admissible chains of $i$-boxes} \label{subsec: Prop of monoidal seed} Let
$(\Dd,\hhw)$ be a \pbw. Throughout this subsection we consider
admissible chains $\frakC$ of $i$-boxes associated with
$(\Dd,\hhw)$.

Let $\frakC=\seq{\ci_k}_{1\le k\le l}$ be an admissible chain of
$i$-boxes with a range $[a,b]$. We define \eq &&\ba{rl}
\K(\frakC)&\seteq[1,l],\\
\Kfr(\frakC)&\seteq\st{s\in\K(\frakC)\mid
\text{$\ci_s=[a(\im)^+,b(\im)^-]$ for some $\im\in\Iff$}},\\
\Kex(\frakC)&\seteq\K(\frakC)\setminus \Kfr(\frakC),\\
\cM[](\frakC)&\seteq\st{M(\ci_k)}_{k\in\K(\frakC)}. \ea
\label{def:chainseed} \eneq Here, if $l=\infty$, we understand that
$\Kfr(\frakC)$ is the empty set. Recall that $\cM[](\frakC)$ is a
commuting family of real simple modules (see Theorem~\ref{thm:
admissible chain commuting}).

\smallskip

We shall first prove the following lemma that assures that
Proposition~\ref{prop:mutation must be}  is applicable to
$\seed(\frakC)$ for any admissible chain $\frakC$ of $i$-boxes with
a finite range.

\Lemma\label{lem:smallfrozen} Let $\frakC=(\ci_k)_{1\le k\le l}$ be
an admissible chain of $i$-boxes associated with $\hhw$ and
assume that its range $[a,b]$  is finite.  Then we have
$$\dim\Bigl(\sum_{1\le k\le l}\Q\rE\bl M(\ci_k)\br\Bigr)=|\Kfr(\frakC)|=\big|\st{\im_s\mid
s\in[a,b]}\big|.$$ \enlemma \Proof By Proposition~\ref{prop:
changing order} and Proposition~\ref{prop: real T-move}, $\sum_{1\le
k\le l}\Q\rE\bl M(\ci_k)\br$ does not change by box moves. Hence, we
have $\rE(\bS{s})\in \sum_{1\le k\le l}\Q\rE\bl M(\ci_k)\br$ for any
$s\in[a,b]$. Therefore,  we have
$$\sum_{1\le k\le l}\Q\,\rE\bl M(\ci_k)\br=
\sum_{s\in [a,b]}\Q\,\rE(\bS{s}) \simeq\sum_{s\in
[a,b]}\Q\,\al_{\im_s},$$ which implies that $\dim\Bigl(\sum_{1\le
k\le l}\Q\rE\bl M(\ci_k)\br\Bigr)=\big|\st{\im_s\mid
s\in[a,b]}\big|$. \QED

The following lemma  says that box moves correspond to
mutations.
\begin{lemma} \label{lem: b is mu}
Let $\mathfrak{C}=(\ci_k)_{1\le k\le l}$ be an admissible chain of
$i$-boxes associated with $\hhw$ and a finite range such that
$\seed\seteq\bl\cM[](\frakC),\tB;\K(\frakC),\Kex(\frakC)\br$ is a
 \Lad  monoidal seed in
$\catCO$ for some exchange matrix $\tB$. If $k_0\in \Kex(\frakC)$
and $\ci_{k_0}$ is a movable $i$-box such that
$\tc_{k_0+1}=\ci_{k_0+1}=[a,b]$, then we have
\eqn
\mu_{k_0}(\seed)
&& =\bl \cM[](B_{k_0}(\frakC)),\mu_{k_0}(\tB);
\K(B_{k_0}(\frakC)),\Kex(B_{k_0}(\frakC))\br \\
&&=  \bl\{ M_i\}_{i\in \K \setminus   \{ k_0 \}   } \sqcup\{ M_{k_0}'  \},
\mu_{k_0}(\tB)\KKC\br,
\eneqn
where
$$ M'_{k_0} \seteq \begin{cases}
M[a,b^-]  & \text{if $\ci_k=[a^+,b]$,} \\
M[a^+,b]  & \text{if $\ci_k=[a,b^-]$.}
\end{cases}
$$
Thus $B_k$ in {\rm Remark~\ref{rmk: combinatorial T-system}}
corresponds to $\mu_k$ in this case and the mutation $\mu_k$
corresponds to T-system in~\eqref{eq: T-system in terms of M[a,b]
again}.
\end{lemma}

\begin{proof}
By Theorem~\ref{thm: commuting ab} and Remark~\ref{rmk:
combinatorial T-system}, the modules $  M[a(\jm)^+,b(\jm)^-]$
($d(\im_a,\jm)=1$) and $M[a^+,b^-]$ commute with $M[a,b]$, and they
are contained in $\cM[](\frakC)$. Thus our assertion follows from
Proposition~\ref{prop: real T-move},
Proposition~\ref{prop:mutation must be}
together with Lemma~\ref{lem:smallfrozen}.
\end{proof}

Together with Proposition~\ref{prop: changing order}, we obtain the
following corollary.
\Cor\label{cor:2chains} For  a finite
interval $[a,b]$, let $\frakC$ and $\frakC'$ be admissible chains of
$i$-boxes associated with the same $\hhw$ and the same range
$[a,b]$. Assume that $(\cM[](\frakC),\tB\KKC)$ is a \ca
monoidal seed in $\catCO$ for some exchange matrix $\tB$.
Then, $(\cM[](\frakC'),\tB'\KKCp)$ is also a \ca monoidal seed in $\catCO$ for
some exchange matrix $\tB'$. \encor

\Prop\label{prop:seedab}
 Let $(\Dd,\hhw)$ be a complete \pbw and
let $\seed=\bl\st{ M_i}_{i\in\K},\tB\KK)$ be a monoidal seed in
$\Ca^{[a,b],\Dd,\hhw}_\g$. If $\seed$ is \rmo completely\rmf
$\Uplambda$-admissible in $\catCO$, then it is \rmo completely\rmf
$\Uplambda$-admissible in $\Ca^{[a,b],\Dd,\hhw}_\g$.
\enprop
\Proof
It is enough to show that the mutation $M'_k$ of $\seed$ at a vertex
$k\in\Kex$ belongs to $\Ca^{[a,b],\Dd,\hhw}_\g$. Since $M_k$,
$M_k\hconv M'_k$ and $M'_k\hconv M_k$ belong to
$\Ca^{[a,b],\Dd,\hhw}_\g$, our assertion follows from
Corollary~\ref{cor:MX}. \QED

\subsection{An example of   \Lad monoidal seeds}

Let $(\Dd,\hhw)$ be a \pbw.
Let $\GLS(\hhw)$  be the quiver structure on $\Z$ with two types of
arrows (with the notations in  \eqref{eq: nota +,-}): \eq
&&\hs{2ex}\ba{llll}
\text{vertical arrow}&:\hs{2ex}& s\To[\akew]t\hs{5ex}&\text{if $s^-<t^-<s<t$ and $d(\im_s,\im_t)=1$,}\\[.5ex]
\text{horizontal arrow}\hs{2ex}&:& s\To[\akew]s^-&
\ea\label{quiver:Omegaminus} \eneq

If there is no afraid of confusion, we write  shortly $\GLS$ for
$\GLS(\hhw)$.

For an interval $[a,b]$, we denote by $\GLS^{[a,b]}(\hhw)$ its
induced quiver on $[a,b]$,
 and by $\BGLS^{[a,b],\hhw}$ the exchange matrix associated with
the quiver $\GLS^{[a,b]}(\hhw)$.

\begin{theorem}\label{thm:triangle}
Let $(\Dd,\hhw)$ be a \pbw. For $-\infty\le a\le b<+\infty$, the
monoidal seed in $\cat^{[a,b],\Dd,\hhw}_\g$ \eq
\seed^{[a,b],\Dd,\hhw}_-\seteq\bl\st{M^{\Dd,\hhw}[s,b\}}_{s\in[a,b]},\BGLS^{[a,b],\hhw}\KK\br
\label{def:seedab} \eneq is \Lad . Here $\K=[a,b]$,
$\Kex=\st{s\in[a,b]\mid a\le s^-}$.
\end{theorem}

Note that
$\seed^{[a,b],\Dd,\hhw}_-=\seed(\frakC^{[a,b],\Dd,\hhw}_-)$ with the
admissible chain of $i$-boxes
$$\frakC^{[a,b], \Dd,\hhw}_-\seteq  \seq{[b+1-k,b\}}_{1\le k\le b-a+1}. $$

We devote   this subsection for the proof of
Theorem~\ref{thm:triangle}. We employ the framework of the
proof of \cite[Theorem 11.2.2]{KKKO18}.

Set $\cM=M^{ \Dd,\hhw}[s,b\}$, $\K=[a,b]$ and \eqn
\Kfr&&=\st{s\in[a,b]\mid s^-<a\le s},\\
\Kex&&=\K\setminus\Kfr=\st{s\in\Z\mid a\le s^-<s\le b},\\
\BGLS^{[a,b],\hhw}&&=( b_{ij})_{(i,j)\in\K\times\Kex}. \eneqn

\medskip
In order to prove Theorem~\ref{thm:triangle}, we will show
\eq
\left\{\parbox{75ex}{\begin{enumerate}[{\rm (a)}]
\item \label{item:a}  $\{ \cM\mid s\in \K \}$ is a real
commuting family,
\item \label{item:c} for $s\in\Kex$, there exist a simple module
${\cM'}$ and an exact sequence
$$0 \To   \dtens_{ b_{ts}  >0}  {\cM[t]}^{\tens b_{ts} } \To
\cM \tens {\cM'} \To  \dtens_{ b_{ts}  <0} {\cM[t]}^{\tens
(-b_{st})} \To 0,
$$
\item \label{item:b}$\de(\cM,\cM')=1$,
\item \label{item: admf} $\cM'$ is real and $\cM'$ commutes with $\cM[t]$ for any $a \le t \ne s \le b$.
\end{enumerate}
}\right.\label{eq:adm}
\eneq

Note that, \eqref{item:a} is already proved in Theorem~\ref{thm:
admissible chain commuting}. Thus it suffices to prove
\eqref{item:c},  \eqref{item:b} and~\eqref{item: admf}.

\smallskip

For a vertex $ s\in\K$, the arrows incident to $s$ can be classified
into four types as follows (see \eqref{quiver:Omegaminus}): (i)
horizontal incoming arrows : $s^+\To s$ (ii) horizontal outgoing
arrows $s\To s^-$ (iii) vertically incoming arrow  : $t\To s$ with
$d(\im_s,\im_t) =1$ and $t^-<s^-<t<s$, and (iv) vertically outgoing
arrow : $s\To t$ with $d(\im_s,\im_t) =1$ and $s^-<t^-<s<t$.

\medskip
Let us fix a vertex $s\in \Kex$ and set $\im=\im_s$. We define the
subsets of $\K$ as follows:
\eq\label{eq: Hi Ho}
\ba{rrl}
\Vi(s) &\seteq& \{ t\in \Z \mid \text{$d(\im,\im_t)=1$ and $t^-<s^-<t<s$}\}\\
&=&\st{s^-(\jm)^+\mid d(\im,\jm)=1,\, s^-<s(\jm)^-},\\[1ex]
\Vo(s) &\seteq& \st{ t\in \Z \mid \text{$d(\im,\im_t)=1$ and $s^-<t^-<s<t \le b $}}\\
&=&\st{s(\jm)^+\mid d(\im,\jm)=1,\, s^-<s(\jm)^-,\,s(\jm)^+\le b}.
\ea
\eneq
We set \eqn \cM^\Vi&&\seteq\dtens_{t\in\Vi(s)} \cM[t]\qtq
\cM^\Vo\seteq\dtens_{t\in\Vo(s)} \cM[t]. \eneqn

Then we have
\begin{align*}
\dtens_{ b_{ts}  >0} \cM[t]^{\tens  b_{ts} }  \simeq \cM[s^+] \tens
\cM^{\Vi} \quad \text{ and } \quad \dtens_{b_{ts} <0}
\cM[t]^{\tens (- b_{ts} )}  \simeq  \cM[s^-] \tens \cM^{\Vo}.
\end{align*}

Now let us show the following lemma. \Lemma\label{thm:step1} The
module
\begin{align} \label{eq: M'}
 \cMp \seteq   \cM^\Vo\hconv \bS{s^-}
 \end{align}
satisfies the following properties:
\bnum
\item $\cMp$ is a simple module.
\item \label{-thm a} $\cM\hconv \cMp \simeq
\dtens_{  b_{ts} <0} \cM[t]^{\tens (-b_{ts}  )}\simeq\cM[s^-] \tens
\cM^{\Vo} $.
\item \label{-thm b} $\cMp \hconv \cM\simeq
\dtens_{b_{ts}>0}  \cM[t]^{\tens b_{ts} }\simeq \cM[s^+] \tens
\cM^{\Vi}$.
\end{enumerate}
\enlemma

Note that $ \cMp$ in~\eqref{eq: M'} is simple by Theorem~\ref{thm:
KKKo15 main} because $\cM^\Vo$ is simple.

\Sub \label{lem: - head} We have
$$\cM \hconv \cMp  \simeq\cM[s^-] \tens \cM^{\Vo} .$$
\ensub

\begin{proof}
We have
$$
\cM[s] \hconv \bS{s^-}\simeq \cM[s^-].$$

On the other hand, we have
\begin{align*}
 & \cM \tens \left( \cM^{\Vo} \hconv \bS{s^-} \right) \monoto
 \cM \tens \bS{s^-} \tens \cM^{\Vo}
\epito \cM[s^-]\tens \cM^{\Vo},
\end{align*}
whose composition does not vanish by
Proposition~\ref{prop:rcomp}~\eqref{item1}. Thus our assertion
follows from Theorem~\ref{thm: KKKo15 main}.
\end{proof}

 Note that we have
$$M[s,b\}\simeq M[s^+,b\}\hconv \bS{s}\monoto \bS{s}\tens M[s^+,b\}.
$$
Thus we have the following monomorphism:
\begin{align*}
\big( \cM^{\Vo} \hconv  \bS{s^-} \big) \tens \cM & \monoto
 \big( \cM^{\Vo} \hconv  \bS{s^-} \big) \tens \bS{s} \tens \cM[s^+].
\end{align*}

Hence, in order to prove~\eqref{-thm b} in Lemma~\ref{thm:step1}, it
is suffices to show that there exists an epimorphism \eq  \big(
\cM^{\Vo} \hconv  \bS{s^-} \big) \tens \bS{s}
  \epito  \cM^{\Vi},\label{eq: epi}
\eneq by Proposition~\ref{prop:rcomp}~\eqref{item1}.

\smallskip

Note that we have a surjective homomorphism
\begin{equation*}%
\cM^{\Vo} \tens  \bS{s^-} \tens \bS{s} \epito  \big( \cM^{\Vo}
\hconv \bS{s^-} \big) \tens \bS{s}.
\end{equation*}

\Sub \label{prop: normal for the seed} The sequence $    \left(
\cM^{\Vo} ,  \bS{s^-}, \bS{s} \right)$  is normal. In particular,
$\hd\bl  \cM^{\Vo}\tens \bS{s^-}\tens \bS{s}\br$ is simple. \ensub

\begin{proof}
Note that if $t\in\Vo(s)$, then $t>s$. Hence we have
$$ \de( \cM^{\Vo}, \D^{-1} \bS{s} )=0,$$
by ~\eqref{it: unmixed} in Proposition~\ref{prop: lemma of KKOP20C},
i.e. $( \cM^{\Vo},\bS{s})$ is unmixed. Hence,  our assertion follows
from Proposition~\ref{prop: L N* normal}.
\end{proof}
Hence we have
$$\hd\left( \cM^{\Vo} \tens \bS{s^-} \tens  \bS{s}\right)  \simeq \big( \cM^{\Vo} \hconv \bS{s^-} \big) \hconv \bS{s}.$$

\Sub\label{lem: - socle} For $s\in\Kex$, we have
$$\hd\left( \cM^{\Vo} \tens \bS{s^-} \tens  \bS{s}\right)  \simeq \cM^\Vi. $$
\ensub

\begin{proof}
By T-system described in~\eqref{eq: T-system in terms of M[a,b]
again}, we have
\begin{align*}
 \cM^{\Vo} \tens \bS{s^-} \tens \bS{s} & \twoheadrightarrow
\cM^\Vo\tens \left( \dtens_{d(\im,\jm)=1,\,s^-<s(\jm)^-}
M[s^-(\jm)^+,s(\jm)^-]  \right) \simeq X\tens Y\tens Z,
\end{align*}
where \eqn
X= \dtens_{\substack{d(\im,\jm)=1,\, s^-<s(\jm)^-,\, s(\jm)^+\le b}}M[s(\jm)^+,b\},\\
Y=\dtens_{\substack{d(\im,\jm)=1,\, s^-<s(\jm)^-,\, s(\jm)^+\le
b}} M[s^-(\jm)^+,s(\jm)^-] ,
 \\
Z=\dtens_{\substack{d(\im,\jm)=1,\, s^-<s(\jm)^-,\, s(\jm)^+> b}}
M[s^-(\jm)^+,s(\jm)^-]. \eneqn We have an epimorphism \eqn X\tens
Y&&\epito \dtens_{\substack{d(\im,\jm)=1,\, s^-<s(\jm)^-,\,
s(\jm)^+\le b}}
\bl M[s(\jm)^+,b\}\hconv  M[s^-(\jm)^+,s(\jm)^-]\br\\
&&\simeq \dtens_{\substack{d(\im,\jm)=1,\, s^-<s(\jm)^-,\,
s(\jm)^+\le b}}
 M[s^-(\jm)^+,b\}.
\eneqn by \cite[Lemma 3.2.22]{KKKO18}. On the other hand we have
$$Z\simeq
\dtens_{\substack{d(\im,\jm)=1,\, s^-<s(\jm)^-,\, s(\jm)^+> b}}
M[s^-(\jm)^+,b\}.
$$
Finally we obtain epimorphisms \eqn \cM^{\Vo} \tens \bS{s^-} \tens
\bS{s}
&&\epito X\tens Y\tens Z\\
&&\epito \Bigl(\dtens_{\substack{d(\im,\jm)=1,\\ s^-<s(\jm)^-,\,
s(\jm)^+\le b}}
 M[s^-(\jm)^+,b\}\Bigr)
\tens
\Bigl(\dtens_{\substack{d(\im,\jm)=1,\\ s^-<s(\jm)^-,\, s(\jm)^+> b}} M[s^-(\jm)^+,b\}\Bigr)\\
&&\simeq \dtens_{\substack{d(\im,\jm)=1,\, s^-<s(\jm)^-}}
 M[s^-(\jm)^+,b\}\simeq\cM^\Vi. \qedhere
\eneqn
\end{proof}
Hence we have shown the existence of an epimorphism \eqref{eq: epi},
and obtain ~\eqref{-thm b} in Lemma~\ref{thm:step1}. Thus we complete the proof
of Lemma~\ref{thm:step1}.

Now let us show \eqref{eq:adm}\;\eqref{item:b}.
\Lemma
For any $s\in\Kex$, we have
$$\de(\cM,\cMp) = 1.$$
\enlemma

\begin{proof}
Since the set of real modules $\{\cM[t]\}_{t\in \K}$ is a real
commuting family and $\cMp\simeq \cM^{\Vo} \hconv \cuspS{s^-}$,
Proposition~\ref{pro:subquotient}, Theorem~\ref{thm: admissible
chain commuting} and Proposition~\ref{prop: de1}  imply that
\begin{align*}
   \de(\cM,\cMp) & \le \de(\cM,\cM^{\Vo}) + \de(\cM,\bS{s^-}  ) \\
& =  \de(\cM ,\bS{s^-}) = 1.
\end{align*}
On the other hand,  we have
$$  \de(\D^{-1} \bS{s^{-}},\cM[s^-] \tens \cM^{\Vo})= 1 \quad \text{ and } \quad \de(\D^{-1} \bS{s^{-}}, \cM[s^+] \tens \cM^{\Vi} )=0,$$
by  Proposition~\ref{prop: lemma of KKOP20C},~\eqref{eq: Hi Ho} and
Lemma~\ref{lem:unique}. Since $\cM\hconv \cMp \simeq\cM[s^-] \tens
\cM^{\Vo} $ and $\cMp \hconv \cM\simeq\cM[s^+] \tens \cM^{\Vi}$, we
have
$$
\cM\hconv \cMp \not\simeq \cMp \hconv \cM.
$$
Thus, we have the desired result from Theorem~\ref{thm: KKKo15
main}~\eqref{i4}.
\end{proof}

Now we shall show \eqref{eq:adm}\;\eqref{item: admf}.

\Lemma\label{lem:5}
For $s<b$, $\cM'$ commutes with $\cusp_b$.
\enlemma
\Proof
Sublemma~\ref{lem: - socle} tells that
\eqn
\cM'\htens \cusp_s\simeq\dtens_{d(\im,\im_t)=1,t^-<s^-<t<s} M[t,b\}.
\eneqn
Hence,  $\cM'\htens \cusp_s$ commutes with $\cusp_b$.
Since $\D^{-1}\cusp_s$ commutes with $\cusp_b$,
we conclude that
$\cM'\simeq (\D^{-1}\cusp_s)\htens (\cM'\htens \cusp_s)$
commutes with $\cusp_b$.
\QED

\Prop\label{prop: admf} Let $s\le b$.
\bnum
\item \label{it: reqlsimple} $\cM'$ is real simple.
\item \label{it: commuting}
$\cM'$ commutes with $M[k,b\}$ for any $k\in\Z$
such that $k\le b$ and $k\not=s$.
\ee
\enprop

\Proof
We argue by induction on $b\ge s$.

If $b=s$, then
$\cM'=\cusp_{s^-}$ is real and commutes with $M[k,b\}$ for $k<b$ by Lemma \ref{lem:s=b}.

Now assume that $b>s$. Assuming that ~\eqref{it: reqlsimple} and ~\eqref{it: commuting} hold
when we replace $b$ with $b-1$,
we shall prove ~\eqref{it: reqlsimple} and ~\eqref{it: commuting}.

Set
\eqn
Y&&=\hs{-5ex}\dtens_{\substack{  t\in[a,b-1];\;d(\im_s,\im_t)=1, \\
\akew[8ex]s^-<t^-<s<t\le b-1}} M[t,b-1\},\\
X&&=Y\hconv \cusp_{s^-}.
\eneqn
Then by the induction hypothesis,
$X$ is real and commutes with $M[k,b-1\}$ if $k\le b$ and $k\not=s$.
(Here, we understand $M[k,b-1\}=\one$ if $k=b$.)
\eqn
&&\cM^\Vo\simeq \bc
\cusp_b\htens Y&\text{if $d(\im_s,\im_b)=1$ and $s^-<s(\im_b)^-$,}\\
Y&\text{otherwise.}\ec
\eneqn
Since $(\cusp_b,\cusp_{s^-})$ is unmixed we have
\eqn
&&\cM' \seteq   \cM^\Vo\hconv \bS{s^-} \simeq \bc
\cusp_b\htens X&\text{if $d(\im,\im_b)=1$ and $s^-<s(\im_b)^-$,}\\
X&\text{otherwise.}\ec
\eneqn

We have also
\eqn
&&M[k,b\}\simeq \bc
\cusp_b\htens M[k,b-1\}&\text{if $\im_k=\im_b$,}\\
M[k,b-1\}&\text{otherwise.}\ec
\eneqn

We have
$\de(\dual \cusp_b,X)=\de(\dual \cusp_b,M[k,b-1\})=0$.
By Lemma~\ref{lem:5},
we have
$$\de(\cusp_b,\cM')=0 \quad \text{and} \quad \de(\cusp_b,M[k,b\})=0.$$
Now we shall apply Lemma~\ref{lem:3} with $M_1=X$, $M_2=M[k,b-1\}$,
and $L_1,L_2=\cusp_b$ or $\one$. Then we have
$$L_1\htens M_1\simeq\cM'  \quad \text{and} \quad  L_2\htens M_2\simeq M[k,b\}.$$
Since $M_1$ and $M_2$ are real and commute by the induction hypothesis,
$L_1\htens M_1\simeq\cM'$ and $L_2\htens M_2\simeq M[k,b\}$
are real and commute.
\QED

Thus we complete the proof of Theorem~\ref{thm:triangle}.


\subsection{The cluster algebra structure on $K(\Ca_\g^{<\xi})$}\label{subsec:Cminus}
We take a quantum affine algebra $\uqpg$ and the associated data
$(\Dynkin,\sigma)$. We freely use terminologies in \S\,\ref{sec:
subcategories}.

{\em In the sequel, we choose an arbitrary $\hI_\g$, and we consider
only $\rmQ$-data $\Qd$ and height
functions $\xi$ such that $\hI_\g=\hI_\Qd$} (see Remark~\ref{rmk: Up to constant}).
Note that the choice of $\hI_\g$ determines $\catCO$. Indeed,
$\catCO$ is the smallest full subcategory of $\catg$ which contains
all $V(\im,p)$ with $(\im,p)\in\hI_\g$ and is stable under taking
tensor products, subquotients and extensions.

Recall that $\hIg^{<\xi}\seteq\st{(\im,p)\in\hIg\mid p<\xi_\im}$ and
the category $\Ca_\g^{< \xi}$ is the smallest full subcategory of
$\catg$ containing all $V(\im,p)$ with $(\im,p)\in\hI_\g^{<\xi}$ and
stable under taking tensor products, subquotients and extensions
(see~\eqref{eq: V(i;p)} and~\eqref{eq: hI<xi}).

\Def
We say that an admissible sequence $\hw=\seq{(\im_k,p_k)}_{k\in\Z}$
in $\hI_\g$ is {\em $\xi$-adapted} if
$\hI_\g^{<\xi}=\set{(\im_k,p_k)}{k\in\Z_{\le0}}$.
\edf
As seen in \S\,\ref{subsec:subcategory}, there exists a
$\xi$-adapted admissible sequence.

\smallskip

\Def We define the quiver $\HL$
whose set of vertices is $\If\times \Z$
and the arrows  are assigned as follows (cf.~\eqref{eq: Psi_g}):
\begin{align*}
 (\im,x) \to (\jm,y) \quad & \  \quad \text{{\rm (i)} if $d(\im,\jm)=1$ and $x=y-2\fd_\jm+\min(\fdi,\fdj)$,}\\
                           & \text{or {\rm (ii)}  $\im=\jm$ and $x=y+2\fdi$. }
\end{align*}
We denote by $\HL^{<\xi}$ the quiver on $\hI_\g^{<\xi}$ induced by
$\HL$.
\edf

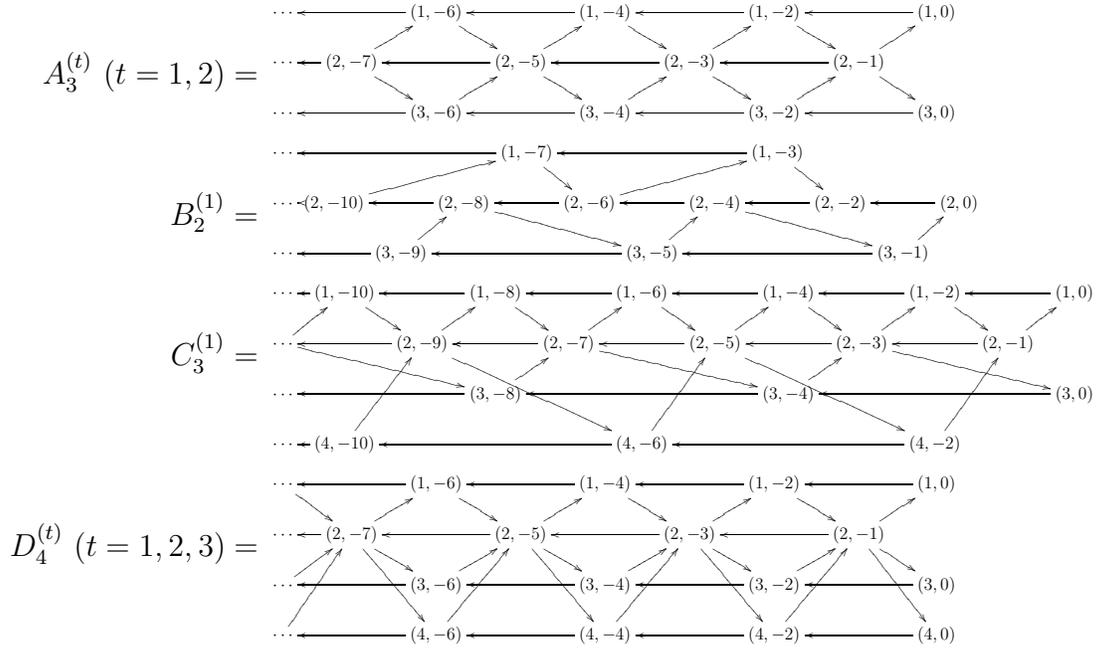
\begin{figure}
\begin{gather*}
\begin{split}
{\tiny A^{(t)}_{3}} \ (t=1,2) &=
\raisebox{2.3em}{\scalebox{0.55}{\xymatrix@C=3ex@R=3ex{
\cdots  &&(1,-6)  \ar[dr] \ar[ll] && (1,-4)  \ar[dr]\ar[ll] && (1,-2)  \ar[dr]\ar[ll] &&  \ar[ll] (1,0)\\
\cdots &\ar[l]\ar[dr]\ar[ur](2,-7)&&(2,-5)  \ar[ur] \ar[dr] \ar[ll]&&  \ar[ur] \ar[dr]\ar[ll] (2,-3) &&  \ar[ur] \ar[ll]\ar[dr] (2,-1) \\
\cdots  &&(3,-6)  \ar[ur] \ar[ll] && (3,-4)  \ar[ur]\ar[ll] &&
(3,-2)  \ar[ur]\ar[ll] &&  \ar[ll]  (3,0)
}}} \\
{\tiny B^{(1)}_{2}} &=
\raisebox{2.3em}{\scalebox{0.55}{\xymatrix@C=0.1ex@R=3ex{
\cdots&&&& \ar[dr] \ar[llll](1,-7)  &&&& \ar[dr] (1,-3)  \ar[llll]  \\
\cdots&(2,-10) \ar[l]\ar[urrr] &&   (2,-8) \ar[ll] \ar[drrr] && \ar[urrr]\ar[ll] (2,-6) && \ar[ll]\ar[drrr] (2,-4) && (2,-2)\ar[ll]  && \ar[ll](2,0)  \\
\cdots&& \ar[ll]\ar[ur] (3,-9)   &&&& \ar[ur] (3,-5) \ar[llll] &&&&
\ar[ur] (3,-1) \ar[llll]
}}}\\
{\tiny C^{(1)}_{3}} &=
\raisebox{2.3em}{\scalebox{0.55}{\xymatrix@C=1.5ex@R=3ex{
\cdots&\ar[l](1,-10) \ar[dr]&& \ar[dr]\ar[ll] (1,-8) && \ar[dr]\ar[ll] (1,-6) && \ar[dr]\ar[ll] (1,-4) && \ar[dr] (1,-2) \ar[ll] && \ar[ll] (1,0)\\
\cdots\ar[ur]\ar[drrr]&&\ar[ll] \ar[ur] (2,-9) \ar[ddrrr]&&  \ar[ur]\ar[ll] (2,-7) \ar[drrr]\ar[ll] && \ar[ur]\ar[ddrrr]\ar[ll]  (2,-5)  && \ar[ur]\ar[ll]\ar[drrr]  (2,-3)  && \ar[ur] \ar[ll](2,-1) \\
\cdots&&& \ar[lll](3,-8)  \ar[ur]&&&& (3,-4) \ar[llll]\ar[ur]  &&&& \ar[llll] (3,0)   \\
\cdots&\ar[l](4,-10) \ar[uur]  &&&& (4,-6)\ar[uur] \ar[llll]  &&&&
(4,-2)\ar[uur] \ar[llll]
}}} \\
{\tiny D^{(t)}_{4}} \ (t=1,2,3) &=
\raisebox{2.3em}{\scalebox{0.55}{\xymatrix@C=3ex@R=3ex{
  \cdots \ar[dr]  && (1,-6)  \ar[ll]\ar[dr] && (1,-4)  \ar[ll]\ar[dr] && \ar[dr]\ar[ll] (1,-2)  && (1,0)\ar[ll]\\
\cdots & \ar[ur]\ar[l] \ar[ddr]\ar[dr] (2,-7) && \ar[ll] \ar[ur] \ar[ddr]\ar[dr] (2,-5) && \ar[ll] \ar[ur] \ar[ddr]\ar[dr] (2,-3) && \ar[ll] \ar[ur] \ar[ddr]\ar[dr] (2,-1) \\
 \cdots \ar[ur] && (3,-6) \ar[ll] \ar[ur] && (3,-4) \ar[ll] \ar[ur] &&  \ar[ur]  (3,-2) \ar[ll] && (3,0)\ar[ll]\\
\cdots \ar[uur] && (4,-6) \ar[ll]  \ar[uur] && (4,-4) \ar[ll] \ar[uur] &&  \ar[uur] (4,-2) \ar[ll] && \ar[ll] (4,0) \\
}}}
\end{split}
\end{gather*}
\caption{Some examples of the quivers $\HL^{<1}$.}
\label{fig:quiver_examples Uppsi}
\end{figure}

The following proposition gives a relation between
 $\HL$ and $\GLS$
given in \eqref{quiver:Omegaminus}.

\Prop\label{prop:HL=GLS} Let $\hw=\seq{(i_s,p_s)}_{s\in\Z}$ be an
admissible sequence in $\hIg$. Then the quiver $\GLS(\hw)$ is
isomorphic to the quiver $\HL$ by the bijection $\Z\isoto \hI_\g $
given by $\Z\ni s\mapsto (\im_s,p_s)\in \hI_\g $. \enprop

\Proof
Set $(\im,x)=(\im_s,p_s)$ and $(\jm,y)=(\im_t,p_t)$. By
Definition~\ref{def: admissible sequence} \eqref{it:ssp}, it is
obvious that
$$\im=\jm \ \ \text{ and } \ \ x=y+2\fd_{\im} \Longleftrightarrow
t=s^-.$$

Hence it remains to prove the equivalence under the hypothesis
$d(\im,\jm)=1$:
$$p_s=p_t-2\fd_\jm+\min(\fd_\im,\fd_\jm)
\Longleftrightarrow s^-<t^-<s<t.
$$
Since
\begin{align*}
p_{t^-}+\min(\fd_{\im_{t^-}},\fd_{\im_s})
=p_t-2\fd_{\im_t}+\min(\fd_{\im_{t}},\fd_{\im_s}),
\end{align*}
it is enough to show (by replacing $t$ with $t^+$)
$$p_s=p_{t}+\min(\fd_{\im_{t}},\fd_{\im_s}) \qt{ if and only if }s^-<t<s<t^+,$$
under the hypothesis that $d(\im_s,\im_t)=1$.

\snoi If $d(\im_s,\im_t)=1$ and $s^-<t< s < t^+$, then
Definition~\ref{def: admissible sequence}~\eqref{it:tsts} implies
$p_s=p_{t}+\min(\fd_{\im_{t}},\fd_{\im_s})$.

\snoi Conversely assume that
$p_s=p_{t}+\min(\fd_{\im_{t}},\fd_{\im_s})$. Then  we have \eqn
p_{s^-}=p_s-2\fdi=p_t-2\fdi+\min(\fdi,\fdj)<p_t<p_s=p_t+\min(\fdi,\fdj)
<p_t+2\fdj=p_{t^+}. \eneqn Thus we have $s^-<t^-<s<t$ by
Lemma~\ref{lem:tendency}\;\eqref{it: condi}. \QED

\Def\label{def:Mxi} For $(\im,p)\in\hI_\g$, we define
\begin{align*}
\M^{< \xi}_{(\im,p)}=\hd\bl V(\im, p+ r\cdot(2\fdi))\tens
V(\im,p+(r-1)(2\fdi))
\tens \cdots\tens V(\im, p)\br,
\end{align*}
where $r$ is the largest integer such that
$p+r\cdot(2\fdi)<\xi_\im$. It is a KR-module.
\end{definition}
\begin{remark}
Let us take a $\xi$-adapted admissible sequence $\hw$ in $\hIg$
and  $a\in\Z$ such that $(\im,p)=(\im_{a},p_{a})$. Then we have
\begin{align} \label{eq: def Mi-k}
\M^{< \xi}_{(\im,p)}  \simeq M^{\hw}[a,0\},
\end{align}
\end{remark}

For $t\in\Z$, we set $\hIg^{<t}\seteq\st{(\im,p)\in\hIg\mid p<t}$,
and let $\catg^{<t}$ be the smallest full subcategory of $\catg$
which contains $V(\im,p)$ for $(\im,p)\in \hIg^{<t}$ and stable
under taking tensor products, subquotients and extensions. Then
there exists a unique height function $\xi^{t}$ such that
$\hIg^{<t}=\hIg^{<\xi^{t}}$, namely, $\xi^{t}$ satisfying
$(\xi^{t})_\im\in[t,t+2\fdi-1]$ for all $\im\in\If$.
 We write $\M^{<t}_{(\im,p)}$ and $\HL^{<t}$ instead of
$\M^{<\xi^{t}}_{(\im,p)}$ and $\HL^{<\xi^{t}}$. Hence we have
\eq\M^{<t}_{(\im,k)}= \hd\Bigl( V\bl \im, p+ r\cdot(2\fdi)\br\tens
V\bl\im,p+(r-1)(2\fdi)\br
\tens \cdots\tens V\bl\im, p\br\Bigr),\label{def:M<t} \eneq where $r$ is the largest integer
such that $p+r\cdot(2\fdi)<t$.

In \cite{HL16}, Hernandez-Leclerc proved that the seed arising from
$\{\M^{<t}_{(\im,p)} \}_{(\im,p)\in\hIg^{<t}}$ and $\HL^{<t}$ gives
a cluster algebra structure on the Grothendieck ring
$K(\Ca_\g^{<t})$ as follows: set $z_{\im,p}=[ \M^{<t}_{(\im,p)} ]\in
K(\Ca_\g^{<t})$ for $(\im,p)  \in  \hI^{<t}_\g  $ and let
$\BHL^{<t}$ be the exchange matrix associated with $\HL^{<t}$.

\begin{theorem}[{\cite[Theorem 5.1]{HL16}}]\label{th:HL16}
Take an arbitrary $\hI_\g$ and an integer $t$.  Let
$\Seed^{<t}\seteq   \bigl( \st{z_{\im,p}}_{  (\im,p) \in \hI^{<t}_\g
},\BHL^{<t} \bigr) $ be a seed in $K(\Ca_\g^{<t})$ with the empty
set of frozen variables. Then, we have
$$\A(\Seed^{<t}) \simeq K(\Ca_\g^{<t}).$$
\end{theorem}

\vskip 2em
\section{Monoidal categorification of cluster algebras} \label{sec: main}

In this section, we shall state and prove our main theorem.  We fix
$\hI_\g$ associated with $\catCO$. Recall that
 $\catCO$ is the smallest full subcategory of $\catg$
which contains $V(\im,p)$ ($(\im,p)\in\hI_\g$) and is stable under
taking tensor products, subquotients and extensions. We only treat
height functions $\xi$ such that $(\im,\xi_\im)\in\hIg$ for every
$\im\in\If$, and admissible sequences $\hw$ in $\hI_\g$.

\subsection{Statement of the main theorem} \label{sec: Statement of main theorem}
The purpose of this section is to prove the following main theorem
of this paper.
\Th\label{th:Main}
Let  $(\Dd,\hhw)$ be a  \pbw with $\Dd=\Dd_\Qd$ for some $\rmQ$-datum $\Qd$ and an arbitrary reduced expression $\tw_0$ of $w_0$,
and let $\frakC$ be an admissible chain of $i$-boxes with a
range $[a,b]$. Then, there is an exchange matrix
$\tB(\frakC)=(b_{s,t})_{(s,t)\in\K(\frakC)\times\Kex(\frakC)}$
satisfying the following properties. \bna
\item\label{def:seedC}
$\seed(\frakC)\seteq\bl\cM[](\frakC),\tB(\frakC);\K(\frakC),\Kex(\frakC)\br$
is a \ca monoidal seed in $\Ca_\g^{[a,b],\Dd,\hhw}$
\rmo see \eqref{def:chainseed} for the notations in $\seed(\frakC)$\rmf.
\item $\A\bl[\seed(\frakC)]\br\simeq K(\Ca_\g^{[a,b], \Dd,\hhw})$,
\ee Namely, the category $\Ca_\g^{[a,b],\Dd,\hhw}$ provides a
$\Uplambda$-monoidal categorification of the cluster algebra
$K(\Ca_\g^{[a,b],\Dd,\hhw})$ with the initial monoidal seed
$\seed(\frakC)$.
\enth
Note that $K(\Ca_\g^{[a,b],\Dd,\hhw})$ is the polynomial algebra generated by
  $\st{[\cusp_s]\mid s\in[a,b]}$  as an algebra (see \cite[\S\,6.3]{KKOP20C}).

This section is devoted to the proof of this theorem. Since the
proof is intricate, we afford the plan of the proof.

For $t\in\Z$ and a height function $\xi$, we define the monoidal
seed $\seed^{<t}$ in $\catg^{<t}$ and the monoidal seed
$\seed^{<\xi}$ in $\catg^{<\xi}$ as follows: \eqn
\seed^{<t}&\seteq&\bl\st{\M^{<t }_{(\im,p)}}_{(\im,p)\in\hI_\g^{<t}},
\BHL^{<t};\hI_\g^{<t},\hI_\g^{<t}\br,\\
\seed^{<\xi}&\seteq&\bl\st{\M^{<
\xi}_{(\im,p)}}_{(\im,p)\in\hI_\g^{<\xi}},
\BHL^{<\xi};\hI_\g^{<\xi},\hI_\g^{<\xi}\br, \eneqn where $\BHL^{<t}$
(resp.\ $\BHL^{<\xi}$) is the exchange matrix associated with the
quiver on $\hI_\g^{<t}$ (resp.\ $\hI_\g^{<\xi}$) induced by $\HL$ (see
\eqref{def:M<t} for $\M^{< t}_{(\im,p)}$ and
Definition~\ref{def:Mxi} for $\M^{< \xi}_{(\im,p)}$).

\smallskip
Remark that, once $\hI_\g$ is fixed,
$\seed^{<t}$ depends only on $t\in\Z$ and $\seed^{<\xi}$ depends only on $\xi$.
The monoidal seed $\seed^{<t}$ is a special case of monoidal seeds
$\seed^{<\xi}$. Indeed, if $\xi$
is a unique height function  satisfying the following condition:
\eq &&\text{$\xi{}_\im\in[t,t+2\fdi-1]$ for every
$\im\in\If$,}\label{cond:standard} \eneq then we have
$\seed^{<\xi}=\seed^{<t}$. We denote by $\xi^t$ the height function satisfying~\eqref{cond:standard}.

Recall that $\hw=\seq{(\im_s,p_s)}_{s\in\Z}$ is $\xi$-adapted if
$\hIg^{<\xi}=\st{(\im_s,p_s)\mid s\le0}$.
We have
$$\seed^{<\xi}=\seed_-^{[-\infty,0],\Dd_\Qd,\hhw}=\seed(\frakC^{[-\infty,0],\hw}),$$
 where $\Qd$ is the Q-datum associated with $\xi$, $\hhw$
is associated with a $\Qd$-adapted reduced expression $\tw_0$ of $w_0$,
$\hw$ is a $\xi$-adapted admissible sequence in $\hIg$ and
$\frakC^{[-\infty,0],\hw}_-$ is the admissible chain of $i$-boxes
  $([1-k,0\})_{1\le k}$ (see \eqref{def:seedab} and Theorem~\ref{th:Main}\;\eqref{def:seedC} for the notations).

\medskip
Consider  the following statements, where
$\xi$ is a height function:
$$\left\{
\parbox{80ex}{
\ben[(I)${}_\xi$]
\renewcommand{\labelenumi}{(\Roman{enumi})\akew[1.2ex]}
\item \label{item:0} $\seed^{<t}$ is a \ca monoidal seed
in $\catg^{<t}$ for any $t\in\Z$,
\vs{1ex}
\renewcommand{\labelenumi}{(\Roman{enumi})${}_\xi$}
\item \label{item:I}$\seed^{<\xi}$ is a \ca monoidal seed
in $\catg^{<\xi}$,
\vs{1ex}
\item \label{item:II}
for  the $\rmQ$-datum $\Qd$ associated with $\xi$ and any reduced expression $\tw_0$ of $w_0$, the seed $\seed_-^{[-\infty,0],\Dd_\Qd,\hhw}$ is a \ca monoidal seed in $\catg^{<\xi}$,
\vs{1ex}
\item \label{item:III} for the $\rmQ$-datum $\Qd$ associated with $\xi$,
any reduced expression $\tw_0$ of $w_0$, 
and any admissible chain $\frakC$ of $i$-boxes,
$\seed(\frakC)=(\cM[](\frakC),\tB;\K(\frakC),\Kex(\frakC)\br$ is a
\ca monoidal seed in 
$\Ca_\g^{[a,b],\Dd_\Qd,\hhw}$  for a suitable choice of
exchange matrix $\tB$.
\vs{1ex}
\ee } \right.
$$

The first step of the proof of main Theorem~\ref{th:Main}  is to show \eqref{item:0}  by using the result of
Hernandez-Leclerc (Theorem~\ref{th:HL16}).
Hence, \mbox{\eqref{item:I}}${}_{\xi}$ holds as soon as
$\xi$ satisfies condition \eqref{cond:standard}.

\snoi
Then we will prove
\begin{align*}
\mbox{\eqref{item:I}}{}_\xi\Longrightarrow
\mbox{\eqref{item:II}}{}_\xi.
\end{align*}
in Proposition~\ref{prop: mutation equivalent}, by showing the monoidally mutation equivalence between $\seed^{< \xi}$
and  $\seed_-^{[-\infty,0],\Dd_\Qd,\hhw}$.
Then we prove
\begin{align*}
\mbox{\eqref{item:II}}{}_\xi\Longrightarrow
\mbox{\eqref{item:III}}{}_\xi.
\end{align*}

Now, we can see easily that for any $\xi$ there exists a
$\xi$-adapted reduced expression $\tw=s_{\im_1}\cdots s_{\im_r}$
such that $\xi'\seteq  s_{\im_r}\cdots s_{\im_1} \xi$ satisfies
condition \eqref{cond:standard}. On the other hand,
\eqref{item:I}${}_{\xi'}$ holds. 
Hence, in order to see \eqref{item:I}${}_\xi$ for an arbitrary
$\xi$, it is enough to show
\begin{align*}
\mbox{\eqref{item:I}}{}_{s_\im\xi}\Longrightarrow
\mbox{\eqref{item:I}}{}_\xi \qt{if $\im$ is a sink of $\xi$.}
\end{align*}

It is performed  in Proposition~\ref{prop:sixi}.

Theorem~\ref{th:Main} is equivalent to saying that
\eqref{item:III}${}_\xi$ holds for any $\xi$.

\subsection{Monoidal categorification by $\Ca_\g^{<t}$}
 In order to achieve the main goal of this section,
we shall first prove the conjecture suggested by Hernandez-Leclerc
for untwisted affine $\g$.

Recall that $\hI_\g^{<t}\seteq \st{(\im,p)\in\hI_\g\mid p<t}$, and
$\cat^{<t}_\g$ is the smallest full subcategory of $\catg$
containing $V(\im,p)$ for $(\im,p)\in\hI_\g^{<t}$ and stable under
taking tensor products, subquotients and extensions.

\Th[{cf.\ \cite[Conjecture 5.2]{HL16}}]  \label{thm: main cg-} The
category $\Ca_\g^{<t}$ provides a $\Uplambda$-monoidal
categorification of the cluster algebra $K(\Ca_\g^{<t})$ associated
to the initial \ca monoidal seed $\seed^{<t}$. \enth

\Proof
Take an admissible sequence $\hw=\seq{(\im_s,p_s)}_{s\in\Z}$ in $\hI_\g$ such that
$\hI_\g^{<t}\seteq\st{(\im_k,p_k)\mid k\le 0}$. Then we have
$\Ca_\g^{[-\infty,0],\hw}=\Ca_\g^{<t}$ and \eq
&&\seed^{{<t}}=\seed_-^{[-\infty,0],\,\hw} \eneq by
Proposition~\ref{prop:HL=GLS}. Hence $\seed^{{<t}}$ is a  \Lad
monoidal seed by Theorem~\ref{thm:triangle}.  
Since $\A([\seed^{<t}])\simeq K(\Ca_\g^{<t})$ by
Theorem~\ref{th:HL16} due to Hernandez-Leclerc, the assertion
follows from Theorem~\ref{th:main KKOP19C}. \QED

 Thus we have obtained \eqref{item:0}.

\subsection{Mutation equivalence}

\begin{definition}\label{def:mutationequiv}
Let  $\seed=(\{M_i\}_{i\in \K},\widetilde B; \K,\Kex)$ and
$\seed'=(\{M'_i\}_{i\in \K'},\widetilde{B}'; \K',(\K')^\ex)$ be admissible monoidal seeds
in $\shc$.
\ben
\item We say that \emph{$\seed'$ is monoidally mutated from $\seed$}
if the following condition is satisfied: For any finite subset $\mathsf{J}$ of $\K'$,
there exist
\bnum
\item a finite sequence $(k_1,k_2,\ldots,k_r)$ in $\Kex$ such that
$\mu_{k_s} \circ \cdots \circ \mu_{k_1} (\seed)$ is an admissible monoidal seed
for each $1 \le s \le r$,
\item an injective map $\upsigma\col\sfJ \to \K$, depending on the choice of $\sfJ$, such that
\bna
\item $ \upsigma(\mathsf{J}^\ex) \subset \Kex$, where $\mathsf{J}^{\ex} \seteq \mathsf{J} \cap (\K')^\ex$, 

\item $M'_{j}=\upmu(M)_{\upsigma(j)}$  for all $j \in \sfJ$,
\item $ (\tB')_{(i,j)}=\upmu(\tB)_{\upsigma(i),\upsigma(j)}$
for any $(i,j)\in \sfJ \times {\sfJ}^\ex$,
\ee
where $\upmu \seteq  \mu_{k_r} \circ \cdots \mu_{k_1}$.
\ee
\item  We say that the admissible monoidal seeds $\seed$ and $\seed'$ are \emph{\me} if $\seed'$ is monoidally mutated from $\seed$ and $\seed$ is also monoidally mutated from $\seed'$.
\ee
\end{definition}

For a \pbw $(\Dd,\hhw)$, let $\frakC^{[-\infty,0],\Dd,\hhw}_-$ be the admissible chain of $i$-boxes
$([1-k,0\})_{1\le k}$ and recall that
\begin{align}\label{eq: theseed}
\seed_-^{[-\infty,0],\Dd,\hhw}\seteq\bl\{M(\frakC^{[-\infty,0],\Dd,\hhw}_-),\BGLS^{\Dd,\hhw}(\hhw) ;[-\infty,0],\,[-\infty,0] \br.
\end{align}

When $\Qd$ is the $\rmQ$-datum with the height function $\xi^t$ (see~\eqref{cond:standard}) and $\tw_0$ is a $\Qd$-adapted reduced expression of $w_0$,
we have
$$\seed^{{<t}}=
\seed_-^{[-\infty,0],\Dd_\Qd,\hhw}$$
and it is a completely \Lad monoidal seed by Theorem~\ref{thm: main cg-}.

\begin{proposition} \label{prop: mutation equivalent}
Let $(\Dd,\hhw)$ be a \pbw and $\tw_0'$ be another reduced expression of $w_0$. Then the \Lad seeds
$$
\seed_-^{[-\infty,0],\Dd,\hhw} \text{ and }  \seed_-^{[-\infty,0],\Dd,\hhw'} \text{ are \me}.
$$
\end{proposition}

This subsection is devoted to prove the above proposition.

\begin{remark}\label{rmk: braid move cuspidal} Recall the notion of the affine cuspidal modules $\bS{k}^{\Dd,\hhw}$ in Section~\ref{subsec:pbw}.
\bna
\item \label{it: braid}
When $\tw_0'=s_{\jm_1} \ldots s_{\jm_\ell}$ is obtained from
$\tw_0=s_{\im_1} \ldots s_{\im_\ell}$ by a single \emph{braid move};
i.e., for some $k$ such that  $1  <  k <\ell$, we have
$$  \text{$d(\im_{k-1},\im_{k})=1$, $ \im_{k \pm 1} = \jm_k$, $\jm_{k \pm 1} = \im_k$
, and $\im_s = \jm_s$ for any $s \not\in \{ k-1,k, k+1\}$,}$$
the affine cuspidal modules are related by
\begin{align}\label{eq: cusp rel}
\bS{s}^{\Dd,\hhw'} \simeq \bc
\bS{s+2}^{\Dd,\hhw}&\text{if $s\equiv_\ell k-1$,}\\
\bS{s+1}^{\Dd,\hhw}\hconv \bS{s-1}^{\Dd,\hhw}&\text{if $s\equiv_\ell k$,}\\
\bS{s-2}^{\Dd,\hhw}&\text{if $s\equiv_\ell k+1$,}\\
\bS{s}^{\Dd,\hhw}&\text{otherwise.}\ec
\end{align}
Here ~\eqref{eq: cusp rel} follows from 
\cite[Proposition 5.9]{KKOP20C}. Indeed, by this proposition, we can reduce to the case where $k=2$, and it is easy to check this case.
Note that $\bS{s}^{\Dd,\hhw} \simeq \bS{s-1}^{\Dd,\hhw}\hconv \bS{s+1}^{\Dd,\hhw}$
if $s\equiv_\ell k$.

\item \label{it: commutation}
When $\tw_0'=s_{\jm_1} \ldots s_{\jm_\ell}$ is obtained from $\tw_0=s_{\im_1} \ldots s_{\im_\ell}$ by a single \emph{commutation move}; i.e.,
there exists a $1  <  k \le\ell$ such that
$$  d(\im_{k-1},\im_{k})>1, \  \im_k=\jm_{k-1}, \  \im_{k-1}=\jm_{k} \text{ and } \im_s = \jm_s  \text{ for } s \not\in \{ k-1, k\},$$
we have
$$
\bS{s}^{\Dd,\hhw'} \simeq
\bc
\bS{s+1}^{\Dd,\hhw}&\text{if $s\equiv_\ell k-1$,}\\
\bS{s-1}^{\Dd,\hhw}&\text{if $s\equiv_\ell k$,}\\
\bS{s}^{\Dd,\hhw}&\text{otherwise.}\ec
$$
\item \label{it: moves}   All reduced expressions of $w_0$ are connected via braid moves and commutation moves.
\ee
\end{remark}

Set $\mathsf{M}^{\Dd,\hhw}_s\seteq M^{\Dd,\hhw}[s,0\}$ for a reduced expression $\tw_0$ of $w_0$ and $s\in\Z_{\le0}$.

\begin{lemma}
Let $(\Dd,\hhw)$ be a \pbw and $\tw_0'=s_{\jm_1} \cdots s_{\jm_\ell}$ be a reduced expression of $w_0$ obtained from $\tw_0=s_{\im_1} \cdots s_{\im_\ell}$ by a single commutation move as in {\rm Remark~\ref{rmk: braid move cuspidal} ~\eqref{it: commutation}}. %
Then we have
$$
\mathsf{M}^{\Dd,\hhw'}_s = \bc
\mathsf{M}^{\Dd,\hhw}_{s-1} &\text{ if $s \equiv_\ell k$}, \\
\mathsf{M}^{\Dd,\hhw}_{s+1} &\text{ if $s\equiv_\ell k-1$}, \\
\mathsf{M}^{\Dd,\hhw}_{s} & \text{ otherwise}.
\ec
$$
\end{lemma}

\begin{proof}
This is a direct consequence of Remark~\ref{rmk: braid move cuspidal} ~\eqref{it: commutation}
\end{proof}

\begin{lemma}\label{lem: mutataion quiver1}
Let us keep the notations in the previous lemma. Then the quivers $\GLS^{[-\infty,0]}(\hhw)$ and $\GLS^{[-\infty,0]}(\hhw')$ are isomorphic to each other under the index change $\Upphi\col[-\infty,0] \to [-\infty,0]$ given by
$$
\Upphi(s) = \bc
s-1 &\text{ if $s \equiv_\ell k$}, \\
s+1 &\text{ if $s \equiv_\ell k-1$}, \\
s & \text{ otherwise}.
\ec
$$
\end{lemma}

 The proof is straightforward.

\begin{lemma} \label{lem: mutataion variable}
Let $(\Dd,\hhw)$ be a \pbw and let $\tw_0'=s_{\jm_1} \cdots s_{\jm_\ell}$ be a reduced expression of $w_0$ obtained from $\tw_0=s_{\im_1} \cdots s_{\im_\ell}$ by a single braid move as in {\rm Remark~\ref{rmk: braid move cuspidal} ~\eqref{it: braid}}. %
Then we have
\eq
\mathsf{M}^{\Dd,\hhw'}_s\simeq
\bc
\mathsf{M}^{\Dd,\hhw}_{s+1}&\text{if $s\equiv_\ell k-1$,}\\
\mathsf{M}^{\Dd,\hhw}_{s-1}&\text{if $s\equiv_\ell k$,}\\
\bl\mathsf{M}^{\Dd,\hhw}_s\br'&\text{if $s\equiv_\ell k+1$,}\\
\mathsf{M}^{\Dd,\hhw}_{s}&\text{otherwise,}
\ec
\eneq
where
$\bl \mathsf{M}^{\Dd,\hhw}_{s} \br' $ denotes the mutation  $\mu_{s}\bl \mathsf{M}^{\Dd,\hhw}\br_{s} $ of $\mathsf{M}^{\Dd,\hhw}_{s}$ described in~\eqref{eq: M'}.
\end{lemma}
Note that, when $s\equiv_\ell k+1$,
there exists at most one \emph{vertically} outgoing arrow starting from $s$ in $\GLS(\hhw)$, which is
$s \to (s-1)^+_\im$ if $(s-1)_\im^+\le0$.
Since $(s)_\im^-=s-2$ (see~\eqref{eq: a(im)^pm_jm} below), we have
$$\bl \mathsf{M}^{\Dd,\hhw}_{s} \br'\simeq\mathsf{M}^{\Dd,\hhw}_{(s-1)_\im^+}\hconv
\mathsf{S}_{s-2}^{\Dd,\hhw}\qt{(when $s\equiv_\ell k+1$).}$$

\begin{proof}
We shall argue by the descending induction on $s$.
If $k+1-\ell<s\le 0$, the assertion is obvious by~\eqref{eq: cusp rel}.

In the course of the proof, we denote
\begin{equation} \label{eq: a(im)^pm_jm}
\begin{aligned}
&(a)_\im^+=\max\st{k>a\mid \im_k=\im_a}, \qquad
(a)_\jm^+=\max\st{k>a\mid \jm_k=\jm_a}  \\
&(a)_\im^-=\max\st{k<a\mid \im_k=\im_a}, \qquad
\end{aligned}
\end{equation}
Note that we have 
$$ \mathsf{M}^{\Dd,\hhw}_{s}\simeq\mathsf{M}^{\Dd,\hhw}_{(s)_\im^+}\hconv \mathsf{S}_{s}^{\Dd,\hhw}\qtq\mathsf{M}^{\Dd,\hhw'}_{s}\simeq\mathsf{M}^{\Dd,\hhw'}_{(s)_\jm^+}\hconv \mathsf{S}_{s}^{\Dd,\hhw'}.$$
Here we understand
$$ \mathsf{M}^{\Dd,\hhw}_{s}\simeq\mathsf{M}^{\Dd,\hhw'}_{s}\simeq \one
\qt{if $s>0$.}$$

Note that if we set $t=(s)_\jm^+>s$, then we have
$$ \mathsf{M}^{\Dd,\hhw'}_{t}
=\bc
 \mathsf{M}^{\Dd,\hhw}_{t+1}&\text{if $t\eqv k-1$,}\\
 \mathsf{M}^{\Dd,\hhw}_{t-1}&\text{if $t\eqv k$,}\\
 \mathsf{M}^{\Dd,\hhw}_{t}&\text{if $t\not\eqv k,k\pm1$.}
\ec$$

\bnum
\item \label{item:sk}
Let us consider the case $s\equiv_\ell k$. Then $\jm_s=\im_{s+1}$.
Set $t=(s)_\jm^+$.
First, let us show
\eq
\mathsf{M}^{\Dd,\hhw'}_t\simeq\mathsf{M}^{\Dd,\hhw}_{(s+1)_{\im}^+}.
\label{eq:k} \eneq

\bna
\item Case $t\eqv k-1$. Then $\jm_s=\im_{s+1}=\jm_t=\im_{t+1}$, and hence
$(s+1)_\im^+=t+1$.
Hence we have \eqref{eq:k}.
\item Case $t\eqv k$. Then $\jm_s=\im_{s+1}=\jm_t=\im_{t-1}$, and hence
$(s+1)_\im^+=t-1$. Therefore we have \eqref{eq:k}.
\item The case $t\eqv k+1$ does not occur.
\item Case $t\not\eqv k,k\pm1$. Then $t=(s+1)_\im^+$ and
we have \eqref{eq:k}.
\ee
Thus we complete the proof of \eqref{eq:k}.

Then, we have
\eqn 
&&\mathsf{M}^{\Dd,\hhw'}_s\simeq\mathsf{M}^{\Dd,\hhw'}_{t}
\hconv (  \mathsf{S}_s^{\Dd,\hhw'}  )
\simeq \mathsf{M}^{\Dd,\hhw}_{(s+1)_\im^+} \hconv (  \mathsf{S}_{s+1}^{\Dd,\hhw} \hconv \mathsf{S}_{s-1}^{\Dd,\hhw}  ) \simeq \mathsf{M}^{\Dd,\hhw}_{ s-1}.
\eneqn

\item 
Consider the case $s\equiv_\ell k-1$.
In this case, we have $(s)^+_\jm=s+2$. Set $t=(s+2)_\jm^+$.
As in (i), let us show
\eq
\mathsf{M}^{\Dd,\hhw'}_t\simeq\mathsf{M}^{\Dd,\hhw}_{(s+1)_{\im}^+}.
\label{eq:k2}
\eneq
\bna
\item Case $t\eqv k-1$. Then $\jm_s=\im_{s+1}=\jm_t=\im_{t+1}$, and hence
$(s+1)_\im^+=t+1$.
Hence we have \eqref{eq:k2}.
\item Case $t\eqv k$. Then $\jm_s=\im_{s+1}=\jm_t=\im_{t-1}$, and hence
$(s+1)_\im^+=t-1$. Therefore we have \eqref{eq:k2}.
\item The case $t\eqv k+1$ does not occur.
\item Case $t\not\eqv k,k\pm1$. Then $t=(s+1)_\im^+$ and
we have \eqref{eq:k2}.
\ee 
Thus we complete the proof of \eqref{eq:k2}.

Hence, we have
\eqn
  \mathsf{M}^{\Dd,\hhw'}_{s}&&\simeq  \mathsf{M}^{\Dd,\hhw'}_{(s+2)_\jm^+}
\hconv \mathsf{S}_{s+2}^{\Dd,\hhw'}\hconv \mathsf{S}_{s}^{\Dd,\hhw'}  \simeq  \mathsf{M}^{\Dd,\hhw}_{(s+1)_\im^+} \hconv (  \mathsf{S}_{s}^{\Dd,\hhw} \hconv \mathsf{S}_{s+2}^{\Dd,\hhw}  )  \\
&&  \simeq  \mathsf{M}^{\Dd,\hhw}_{(s+1)_\im^+} \hconv (  \mathsf{S}_{s+1}^{\Dd,\hhw}  ) \simeq \mathsf{M}^{\Dd,\hhw}_{s+1},
\eneqn
which implies the assertion for this case.

\item 
Let us consider $s\equiv_\ell k+1$.
Similarly to the proof of \eqref{eq:k}, we can prove
\eqn
\mathsf{M}^{\Dd,\hhw'}_{(s)_\jm^+}\simeq\mathsf{M}^{\Dd,\hhw}_{(s-1)_{\im}^+}
\eneqn

Hence, we have
\eqn
&&  \mathsf{M}^{\Dd,\hhw'}_{s} \simeq \mathsf{M}^{\Dd,\hhw'}_{(s)^+_\jm}\hconv\mathsf{S}_{s}^{\Dd,\hhw'}   \simeq  \mathsf{M}^{\Dd,\hhw}_{(s-1)_\im^+} \hconv  \mathsf{S}_{s-2}^{\Dd,\hhw} \simeq (\mathsf{M}^{\Dd,\hhw}_{s})',
\eneqn
which implies the assertion in this case.

\item 
Finally, assume that $s\not\equiv_\ell k,k\pm1$.
Similarly to the proof of \eqref{eq:k}, we can prove
\eqn
\mathsf{M}^{\Dd,\hhw'}_{(s)_\jm^+}\simeq\mathsf{M}^{\Dd,\hhw}_{(s)_{\im}^+}
\eneqn
Hence, we have
\eqn
&&  \mathsf{M}^{\Dd,\hhw'}_{s} \simeq  \mathsf{M}^{\Dd,\hhw'}_{(s)^+_\jm}\hconv\mathsf{S}_{s}^{\Dd,\hhw'}   \simeq  \mathsf{M}^{\Dd,\hhw}_{(s)_\im^+} \hconv  \mathsf{S}_{s}^{\Dd,\hhw} \simeq \mathsf{M}^{\Dd,\hhw}_{s}. \qquad \qquad \qedhere
\eneqn
\ee
\QED

\begin{lemma}\label{lem: mutataion quiver2}
Let us keep the notations in the previous lemma. Then there exists a quiver isomorphism
$$ \Upphi: \mu_{k+1-\ell\Z_{>0}} \bl \GLS^{[-\infty,0]}(\hhw) \br \to \GLS^{[-\infty,0]}(\hhw'),$$
where
$$\mu_{k+1-\ell\Z_{>0}}(\GLS^{[-\infty,0]}(\hhw))\seteq     \cdots  \circ \mu_{k+1-2\ell} \circ \mu_{k+1-\ell} (\GLS^{[-\infty,0]}(\hhw)).$$
\end{lemma}

\begin{proof} Note that $\mu_{k+1-t\ell}  \circ \mu_{k+1-s\ell} \bl \GLS^{[-\infty,0]}(\hhw) \br= \mu_{k+1-s\ell}  \circ \mu_{k+1-t\ell} \bl \GLS^{[-\infty,0]}(\hhw) \br$
for any $t,s \in \Z_{>0}$. For $s \in \Z_{\le0}$, define the index change map $\Upphi:[-\infty,0] \to [-\infty,0] $ as follows:
$$\Upphi(s)
= \bc
s & \text{ if } s \not\equiv_\ell k, k-1 \\
s-1 & \text{ if } s \equiv_\ell k,\\
s+1 & \text{ if } s \equiv_\ell k-1.
\ec
$$

For each $t \in \Z_{\ge0}$ and $\im \in \If$, set  $ u^t_\im \seteq \min\st{ x \in [-(t+1)\ell+1,-t\ell]  \mid \im_x =\im } $.
Then, there is no arrow between
$$\text{$[-(t_1+1)\ell+1,-t_1\ell] \setminus \{ u^{t_1}_\im \ | \ \im \in \If \}$ and
$[-(t_2+1)\ell+1,-t_2\ell]$} $$ in $\GLS^{[-\infty,0]}(\hhw)$ when $t_1 > t_2 \in \Z_{\ge0}$.
Thus it is enough to show that the restriction of $\Phi$ to $[a,0]$ is a quiver isomorphism if $a$ is sufficiently small. Then, our assertion follows from
the lemma above together with
Proposition \ref{prop:uniqueexchange} and Lemma \ref{lem:smallfrozen}.
\end{proof}

\begin{proof}[Proof of Proposition~\ref{prop: mutation equivalent}]
From the above four lemmas, we can conclude the followings:
\bnum
\item If $\tw_0'$ can be obtained from $\tw_0$ by a single commutation move $\seed_-^{[-\infty,0],\Dd,\hhw}$ and  $\seed_-^{[-\infty,0],\Dd,\hhw'}$ are
monoidally mutation equivalent.
\item If $\tw_0'$ can be obtained from $\tw_0$ by a single braid move, the  sequence  $\mu_{k+1-\ell\Z_{>0}}$  gives a monoidally  mutation equivalence between
$\seed_-^{[-\infty,0],\Dd,\hhw}$ and  $\seed_-^{[-\infty,0],\Dd,\hhw'}$ (see Definition~\ref{def:mutationequiv}).
\ee
Thus our assertion follows from Remark~\ref{rmk: braid move cuspidal}~\eqref{it: moves}.
\end{proof}

As a corollary of Proposition~\ref{prop: mutation equivalent}, we obtain the following assertion.
\Cor\label{cor: mutation equivalent}
Let $(\Dd,\hhw)$ be a \pbw and $\tw_0'$ be another reduced expression of $w_0$. If the \Lad seed
$\seed_-^{[-\infty,0],\Dd_\Qd,\hhw}$ is a \ca monoidal seed in $\catg^{<\xi}$,
then so is $\seed_-^{[-\infty,0],\Dd_\Qd,\hhw'}$.
\encor

Thus we obtain \eqref{item:II}${}_\xi$ for $\xi$ satisfying
condition \eqref{cond:standard}.

\subsection{Proof of the main theorem}
By the result of the preceding subsection, we have already proved
that the monoidal seed $\seed_-^{[-\infty,0],\Dd_\Qd,\hhw}$
is a \ca monoidal seed when  $\Qd$ is a $\rmQ$-datum with the height function $\xi^t$ and $\tw_0$ is an arbitrary reduced expression of $w_0$.
 The next step
is to generalize this result to an arbitrary $\rmQ$-datum (with an arbitrary height function $\xi$) and an arbitrary interval $[a,b]$.
\Prop\label{prop: main}
Let $t\in\Z$ and assume that $\seed^{[-\infty,t],\,\Dd,\,  \hhw }_-$ is
a \ca monoidal seed in $\catCO$. Then, for any admissible chain
$\frakC=(\ci_k)_{k\in[1,l]}$ of $i$-boxes with a range
$[a,b]$, $\seed(\frakC)=\bl  \cM[](\frakC), \tB(\frakC)\KKC\br$
is a \ca monoidal seed in $\catg^{[a,b],\,\,\Dd,\,  \hhw }$ with a suitable
exchange matrix $\tB(\frakC)$.
\enprop

\Proof We shall prove first the case $-\infty<a\le b=t$. Set $\K_0=[a,t]$ and
$\Kex_0=\st{s\in\K_0\mid a\le s^-}$. Then there is no arrow in
$\GLS^{[-\infty,t],\Dd,  \hhw }$ between $\Kex_0$ and
$[-\infty,t]\setminus\K_0$. Hence 
$\seed^{[a,t],\, \Dd, \, \hhw}_-=\Bigl(\seed^{[-\infty,t], \,\Dd, \hhw}_-\Bigr)\big\vert_{(\K_0,\,\Kex_0)}$
is a \ca monoidal seed in $\catCO$  by Lemma~\ref{lem:restrseed}.
Hence $\seed^{[a,t],\,\Dd,  \hhw }_-$ is a \ca monoidal seed in
$\cat^{[a,t], \, \Dd,  \hhw }_\g$ by Proposition~\ref{prop:seedab}. Then,
Corollary~\ref{cor:2chains} implies that, for any admissible chain
$\frakC$ of $i$-boxes with the range $[a,t]$,
$\seed(\frakC)=\bl \cM[](\frakC),\tB(\frakC)\KKC\br$ is also a
\ca monoidal seed in $\cat^{[a,t],  \, \Dd,  \hhw}_\g$ for some exchange matrix
$\tB(\frakC)$.

\medskip
Now we treat the case $-\infty<a\le b\le t$. First we treat the
special case $\frakC_0=\frakC^{[a,b]}_- \seteq \seq{[b-k+1,b\}}_{1\le
k\le b-a+1}$. Then the monoidal seed
$$\seed(\frakC_0)=
\seed^{[a,b],  \, \Dd,  \hhw}_-\seteq\bl\st{M^{ \Dd,  \hhw }[s,b\}}_{s\in[a,b]},
\BGLS^{[a,b],  \, \Dd,  \hhw};\K(\frakC_0),\Kex(\frakC_0)\br$$ is
$\Uplambda$-admissible by 
Theorem~\ref{thm:triangle}.
Then we can enlarge $\frakC_0$ to an admissible chain $\frakC'
=\seq{\ci_k}_{1\le k\le l'}$ of $i$-boxes with the range $[a,t]$. For
example, we can take $l'=l+(t-b)$, $\ci_k=\{a,b+k-l]$ for $l<k\le
l'$. Since $\seed(\frakC_0)$ is $\Uplambda$-admissible,
Lemma~\ref{prop:uniqueexchange} implies that condition
\eqref{eq:inout} is satisfied for $\K(\frakC')$, $\K(\frakC_0)$ and
$\Kex(\frakC_0)$. Moreover it implies
$\seed(\frakC_0)=\seed(\frakC')\vert_{(\K(\frakC_0),\;\Kex(\frakC_0))}$.
We know already that $\seed(\frakC')$ is a \ca monoidal seed with a
suitable exchange matrix. Hence we can apply
Lemma~\ref{lem:restrseed} to conclude that $\seed(\frakC_0)$ is also
a \ca monoidal seed. Finally, Proposition~\ref{prop:seedab} implies
that $\seed(\frakC)$ is a \ca monoidal seed in $\catg^{[a,b],  \, \Dd,  \hhw}$
for any admissible chain $\frakC$ with a finite range $[a,b]$
provided that  $-\infty<a\le b\le t$.

\medskip
Now let $[a,b]$ be an arbitrary {\em finite} interval, and let $\frakC$ be an
admissible chain of $i$-boxes with the range $[a,b]$. There is
$m\in\Z_{>0}$ such that $b-2m\ell< t$. Here $\ell$ is the length of
the longest element of the Weyl group of $\gf$. Define the
automorphism $\delta$ of $\Z$ by $\delta(s)=s+\ell$. Then we have
$\cuspS{\delta(s)}\simeq\D \cuspS{s}$ by Definition \eqref{def:pbw}.
We set $\delta^{-2m}\frakC=\st{\delta^{-2m}\ci_k}\vert_{1\le k\le l}$. Then
the monoidal autofunctor $\D^{-2m}$ transforms $\seed(\frakC)$ to
$\seed(\delta^{-2m}\frakC)$ with the range $[a-2m\ell,b-2m\ell]$.
Since $\seed(\delta^{-2m}\frakC)$ is a \ca monoidal seed in
$\catg^{[a-2m\ell,b-2m\ell], \Dd,  \hhw}$ with a suitable exchange matrix,
we conclude that $\seed(\frakC)$ is a \ca monoidal seed in
$\catg^{[a,b], \Dd,  \hhw }$.

Now we take an arbitrary interval $[a,b]$,  possibly infinite. Let
$\frakC=\seq{\ci_k}_{1\le k\le l}$ with envelopes
$\tc_k=[\ta_k,\tb_k]$. For each $m$ such that $m\le l$, set $\frakC_{\le
m}=\seq{\ci_k}_{1\le k\le m}$ with the range $[\ta_m  ,  \tb_m]$. Then
$\seed(\frakC_{\le m})$ is a \ca monoidal seed in
$\catg^{[a,b],\Dd,  \hhw}$ for any $m$, and hence $\seed(\frakC)$ is also a
\ca monoidal seed in $\catg^{[a,b], \Dd,  \hhw}$. \QED

\Prop\label{prop:sixi} Let $\xi$ be a height function, and let $\im$
be a sink of $\xi$. If $\seed^{<s_\im\xi}$ is a \ca monoidal seed in
$\catCO$, then $\seed^{<\xi}$ is a \ca monoidal seed in $\catCO$
\enprop

\Proof Take a $\xi$-adapted admissible sequence $\hw$ in $\hI_\g$
such that $\im_1=\im$. Then we have
$\seed^{<s_\im\xi}=\seed_-^{[-\infty,1],\hw}$ and
$\seed^{<\xi}=\seed_-^{[-\infty,0],\hw}$. Then the result follows
from Proposition~\ref{prop: main}. \QED

\medskip
\Proof[Proof of Theorem~\ref{th:Main}] Let $\xi$ be the height
function associated with
 the admissible sequence $\hw$.
Take $t\in\Z$ such that $\xi_\im<t$ for any $\im\in\If$. Then, there
is a $\xi$-adapted reduced expression $\tw=s_{\im_1}\cdots
s_{\im_r}$ such that $\xi'\seteq s_{\im_r}\cdots s_{\im_1}\xi$
satisfies $\xi'_\im\in[t,t+2\fdi-1]$ for all $\im\in\If$. Hence
$\hI^{<\xi'}=\hI^{<t}$. Then Theorem~\ref{thm: main cg-} implies
that $\seed^{<\xi'}=\seed^{<t}$ is a \ca monoidal seed in $\catCO$.
Hence by a successive application of Proposition~\ref{prop:sixi},
the monoidal seed $\seed^{<\xi}$ is  \ca in $\catCO$.

By Corollary~\ref{cor: mutation equivalent}, the result for arbitrary $\seed^{<\xi}$ in Proposition~\ref{prop:sixi} implies that $\seed_-^{[-\infty,0],\Dd_\Qd,\hhw}$
is also a \ca monoidal seed for any $\rmQ$-datum $\Qd$ and any reduced expression $\tw_0$ of $w_0$.

Hence we conclude that for any admissible chain  $\frakC$ of
$i$-boxes with the range $[a,b]$, $\seed\seteq\seed(\frakC)$ is a
\ca monoidal seed in $\catg^{[a,b],\, \Dd_\Qd, \hhw}$ by Proposition~\ref{prop:
main}.

\medskip
It remains to prove that $\A([\seed])\simeq K(\catg^{[a,b],\, \Dd_\Qd, \hhw})$.
By truncating $\frakC$, we may assume that $[a,b]$ is a finite
interval. Set $\seed=\bl\st{M_k}_{1\le k\le l},\tB,\KK)$.

Let $\st{X_k}_{1\le k\le l}$ be the cluster variables, and set
$f\seteq\prod_{k=1}^l[M_k]\in K(\catg^{[a,b],\, \Dd_\Qd, \hhw })$. Then we have
$$\A([\seed])\subset \Z[X_1^{\pm1},\ldots, X_l^{\pm1}]\to
K(\catg^{[a,b],\, \Dd_\Qd, \hhw })[f^{-1}]$$ by $X_k\mapsto[M_k]$. Since $\seed$
is \ca in $\catg^{[a,b],\, \Dd_\Qd, \hhw}$, it induces a ring homomorphism
$$\A([\seed])\To K(\catg^{[a,b],\, \Dd_\Qd, \hhw }).$$
For any $s\in[a,b]$, after successive box moves, the moved $\frakC$
contains $\st{[s]}$. Hence, the image of  $\A([\seed])$ contains
$[\cuspS{s}]$. Since $K(\catg^{[a,b],\,\Dd_\Qd, \hhw })$ is the polynomial ring
with the system of generators $\st{[\cuspS{s}]}_{s\in[a,b]}$,
$\A([\seed])\to K(\catg^{[a,b],\, \Dd_\Qd, \hhw})$ is surjective. Since their
dimensions are equal to $| [a,b]|$, it is an isomorphism. \QED

\subsection{Conjecture}

We give the following conjecture which asserts that Theorem~\ref{th:Main}
holds for an arbitrary \pbw $(\Dd,\hhw)$.

\begin{conjecture} \label{conj: Lambda-ad} \hfill
\bnum
\item
For any \pbw $(\Dd,\hhw)$ and   any admissible chain $\frakC=\seq{\ci_k}_{1\le k\le l}$ of $i$-boxes with an arbitrary  range $[a,b]$, the monoidal category
$\cat^{[a,b],\Dd,\hhw}_\g$ is a $\Uplambda$-monoidal
categorification of $K(\cat^{[a,b],\Dd,\hhw}_\g)$ with an initial \Lad
 monoidal seed
 $(\cM[](\frakC),\tB\KKC)$
for some exchange matrix $\tB$.
\item If $-\infty\le a\le b<+\infty$,
the monoidal category $\cat^{[a,b],\Dd,\hhw}_\g$ is a $\Uplambda$-monoidal
categorification of $K(\cat^{[a,b],\Dd,\hhw}_\g)$ with an initial \Lad
monoidal seed
$$
\seed^{[a,b],\Dd,\hhw}_-\seteq\bl\st{M^{\Dd,\hhw}[s,b\}}_{s\in[a,b]},\BGLS^{[a,b],\Dd,\hhw}\KK\br.
$$
\ee
\end{conjecture}

\Prob For any \pbw and an admissible chain
$\frakC$ of $i$-boxes, determine explicitly the exchange matrix of
$\seed(\frakC)$. \enprob

\end{document}